\documentclass[smallextended]{svjour3}
\usepackage{amsmath}
\usepackage{stmaryrd}
\usepackage{bbding}
\usepackage{dcolumn}
\usepackage{graphicx}
\usepackage{amsfonts}
\usepackage{amssymb}
\usepackage{psfrag}
\usepackage{wrapfig}
\usepackage{subfigure}
\usepackage{makeidx}
\usepackage{bm}
\usepackage{epsf}
\usepackage{epsfig}
\usepackage{setspace}
\usepackage{color}
\usepackage{algorithm}
\usepackage{algorithmic}
\usepackage{setspace}

\definecolor{SpringGreen4}{RGB}{0,139,69}

\begin{document}
\title{Three-dimensional third-order gas-kinetic scheme on hybrid unstructured meshes for Euler and Navier-Stokes equations}

\author{Yaqing Yang \and Liang Pan  \and Kun Xu}
\institute{ Yaqing Yang
\at Laboratory of Mathematics and Complex Systems, School of Mathematical Sciences, Beijing Normal University, Beijing, China \\
\email{yqyangbnu@163.com}
\and
Liang Pan
\at Laboratory of Mathematics and Complex Systems, School of Mathematical Sciences, Beijing Normal University, Beijing, China \\
\email{panliang@bnu.edu.cn}
\and
Kun Xu \at
Department of Mathematics, Hong Kong University of Science and Technology, Clear Water Bay, Kowloon, Hong Kong\\
Shenzhen Research Institute, Hong Kong University of Science and Technology, Shenzhen, China\\
\email{makxu@ust.hk}}

\date{Received: date/Accepted: date}

\maketitle
\begin{abstract}
In this paper, a third-order gas-kinetic scheme is developed on the
three-dimensional hybrid unstructured meshes for the numerical
simulation of compressible inviscid and viscous flows. In the
classical weighted essentially non-oscillatory (WENO) scheme, the
high-order spatial accuracy is achieved by the non-linear
combination of lower order polynomials. However, for the hybrid
unstructured meshes, the procedures, including the selection of
candidate stencils and calculation of linear weights, become
extremely complicated, especially for three-dimensional problems. To
overcome the drawbacks, a third-order WENO reconstruction is
developed on the hybrid unstructured meshes, including tetrahedron,
pyramid, prism and hexahedron. A simple strategy is adopted for the
selection of big stencil and sub-stencils, and the topologically
independent linear weights are used in the spatial reconstruction. A
unified interpolation is used for the volume integral of different
control volumes, as well as the flux integration over different cell
interfaces. Incorporate with the two-stage fourth-order temporal
discretization, the explicit high-order gas-kinetic schemes are
developed for unsteady problems. With the lower-upper symmetric
Gauss-Seidel (LU-SGS) methods, the implicit high-order gas-kinetic
schemes are developed for steady problems. A variety of numerical
examples, from the subsonic to supersonic flows, are presented to
validate the accuracy and robustness for both inviscid and viscous
flows. To accelerate the computation, the
current scheme is implemented with the graphics processing unit (GPU) using
compute unified device architecture (CUDA).  The speedup of GPU
code suggests the potential of high-order gas-kinetic schemes for the
large scale computation.  

\end{abstract}

\maketitle
\section{Introduction}
The unstructured meshes are widely used in the computational fluid
dynamics (CFD) methods with complex geometries. Various high-order
numerical methods on unstructured meshes have been developed in the
past decades, including discontinuous Galerkin (DG)
\cite{DG-A,DG-B,DG-C}, spectral volume (SV) \cite{SV-A}, correction
procedure using reconstruction (CPR) \cite{CPR-A},  the
finite-volume type essential non-oscillatory (ENO) \cite{ENO-un},
weighted essential non-oscillatory (WENO)
\cite{WENO-un1,WENO-un2,WENO-un3,WENO-un4} and Hermite WENO (HWENO)
\cite{HWENO} methods.

In the framework of finite volume methods, the WENO schemes have
been successfully developed for the compressible flows, the
WENO-type high-order methods have received the most attention in
recent years. For the one-dimensional WENO schemes \cite{WENO-Liu,
WENO-JS, WENO-Z}, the high-order of accuracy is obtained by the
non-linear combination of lower order polynomials from the candidate
stencils. For the multi-dimensional structured meshes, the WENO
scheme is performed dimension-by-dimension. On the two-dimensional
unstructured meshes, the WENO schemes were also developed with same
idea \cite{WENO-un1,WENO-un2,WENO-un3}. However, the linear weights
need to be computed and restored for each cell, and the appearance
of negative weights also affect the performance of WENO schemes. The
central/compact WENO (CWENO) schemes were developed when facing
distorted local mesh geometries or degenerate cases
\cite{CWENO1,CWENO2,CWENO3}. Following the original idea of
classical CWENO schemes, two types of WENO scheme is developed,
i.e., the WENO schemes with adaptive order
\cite{WENO-ao-1,WENO-ao-2} and the simple WENO schemes
\cite{WENO-simple-1,WENO-simple-2}. The linear weights are topology
independent and artificially set to be positive numbers, and the
non-linear weights are chosen to achieve the optimal order of
accuracy in the smooth region and suppress the oscillations near the
discontinuous region. Most of efforts are spent in the unstructured
triangular meshes for two-dimensional problems and tetrahedral
meshes for three-dimensional problems. However, in practical
applications, such as the viscous flows around or inside complex
geometries, the hybrid unstructured meshes are usually adopted for
accuracy and efficiency. More recently, the WENO schemes on
arbitrary unstructured meshes has been developed for inviscid flows
\cite{WENO-un4,WENO-un5,WENO-un6}, and also extended for viscous
flows, including laminar, transitional and turbulent problems
\cite{WENO-un7}.

In the past decades, the gas-kinetic schemes (GKS) based on the
Bhatnagar-Gross-Krook (BGK) model \cite{BGK-1,BGK-2} have been
developed systematically for the computations from low speed flows
to supersonic ones \cite{GKS-Xu1,GKS-Xu2}. The gas-kinetic scheme
presents a gas evolution process from the kinetic scale to
hydrodynamic scale, and  both inviscid and viscous fluxes can be
calculated in one framework. Recently, a time-dependent gas
distribution function can be constructed at a cell interface, which
is important for the construction of high-order scheme. With the
two-stage temporal discretization, which was originally developed
for the Lax-Wendroff type flow solvers \cite{GRP-high-1,GRP-high-2},
a reliable framework was provided to construct gas-kinetic scheme
with fourth-order and even higher-order temporal accuracy
\cite{GKS-high-1,GKS-high-2}. More importantly, the high-order
scheme is as robust as the second-order one and works perfectly from
the subsonic to hypersonic flows. The implicit methods for GKS and
unified GKS have also been constructed
\cite{GKS-implicit-1,GKS-implicit-2}, and the implicit temporal
methods provide efficient techniques for speeding up the convergence
of steady problems. Recently, with the simple WENO type
reconstruction, the third-order and fourth-order gas-kinetic schemes
are developed on the three-dimensional structured meshes, in which a
simple strategy of selecting stencils for reconstruction is adopted
and the topology independent linear weights are used
\cite{GKS-high-4,GKS-high-5}. Based on the spatial and temporal
coupled property of GKS solver and the HWENO reconstruction
\cite{HWENO}, the high-order compact gas-kinetic schemes are also
developed \cite{GKS-high-6,GKS-high-7}.

In this paper, a third-order gas-kinetic scheme is developed on the
three-dimensional hybrid unstructured meshes for the compressible
inviscid and viscous flows. Due to the complex mesh topology of
hybrid unstructured meshes, the procedures of classical WENO scheme,
including the selection of candidate stencils and calculation of
linear weights, become extremely complicated. A simple WENO
reconstruction is extended to the unstructured meshes with
tetrahedral, pyramidal, prismatic and hexahedral cells. A large
stencil is selected with the neighboring cells and the neighboring
cells of neighboring cells, and a quadratic polynomial can be
obtained. For the tetrahedral and pyramidal cells, the centroids of
control volume and neighboring cells might be coplanar, and the
coefficient matrix might become singular. A robust selections of
candidate sub-stencils are also given, such that the linear
polynomials are solvable for each sub-stencil. The trilinear
interpolation is used for each cell, and the Gaussian quadrature
over a plane cell interface is used to achieve the spatial accuracy.
For the cell interface, the trilinear interpolation degenerates to a
bilinear interpolation, and a unified formulation can be used to
calculate the numerical fluxes over both triangular and
quadrilateral interfaces.  Incorporate with the two-stage
fourth-order temporal discretization, the explicit high-order
gas-kinetic schemes are developed for unsteady problems. With the
lower-upper symmetric Gauss-Seidel (LU-SGS) methods, the implicit
high-order gas-kinetic schemes are developed for steady problems.
Various three-dimensional numerical experiments, including unsteady
and steady problems, are presented to to validate the accuracy and
robustness of WENO scheme. To accelerate the computation, the current scheme
is implemented to run on graphics processing unit (GPU) using
compute unified device architecture (CUDA). The computational
efficiency using single Nvidia TITAN RTX GPU is demonstrated. Obtained results are compared with those obtained by
an octa-core Intel i7-11700 CPU in terms of calculation time.
Compared with the CPU code, 6x speedup is achieved for GPU code.  In the future, more challenging
compressible flow problems will be investigated with
multiple GPUs.

This paper is organized as follows. In Section 2, the
three-dimensional WENO reconstruction on the hybrid unstructured
meshes will be introduced. The high-order gas-kinetic scheme for
inviscid and viscous flows will be presented in Section 3. Numerical
examples are included in Section 4. The last section is the
conclusion.

\begin{figure}[!h]
\centering
\includegraphics[width=0.65\textwidth]{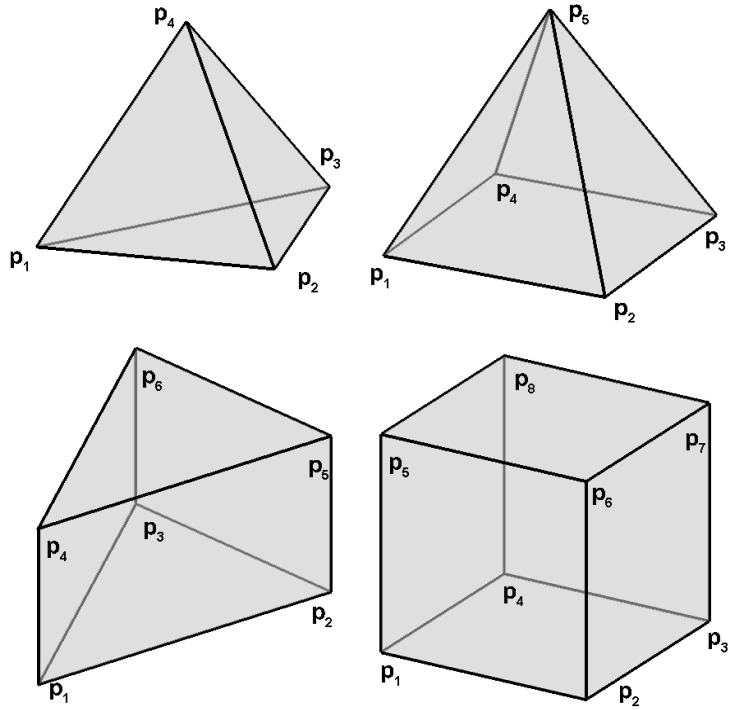}
\caption{\label{schematic-cell} Schematic for tetrahedron, pyramid,
prism and hexahedron.}
\end{figure}

\section{WENO reconstruction on hybrid unstructured meshes}
In this section, a third-order WENO reconstruction will be presented
on the three-dimensional unstructured hybrid meshes. Similar with
the previous study, the idea of simple WENO reconstruction is
adopted \cite{WENO-simple-1,WENO-simple-2,GKS-high-5,GKS-high-4}.
For the unstructured hybrid meshes, the meshes are consist of the
tetrahedral, pyramidal, prismatic and hexahedral cells. The
schematic of the cells are given in Fig.\ref{schematic-cell} with
the label of vertexes. For the sake of clearness, the faces of
control volume also need to be labeled. For the cell $\Omega_0$, the
faces are labeled as follows
\begin{itemize}
\item For the tetrahedral cell, four faces are denoted as
\begin{align*}
F_{1}=\{p_1p_2p_3\}, F_{2}=\{p_1p_2p_4\},
F_{3}=\{p_2p_3p_4\}, F_{4}=\{p_3p_1p_4\}.
\end{align*}
\item For the pyramidal cell, five faces are denoted as
\begin{align*}
F_{1}=\{p_1p_2p_5\}, F_{2}=\{p_2p_3p_5\},F_{3}=\{p_3p_4p_5\}, F_{4}=\{p_4p_1p_5\}, F_{5}=\{p_1p_2p_3p_4\}.
\end{align*}
\item For the prismatic cell, five faces are denoted as
\begin{align*}
F_{1}=\{p_1p_2p_3\}, F_{2}=\{p_4p_5p_6\},
F_{3}=\{p_1p_2p_6p_4\}, F_{4}=\{p_2p_3p_6p_5\}, F_{5}=\{p_3p_1p_4p_6\}.
\end{align*}
\item For the hexahedral cell, six faces are denoted as
\begin{align*}
F_{1}=\{p_1p_2p_3p_4\}, F_{2}=\{p_1p_2p_6p_5\}, F_{3}=\{p_2p_3p_7p_6\},\\
F_{4}=\{p_3p_4p_8p_7\}, F_{5}=\{p_4p_1p_5p_8\}, F_{6}=\{p_5p_6p_7p_8\}.
\end{align*}
\end{itemize}
With the labeled faces given above, the neighboring cell of
$\Omega_{i}$, which shares the face $F_{p}$, is denoted as
$\Omega_{i_p}$. Meanwhile, the neighboring cells of $\Omega_{i_p}$
are denoted as $\Omega_{i_{pm}}$. To achieve the third-order
accuracy, a big stencil $S_i$ for cell $\Omega_i$  is selected as
follows
\begin{align*}
S_i&=\{\Omega_{i},\Omega_{i_p},\Omega_{i_{pm}}\},
\end{align*}
which is consist of the neighboring cells and neighboring cells of
neighboring cells of $\Omega_{i}$.

Taking the boundary condition into account, the selection of big
stencil becomes more complicated. In the following, the inner cell
represents the cell with no face on the boundary and the boundary
cell is the cell with at least one face on the boundary. For the
periodic boundary condition,  the generated triangular or
quadrilateral meshes on the periodic boundaries should be identical.
If $\Omega_{i}$ is boundary cell,  the neighboring cells of
$\Omega_{i}$ can be found with the periodic boundary condition.
Thus,  for each cell $\Omega_{i}$ in the domain, the big stencil
$S_i$ can be selected. For other boundary conditions, the following
two stages need to be considered. If $\Omega_{i}$ is an inner cell
and $\Omega_{i_p}$ is a boundary cell,  the geometric information of
$\Omega_{i_{pm}}$ can be provided according to the mirror image for
cell interface and the flow variables are given according to the
boundary condition. If $\Omega_{i}$ is a boundary cell, the
geometric information and flow variables of the ghost cell
$\Omega_{i_{p}}$ are provided firstly. For the ghost cell
$\Omega_{i_{p}}$, the geometric information  of its neighboring
cells $\Omega_{i_{pm}}$ are provided by the mirror image of
neighboring cells of $\Omega_{i}$, and the flow variables on the
cell $\Omega_{i_{pm}}$ are given according to the boundary condition
corresponding to neighboring cells of $\Omega_{i}$.

With the procedure above, the candidate stencils can be selected and
the index of big stencil is rearranged as
\begin{align*}
S_i=\{\Omega_0,\Omega_1,...,\Omega_{N},\Omega_{N+1},...,\Omega_{K}\},
\end{align*}
where $K$ is the total number of cells and $N$ is the number of
neighboring cells. Based on the big stencil, a quadratic polynomial
can be constructed as
\begin{equation*}
P_0(\boldsymbol{x})=Q_{0}+\sum_{|\boldsymbol n|=1}^2a_{\boldsymbol n}p_{\boldsymbol n}(\boldsymbol{x}),
\end{equation*}
where $Q_{0}$ is the cell averaged conservative variables over
$\Omega_{0}$, the multi-index $\boldsymbol n=(n_1, n_2, n_3)$ and
$|\boldsymbol n|=n_1+n_2+n_3$. The base function $p_{\boldsymbol
n}(\boldsymbol{x})$ is defined as
\begin{align*}
\displaystyle
p_{\boldsymbol n}(\boldsymbol{x})=x^{n_1}y^{n_2}z^{n_3}-\frac{1}{\left|\Omega_{0}\right|}\iiint_{\Omega_{0}}x^{n_1}y^{n_2}z^{n_3}\text{d}V.
\end{align*}
To determine this polynomial, the following constrains need to be
satisfied for all cells in the big stencil
\begin{align*}
\frac{1}{\left|\Omega_{k}\right|}\iiint_{\Omega_{k}}P_0(\boldsymbol{x})\text{d}V=Q_{k},~\Omega_{k}\in S,
\end{align*}
where $Q_{k}$ is the conservative variable with newly rearranged
index. An over-determined linear system can be generated and the
least square method is used to obtain the coefficients.

To deal with the discontinuity, the linear polynomials are
constructed based on the candidate sub-stencils
\begin{equation}\label{lin-def}
P_m(\boldsymbol{x})=Q_{0}+\sum_{|\boldsymbol n|=1}b_{\boldsymbol n}^mp_{\boldsymbol n}(\boldsymbol{x}),
\end{equation}
where $m=1,...M$ and $M$ is the number of sub-stencils.  To
determine these polynomials, the following constrains need to be
satisfied
\begin{align*}
\frac{1}{\left|\Omega_{m_k}\right|}\iiint_{\Omega_{k}}P_m(\boldsymbol{x})\text{d}V=Q_{k},~\Omega_{m_k}\in S_{i_m},
\end{align*}
where $Q_{m_k}$ is the conservative variable with newly rearranged
index. With the selected sub-stencils $ S_{i_m}$, the the linear
polynomials can be fully determined.
\begin{itemize}
\item For the hexahedral cell, $M=8$ and the sub-candidate stencils are selected as follows
\begin{align*}
S_{i_1}=\{\Omega_{i},\Omega_{i_1},\Omega_{i_2},\Omega_{i_3}\},~S_{i_5}=\{\Omega_{i},\Omega_{i_6},\Omega_{i_2},\Omega_{i_3}\},\\
S_{i_2}=\{\Omega_{i},\Omega_{i_1},\Omega_{i_3},\Omega_{i_4}\},~S_{i_6}=\{\Omega_{i},\Omega_{i_6},\Omega_{i_3},\Omega_{i_4}\},\\
S_{i_3}=\{\Omega_{i},\Omega_{i_1},\Omega_{i_4},\Omega_{i_5}\},~S_{i_7}=\{\Omega_{i},\Omega_{i_6},\Omega_{i_4},\Omega_{i_5}\},\\
S_{i_4}=\{\Omega_{i},\Omega_{i_1},\Omega_{i_5},\Omega_{i_2}\},~S_{i_8}=\{\Omega_{i},\Omega_{i_6},\Omega_{i_5},\Omega_{i_2}\}.
\end{align*}
\item For the prismatic cell, $M=6$ and the sub-candidate stencils are selected as
\begin{align*}
S_{i_1}=\{\Omega_{i},\Omega_{i_1},\Omega_{i_3},\Omega_{i_4}\},~S_{i_4}=\{\Omega_{i},\Omega_{i_2},\Omega_{i_3},\Omega_{i_4}\},\\
S_{i_2}=\{\Omega_{i},\Omega_{i_1},\Omega_{i_4},\Omega_{i_5}\},~S_{i_5}=\{\Omega_{i},\Omega_{i_2},\Omega_{i_4},\Omega_{i_5}\},\\
S_{i_3}=\{\Omega_{i},\Omega_{i_1},\Omega_{i_5},\Omega_{i_3}\},~S_{i_6}=\{\Omega_{i},\Omega_{i_2},\Omega_{i_5},\Omega_{i_3}\}.
\end{align*}
\end{itemize}
For the hexahedral and pyramidal cells, the sub-stencils are consist
of $\Omega_{0}$ and three neighboring cells. The linear polynomials
can be determined.

However, for the tetrahedral and pyramidal cells, the centroids of
$\Omega_{0}$ and three of neighboring cells might becomes coplanar.
If we only use the neighboring cells for the sub-candidate stencils,
the coefficient matrix might become singular and more cells are
needed.
\begin{itemize}
\item For the tetrahedral cell, $M=4$ and the sub-candidate stencils are selected as
\begin{align*}
S_{i_1}=\{\Omega_{i},\Omega_{i_1},\Omega_{i_2},\Omega_{i_3},\Omega_{i_{11}},\Omega_{i_{12}},\Omega_{i_{13}}\},\\
S_{i_2}=\{\Omega_{i},\Omega_{i_1},\Omega_{i_2},\Omega_{i_4},\Omega_{i_{21}},\Omega_{i_{22}},\Omega_{i_{23}}\},\\
S_{i_3}=\{\Omega_{i},\Omega_{i_2},\Omega_{i_3},\Omega_{i_4},\Omega_{i_{31}},\Omega_{i_{32}},\Omega_{i_{33}}\},\\
S_{i_4}=\{\Omega_{i},\Omega_{i_3},\Omega_{i_1},\Omega_{i_4},\Omega_{i_{41}},\Omega_{i_{42}},\Omega_{i_{43}}\}.
\end{align*}
\item For the pyramidal cell, $M=8$ and the sub-candidate stencils are selected as
\begin{align*}
S_{i_1}&=\{\Omega_{i},\Omega_{i_1},\Omega_{i_2},\Omega_{i_3},\Omega_{i_{11}},\Omega_{i_{12}},\Omega_{i_{13}}\},~
S_{i_5}=\{\Omega_{i},\Omega_{i_1},\Omega_{i_2},\Omega_{i_5},\Omega_{i_{11}},\Omega_{i_{12}},\Omega_{i_{13}}\},\\
S_{i_2}&=\{\Omega_{i},\Omega_{i_2},\Omega_{i_3},\Omega_{i_4},\Omega_{i_{21}},\Omega_{i_{22}},\Omega_{i_{23}}\},~
S_{i_6}=\{\Omega_{i},\Omega_{i_2},\Omega_{i_3},\Omega_{i_5},\Omega_{i_{21}},\Omega_{i_{22}},\Omega_{i_{23}}\},\\
S_{i_3}&=\{\Omega_{i},\Omega_{i_3},\Omega_{i_4},\Omega_{i_1},\Omega_{i_{31}},\Omega_{i_{32}},\Omega_{i_{33}}\},~
S_{i_7}=\{\Omega_{i},\Omega_{i_3},\Omega_{i_4},\Omega_{i_5},\Omega_{i_{31}},\Omega_{i_{32}},\Omega_{i_{33}}\},\\
S_{i_4}&=\{\Omega_{i},\Omega_{i_4},\Omega_{i_1},\Omega_{i_2},\Omega_{i_{41}},\Omega_{i_{42}},\Omega_{i_{43}}\},~
S_{i_8}=\{\Omega_{i},\Omega_{i_4},\Omega_{i_1},\Omega_{i_5},\Omega_{i_{41}},\Omega_{i_{42}},\Omega_{i_{43}}\}.
\end{align*}
\end{itemize}
where $\Omega_{i_{pm}}\neq\Omega_{i}, p=1,2,3,4, m=1,2,3$, and the
cells of sub-candidate stencils are consist of the three neighboring
cells and three neighboring cells of one neighboring cell. With such
an enlarged sub-stencils, Eq.\eqref{lin-def} becomes solvable and
the linear polynomials can be determined.

With the reconstructed polynomial $P_m(\boldsymbol{x}), m=0,...,M$,
the point-value $Q(\boldsymbol{x}_{G})$ and the spatial derivatives
$\partial_{x,y,z} Q(\boldsymbol{x}_{G})$  for reconstructed
variables  at Gaussian quadrature point can be given by the
non-linear combination
\begin{equation}\label{weno}
\begin{split}
Q(\boldsymbol{x}_{G})=\overline{\omega}_0(\frac{1}{\gamma_0}P_0(\boldsymbol{x}_{G})-&
\sum_{m=1}^{M}\frac{\gamma_m}{\gamma_0}P_m(\boldsymbol{x}_{G}))+\sum_{m=1}^{M}\overline{\omega}_mP_m(\boldsymbol{x}_{G}),\\
\partial_{x,y,z} Q(\boldsymbol{x}_{G})=\overline{\omega}_0(\frac{1}{\gamma_0}\partial_{x,y,z}
P_0(\boldsymbol{x}_{G})-&\sum_{m=1}^{M}\frac{\gamma_m}{\gamma_0}\partial_{x,y,z}
P_m(\boldsymbol{x}_{G}))+\sum_{m=1}^{M}\overline{\omega}_m\partial_{x,y,z}
P_m(\boldsymbol{x}_{G})),
\end{split}
\end{equation}
where $\gamma_0, \gamma_1,...,\gamma_M$ are the linear weights. In
the computation, $\gamma_i=0.0125, i=1,...,M$ and
$\gamma_0=1-0.0125M$. $\partial_{x,y,z} P_m(\boldsymbol{x}_{G})$ can
be obtained by taking derivatives of the candidate polynomials
directly. The non-linear weights $\omega_m$ and normalized
non-linear weights $\overline{\omega}_m$ are defined as
\begin{align*}
\overline{\omega}_{m}=\frac{\omega_{m}}{\sum_{m=0}^{M} \omega_{m}},~
\omega_{m}=\gamma_{m}\Big[1+\Big(\frac{\tau}{\beta_{m}+\epsilon}\Big)\Big],
\end{align*}
where $\epsilon$ is a small positive number.
To achieve a third-order accuracy, a quadratic polynomial is used
for $P_0(\boldsymbol{x})$ and $\tau$ is chosen as
\begin{align*}
\tau=\sum_{m=1}^{M}\Big(\frac{|\beta_0-\beta_m|}{M}\Big).
\end{align*}
The smooth indicator $\beta_{m}$ is defined as
\begin{align}\label{indicator}
\beta_m=\sum_{|l|=1}^{r_m}|\Omega_{ijk}|^{\frac{2|l|}{3}-1}\int_{\Omega_{ijk}}
\Big(\frac{\partial^lP_m}{\partial_x^{l_1}\partial_y^{l_2}\partial_z^{l_3}}(x,y,z)\Big)^2\text{d}V,
\end{align}
where $r_0=2$ and $r_m=1$ for $m=1,...,M$. It can be proved that
Eq.\eqref{weno} ensures third-order accuracy and more details can be
found in \cite{GKS-high-4}.

For the spatial reconstruction, the integrals over tetrahedron,
pyramid, prism and hexahedron, including volume of cell, need to be
calculated. The volume of tetrahedron can be given by an explicit
formulation. The pyramid, prism and hexahedron can be divided into
several tetrahedrons, and their integrals can be calculated simply.
However, the quadrilateral interface of pyramid, prism and
hexahedron maybe non-coplanar, such kind of method for calculating
volume is inaccurate. For the sake of accuracy, the integrals over a
hexahedron cell is computed by the Gaussian quadrature
\begin{align*}
&\iiint_{\Omega}x^{a}y^{b}z^{c}\text{d}x\text{d}y\text{d}z=
\sum_{l,m,n=1}^3\omega_{lmn}x^{a}y^{b}z^{c}(\xi_l,\eta_m,\zeta_n)\Big|
\frac{\partial(x,y,z)}{\partial(\xi,\eta,\zeta)}\Big|_{(\xi_l,\eta_m,\zeta_n)}\Delta\xi\Delta\eta\Delta\zeta,
\end{align*}
where $\omega_{lmn}$ is the Gaussian quadrature weight and
$(\xi_l,\eta_m,\zeta_n)$ is the quadrature point. A trilinear
interpolation is introduced to parameterize the hexahedron cell
$\Omega_{i}$
\begin{align}\label{trilinear}
\boldsymbol{X}(\xi,\eta,\zeta)=\sum_{m=1}^8\boldsymbol{x}_m\psi_m(\xi,\eta,\zeta),
\end{align}
where $(\xi,\eta,\zeta)\in[-1/2,1/2]^3$, $\boldsymbol{x}_m$ is the
vertex of a hexahedron cell and the base function $\psi_m$ is given
as follows
\begin{align*}
\psi_1=\frac{1}{8}(1-2\xi)(1-2\eta)(1-2\zeta),
\psi_2=\frac{1}{8}(1-2\xi)(1-2\eta)(1+2\zeta),\\
\psi_3=\frac{1}{8}(1-2\xi)(1+2\eta)(1-2\zeta),
\psi_4=\frac{1}{8}(1-2\xi)(1+2\eta)(1+2\zeta),\\
\psi_5=\frac{1}{8}(1+2\xi)(1-2\eta)(1-2\zeta),
\psi_6=\frac{1}{8}(1+2\xi)(1-2\eta)(1+2\zeta),\\
\psi_7=\frac{1}{8}(1+2\xi)(1+2\eta)(1-2\zeta),
\psi_8=\frac{1}{8}(1+2\xi)(1+2\eta)(1+2\zeta).
\end{align*}
With trilinear interpolation Eq.\eqref{trilinear}, the hexahedron
cell  with non-coplanar vertexes can be dealt with accurately. As
shown in Fig.\ref{schematic-cell}, the hexahedron cell denoted as
the sequential  label $\Omega=\{p_1p_2p_3p_4p_5p_6p_7p_8\}$. For the
tetrahedral, pyramidal and prismatic cells, they can be considered
as degenerated hexahedron, and can be relabeled follows
\begin{align*}
\Omega&=\{p_1p_2p_3p_3p_4p_4p_4p_4\} \text{~for tetrahedron},\\
\Omega&=\{p_1p_2p_3p_4p_5p_5p_5p_5\} \text{~for pyramid},\\
\Omega&=\{p_1p_2p_3p_3p_4p_5p_6p_6\} \text{~for prism}.
\end{align*}
Eq.\eqref{trilinear} is also used to parameter the tetrahedral,
pyramidal, prismatic cells and the integrals can be calculated.

\section{High-order gas-kinetic scheme}
\subsection{BGK equation and finite volume scheme}
With the thrid-order WENO reconstruction, the high-order gas-kinetic
scheme (HGKS) will be presented for the three-dimensional flows in
the finite volume framework. The three-dimensional BGK equation
\cite{BGK-1,BGK-2} can be written as
\begin{equation}\label{bgk-3d}
f_t+uf_x+vf_y+wf_z=\frac{g-f}{\tau},
\end{equation}
where $\boldsymbol{u}=(u,v,w)$ is the particle velocity, $f$ is the
gas distribution function, $g$ is the three-dimensional Maxwellian
distribution and $\tau$ is the collision time. The collision term
satisfies the compatibility condition
\begin{equation}\label{compatibility}
\int \frac{g-f}{\tau}\psi \text{d}\Xi=0,
\end{equation}
where
$\displaystyle\psi=(\psi_1,...,\psi_5)^T=(1,u,v,w,\frac{1}{2}(u^2+v^2+w^2+\varsigma^2))^T$,
the internal variables
$\varsigma^2=\varsigma_1^2+...+\varsigma_K^2$,
$\text{d}\Xi=\text{d}u\text{d}v\text{d}w\text{d}\varsigma^1...\text{d}\varsigma^{K}$,
$\gamma$ is the specific heat ratio and  $K=(5-3\gamma)/(\gamma-1)$
is the degrees of freedom for three-dimensional flows. According to
the Chapman-Enskog expansion for BGK equation, the macroscopic
governing equations can be derived \cite{GKS-Xu1,GKS-Xu2}. In the
continuum region, the BGK equation can be rearranged and the gas
distribution function can be expanded as
\begin{align*}
f=g-\tau D_{\boldsymbol{u}}g+\tau D_{\boldsymbol{u}}(\tau
D_{\boldsymbol{u}})g-\tau D_{\boldsymbol{u}}[\tau
D_{\boldsymbol{u}}(\tau D_{\boldsymbol{u}})g]+...,
\end{align*}
where $D_{\boldsymbol{u}}=\displaystyle\frac{\partial}{\partial
t}+\boldsymbol{u}\cdot \nabla$. With the zeroth-order truncation
$f=g$, the Euler equations can be obtained. For the first-order
truncation
\begin{align*}
f=g-\tau (ug_x+vg_y+wg_z+g_t),
\end{align*}
the Navier-Stokes equations can be obtained \cite{GKS-Xu1,GKS-Xu2}.

Taking moments of Eq.\eqref{bgk-3d} and integrating with respect to
space, the semi-discretized finite volume scheme can be expressed as
\begin{align}\label{semi}
\frac{\partial Q_i}{\partial t}=\mathcal{L}(Q_{i}),
\end{align}
where $Q=(\rho, \rho U,\rho V, \rho W, \rho E)$ is the cell averaged
conservative value of $\Omega_{i}$. The operator $\mathcal{L}$ is
defined as
\begin{equation}\label{finite}
\mathcal{L}(Q_{i})=-\frac{1}{|\Omega_{i}|}\sum_{i_p\in N(i)}\mathbb{F}_{i_p}(t)=-\frac{1}{|\Omega_{i}|}\sum_{i_p\in N(i)}\iint_{\Sigma_{i_p}}\boldsymbol{F}(Q,t)\text{d}\sigma,
\end{equation}
where $Q_{i}$ is the cell averaged conservative variables of
$\Omega_{i}$, $|\Omega_{i}|$ is the area of $\Omega_{i}$,
$\Sigma_{i_p}$ is the common cell interface of $\Omega_{i}$ and
$\Omega_{i_p}$ and $N(i)$ is the set of index for neighbor cells of
$\Omega_{i}$.

For the hybrid meshes, the cell interface can be quadrilateral and
triangular. The triangular interface can be considered as
degenerated quadrilateral interface, where two vertexes of
quadrilateral are identical. Thus, both triangular and quadrilateral
interfaces can be calculated in a unified formulation. For the
quadrilateral interface, the four vertexes may be non-coplanar. To
calculate the numerical fluxes accurately, a curved interface need
to be considered. To be consist with the calculation of integrals
over control volume, a bilinear interpolation is used to
parameterize the cell interface as follows
\begin{align*}
\boldsymbol{X}(\eta,\zeta)=\sum_{m=1}^4\boldsymbol{x}_{m}\phi_m(\eta,\zeta),
\end{align*}
where $(\eta,\zeta)\in[-1/2,1/2]^2$, $\boldsymbol{x}_m$ is the
vertex of the interface and $\phi_m$ is the base function
\begin{align*}
\phi_1=\frac{1}{4}(1-2\eta)(1-2\zeta),
\phi_2=\frac{1}{4}(1-2\eta)(1+2\zeta),\\
\phi_3=\frac{1}{4}(1+2\eta)(1-2\zeta),
\phi_4=\frac{1}{4}(1+2\eta)(1+2\zeta),
\end{align*}
which is a degenerated form of trilinear interpolation
Eq.\eqref{trilinear}. With the parameterized cell interface, the
numerical flux can be determined by the following Gaussian
integration
\begin{align}\label{flux-qud}
\mathbb{F}_{i_p}(t)&=\iint_{\Sigma_{i_p}}F(Q, t)\text{d}\sigma=\iint_{\sigma_{i_p}}F(Q(\boldsymbol{X}(\eta,\zeta)),t)\|\boldsymbol{X}_\eta\times\boldsymbol{X}_\zeta\|\text{d}\eta\text{d}\zeta\nonumber\\
&=\sum_{k=1}^4\omega_{G_k}F_{G_k}(t)\|\boldsymbol{X}_\eta\times\boldsymbol{X}_\zeta\|_{G_k}\Delta\eta\Delta\zeta,
\end{align}
where the local orthogonal coordinate for Gaussian quadrature point
of the parameterized cell interface is distinct,
$\sigma_{i_p}=[-1/2,1/2]\times[-1/2,1/2]$, $\omega_{G_k}=1/4$ is
quadrature weight for Gaussian quadrature point
$\boldsymbol{x}_{G_k}, k=1,2,3,4$ and
$\boldsymbol{x}_{G_k}=(\pm\sqrt{3}/3,\pm\sqrt{3}/3)$.  The numerical
flux at Gaussian quadrature point $F_{G_k}(t)$ can be obtained by
taking moments of gas distribution function in the global coordinate
\begin{align}\label{flux-gau}
F_{G_k}(t)=\int\boldsymbol\psi \boldsymbol{u}\cdot (\boldsymbol{n}_1)_{G_k} f(\boldsymbol{x}_{G_k},t,\boldsymbol{u},\varsigma)\text{d}\Xi,
\end{align}
where $F_{G_k}=(F^{\rho}_{G_k},F^{\rho U}_{G_k},F^{\rho
V}_{G_k},F^{\rho W}_{G_k}, F^{\rho E}_{G_k})$ and
$(\boldsymbol{n}_1)_{G_k}$ is the local normal direction.

\subsection{Gas-kinetic solver}
In the computation, the numerical flux is usually obtained by taking
moments of gas distribution function in the local coordinate and
transferred to the global coordinate. In the local coordinate, the
gas distribution function is constructed by the integral solution of
BGK equation Eq.\eqref{bgk-3d}
\begin{align*}
f(\boldsymbol{x}_{G_k},t,\boldsymbol{u},\varsigma)=&\frac{1}{\tau}\int_0^t
g(\boldsymbol{x}',t',\boldsymbol{u},\varsigma)e^{-(t-t')/\tau}\text{d}t'
+e^{-t/\tau}f_0(-\boldsymbol{u}t,\varsigma),
\end{align*}
where the gas distribution function in the local coordinate is also
denoted as $f(\boldsymbol{x}_{G_k},t,\boldsymbol{u},\varsigma)$ for
simplicity, $\boldsymbol{u}=(u,v,w)$ is the particle velocity in the
local coordinate,
$\boldsymbol{x}_{G_k}=\boldsymbol{x}'+\boldsymbol{u}(t-t')$ is the
trajectory of particles, $f_0$ is the initial gas distribution
function, and $g$ is the corresponding equilibrium state. With the
reconstruction of macroscopic variables, the second-order gas
distribution function at the cell interface can be expressed as
\begin{equation}\label{flux-1}
\begin{split}
f(\boldsymbol{x}_{G_k},t,\boldsymbol{u},\varsigma)
&=(1-e^{-t/\tau})g_0+((t+\tau)e^{-t/\tau}-\tau)(\overline{a}_1u+\overline{a}_2v+\overline{a}_3w)g_0\\
&+(t-\tau+\tau e^{-t/\tau}){\bar{A}} g_0\\
&+e^{-t/\tau}g_r[1-(\tau+t)(a_{1}^{r}u+a_{2}^{r}v+a_{3}^{r}w)-\tau A^r)](1-H(u))\\
&+e^{-t/\tau}g_l[1-(\tau+t)(a_{1}^{l}u+a_{2}^{l}v+a_{3}^{l}w)-\tau A^l)]H(u),
\end{split}
\end{equation}
where $g_l, g_r$ are the equilibrium states  corresponding  to the
reconstructed variables $Q_l, Q_r$ at both sides of cell interface.
For the smooth flows, the simplified version of
gas distribution function can be used
\begin{align}\label{flux-2}
f(\boldsymbol{x}_{G_k},t,\boldsymbol{u},\varsigma)=g_0\big(1-\tau(\overline{a}_1u+\overline{a}_2v+\overline{a}_3w) +{\bar{A}} t\big).
\end{align}
The coefficients can be obtained by the reconstructed directional
derivatives and compatibility condition
\begin{align*}
\langle a_{m}^{l,r}\rangle&=\frac{\partial Q_{l,r}}{\partial\boldsymbol{n}_m}, ~\langle a_{1}^{l,r}u+a_{2}^{l,r}v+a_{3}^{l,r}w+A^{l,r}\rangle=0,
\end{align*}
where $m=1,2,3$, the spatial derivatives
$\displaystyle\frac{\partial Q_{l,r}}{\partial { \boldsymbol{n}_m}}$
can be determined by spatial reconstruction and $\langle...\rangle$
are the moments of the equilibrium $g$ and defined by
\begin{align*}
\langle...\rangle=\int g (...)\psi \text{d}\Xi.
\end{align*}
The equilibrium state $g_{0}$ and corresponding conservative
variables $Q_{0}$ are given by the compatibility condition
Eq.\eqref{compatibility}
\begin{align*}
\int\psi g_{0}\text{d}\Xi=Q_0=\int_{u>0}\psi g_{l}\text{d}\Xi+\int_{u<0}\psi g_{r}\text{d}\Xi.
\end{align*}
To avoid the extra reconstruction for the equilibrium part, the
spatial derivatives for equilibrium part can be determined by
\begin{align*}
\langle\overline{a}_m\rangle=\frac{\partial Q_{0}}{\partial {\boldsymbol{n}_m}}=\int_{u>0}\psi
a_{m}^{l} g_{l}\text{d}\Xi+\int_{u<0}\psi a_{m}^{r}  g_{r}\text{d}\Xi,
\end{align*}
and the temporal derivative is also given by
\begin{align*}
\langle\overline{a}_1u+\overline{a}_2v+\overline{a}_3w+\overline{A}\rangle=0.
\end{align*}
Thus, a time dependent numerical flux can be obtained with spatial
reconstruction.

\subsection{Temporal discretization}
In this paper, the unsteady and steady problems will be simulated.
To achieve the high-order temporal accuracy,  the two-stage
fourth-order temporal discretization \cite{GRP-high-1} is used in
the unsteady flows. For the Lax-Wendroff type flow solvers, such as
generalized Riemann problem (GRP) solver and gas-kinetic scheme
(GKS), the time dependent flux function can be provided. With the
temporal derivative of the flux function, a two-stage fourth-order
time accurate method was developed, and more details on the
implementation of  two-stage fourth-order method can be found in
\cite{GKS-high-1,GKS-high-4}.

For steady state problems, the implicit method is usually developed
for increasing computational efficiency. In this section, the LU-SGS
method is employed for the implicit temporal discretization. The
backward Euler method for Eq.(\ref{semi}) at $t^{n+1}$ can be
written by
\begin{equation*}
\frac{1}{\Delta t}\Delta Q_i^n=\mathcal{L}(Q_{i}^{n+1}),
\end{equation*}
where $\Delta Q_i^n=Q_{i}^{n+1}-Q_{i}^{n}$, $\Delta t$ is the time
step.For the implicit scheme,  This equation can be rewritten as
\begin{equation}\label{impEuler}
\frac{1}{\Delta t}\Delta Q_i^n-(\mathcal{L}(Q_{i}^{n+1})- \mathcal{L}(Q_{i}^{n}))=\text{Res}_i^n,
\end{equation}
where the residual is given by
\begin{align*}
 \text{Res}_i^n\triangleq\mathcal{L}(Q_{i}^{n})\approx-\frac{1}{\Delta t|\Omega_i|}\sum_{i_p\in N(i)}\int_{t_n}^{t_n+\Delta t}\mathbb{F}_{i_p}(t)\text{d}t.
\end{align*}
If $(\mathcal{L}(Q_{i}^{n+1})-\mathcal{L}(Q_{i}^{n}))$ is given by
the gas-kinetic solver, extra difficulty is introduced. To simplify
the formulation the flux splitting method is used
\begin{equation}\label{euler-scheme}
\mathcal{L}(Q_{i}^{n})=-\sum_{i_p\in N(i)}\mathbb{F}_{i_p}^{n}=-\sum_{i_p\in N(i)}\frac12(T(Q_{i_p}^{n})+T(Q_{i}^{n})-r_{i_p}(Q_{i_p}-Q_i^n))S_{i_p}
\end{equation}
where $T_i$ and $T_{i_p}$ are the Euler fluxes of cell $\Omega_i$
and $\Omega_{i_p}$,
$S_{i_p}=\sum_{i_p=1}^4\omega_{G_{i_p}}\|\boldsymbol{X}_\eta\times\boldsymbol{X}_\zeta\|_{G_{i_p}}\Delta\eta\Delta\zeta$
is the area of cell interface, and the factor $r_{i_p}$ represents
the spectral radius of the Euler flux Jacobian and ensures a
diagonal dominant matrix system with
\begin{align*}
r_{i_p}\geq U_{i_p}+a_{i_p}.
\end{align*}
Substituting Eq.\eqref{euler-scheme} into Eq.\eqref{impEuler},  we have
\begin{align}\label{jacobi}
\Big( \frac{1}{\Delta t}+\frac{1}{2}\sum_{i_p\in N(i)}r_{i_p}S_{i_p} \Big)\Delta Q_i^{n+1}
+\frac{1}{2}\sum_{i_p\in N(i)}\left[T(Q_{i_p}^{n+1})-T(Q_{i_p}^{n})-r_{i_p}\Delta Q_{i_p}^{n+1}\right]S_{i_p}=\text{Res}_i^n.
\end{align}
The set $N(i)$ can be divided into two parts, i.e. $L(i)$ and
$U(i)$, where $L(i)=\{i_p| i_p<i, i_p\in N(i)\}$ occupying in the
lower triangular area of this matrix and $U(i)=\{i_p| i_p>i, i_p\in
N(i)\}$ occupying in the upper triangular area.  The LU-SGS method,
i.e., a Gauss-Seidel iteration process, can be applied to solve the
equation above by a forward sweep step 
\begin{align*}
\Big( \frac{1}{\Delta t}+\frac{1}{2}\sum_{i_p\in N(i)}r_{i_p}S_{i_p} \Big)\Delta Q_i^{*}
+\frac{1}{2}\sum_{i_p\in L(i)}\left[T(Q_{i_p}^{*})-T(Q_{i_p}^{n})-r_{i_p}\Delta Q_{i_p}^{*}\right]S_{i_p}=\text{Res}_i^n, 
\end{align*}
and a backward sweep step
\begin{align*}
&\Big( \frac{1}{\Delta t}+\frac{1}{2}\sum_{i_p\in N(i)}r_{i_p}S_{i_p} \Big)\Delta Q_i^{n+1}
+\frac{1}{2}\sum_{i_p\in U(i)}\left[T(Q_{i_p}^{n+1})-T(Q_{i_p}^{n})-r_{i_p}\Delta Q_{i_p}^{n+1}\right]S_{i_p} \\
&=\Big( \frac{1}{\Delta t}+\frac{1}{2}\sum_{i_p\in N(i)}r_{i_p}S_{i_p} \Big)\Delta Q_i^{*}.
\end{align*}
More details can be found in \cite{GKS-high-5}. For the  LU-SGS method, the forward sweep step and the backward sweep step 
need to be performed sequentially, which reduce the efficiency of computation. 
In order to make the implicit method more suitable for parallel computation,
Jacobi iteration is applied to solve Eq.(\ref{jacobi}).
The inner iteration is performed by 
\begin{align*}
\Big( \frac{1}{\Delta t}+\frac{1}{2}\sum_{i_p\in N(i)}r_{i_p}S_{i_p} \Big)\Delta Q_i^{(k)}
+\frac{1}{2}\sum_{i_p\in N(i)}\left[T(Q_{i_p}^{(k-1)})-T(Q_{i_p}^{n})-r_{i_p}\Delta Q_{i_p}^{(k-1)}\right]S_{i_p}=\text{Res}_i^n,   
\end{align*}
where $1\le k\le k_{\max}$, $k$ is the number of inner iterations and the initial solution $\Delta Q_i^{(0)}$ is given by
\begin{align*}
\Big( \frac{1}{\Delta t}+\frac{1}{2}\sum_{i_p\in N(i)}r_{i_p}S_{i_p} \Big)\Delta Q_i^{(0)}=\text{Res}_i^n.
\end{align*}
Thus, Eq.(\ref{jacobi}) is solved by $\Delta Q_i^{n+1}=\Delta Q_i^{(k_{\max})}$.
Obviously, the above processes are easy to implement in parallel.

\section{GPU architecture and code design}
The computation of HGKS is mainly based on the central processing unit (CPU) code. 
To improve the efficiency, the OpenMP directives and message passing interface (MPI) are used for parallel computation. 
However, the CPU computation is usually limited in the number of threads which are handled in parallel.
Graphics processing unit (GPU) is a form of hardware acceleration, which is originally developed for graphics 
manipulation and execute highly-parallel computing tasks. Currently,  GPU has gained significant popularity in large-scale scientific computations.  
In the previous study \cite{GKS-GPU},  discontinuous Galerkin-based HGKS was implemented with GPU 
using compute unified device architecture (CUDA), and the efficiency is improved greatly. 
In this work,  the high-order gas-kinetic scheme on hybrid unstructured meshes is implemented with GPU as well.

\begin{figure}[!h]
\centering
\includegraphics[width=0.5\textwidth]{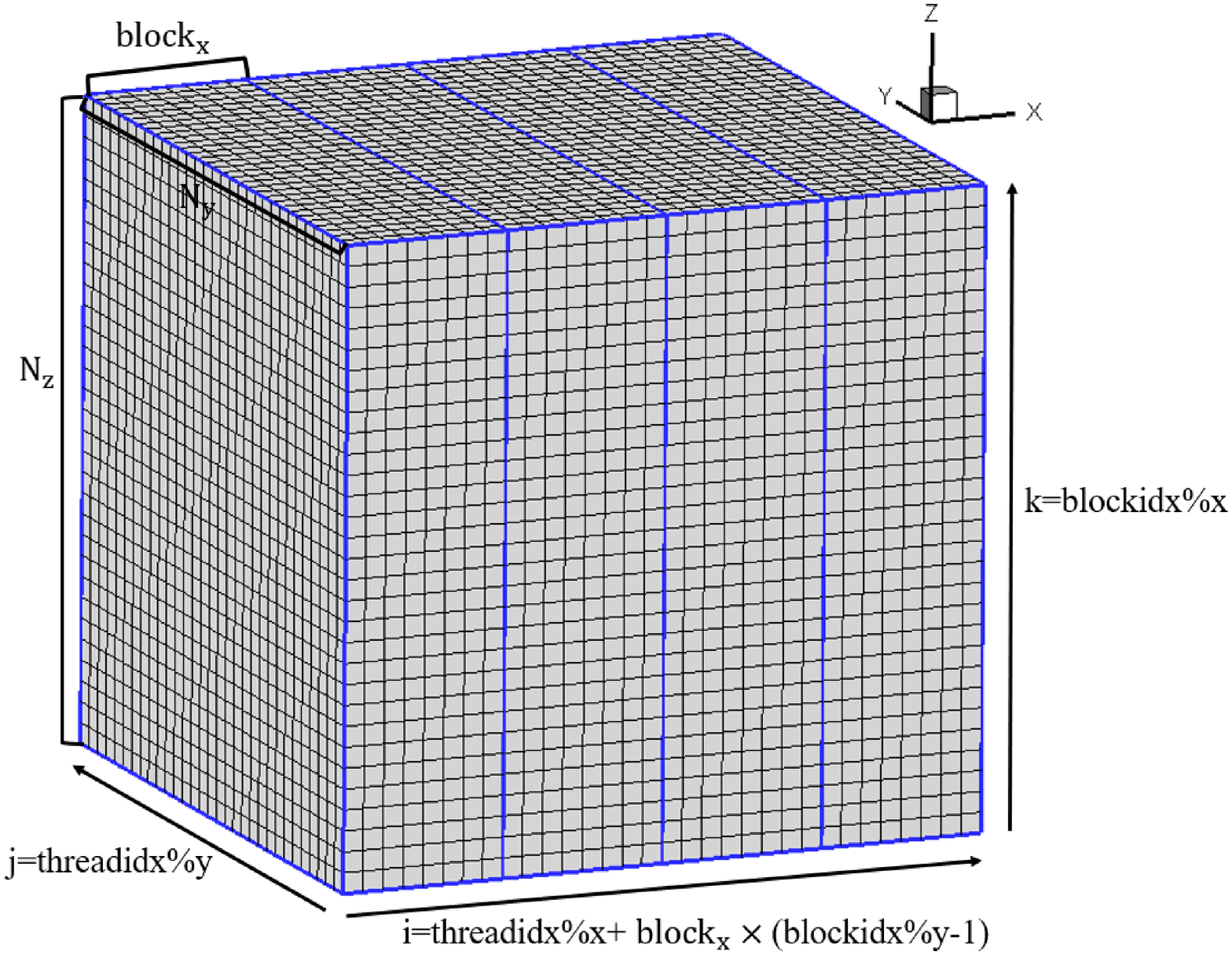}
\caption{\label{thread-block-structured} Connection between threads and structured cells.}
\centering
\includegraphics[width=0.8\textwidth]{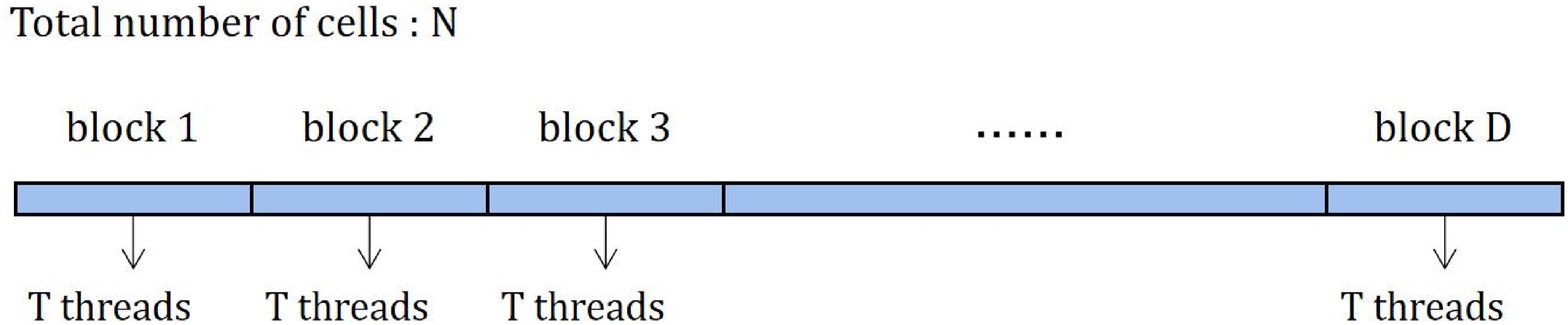} 
\caption{\label{thread-block-unstructured} Connection between threads and unstructured cells.}
\end{figure}

The CUDA threads are organized into thread blocks, and thread blocks constitute a thread grid. 
The thread grid may be seen as a computational structure on GPUs, which can be corresponded to computational cells. 
In the previous study \cite{GKS-GPU}, the three-dimensional structured meshes are used,  and a simple partition of mesh is given in Figure.\ref{thread-block-structured}.
Assume that the total number of meshes is
$N_x\times N_y\times N_z$, and the three-dimensional  computational
domain is divided into $\text{block}_x$ parts in $x$-direction,
where $\text{block}_x$ is an integer defined according to tests. The
variables "dimGrid" and "dimBlock" are defined to set a
two-dimensional grid, which consist of  two-dimensional blocks as
follows
\begin{align*}
{\rm dimGrid}&={\rm dim3}(N_z, N_x/\text{block}_x, 1), \\
 {\rm dimBlock}&={\rm dim3}(\text{block}_x,N_y,1).
\end{align*}
It is natural to assign one thread to complete computations of a cell
$K_{ijk}$, and the one-to-one correspondence of thread block index and
cell index $(i,j,k)$ is simpliy given as follows
\begin{align*}
i&=\text{threadidx}\%\text{x}+\text{block}_x*(\text{blockidx}\%\text{y}-1),\\
j&=\text{threadidx}\%\text{y}, \\
k&= \text{blockidx}\%\text{x}.
\end{align*}
For the unstructured meshes, the data for cells, interfaces and nodes are stored one-dimensionally. 
Compared with the structured meshes, the partition for blocks becomes simpler.  As shown in in Figure.\ref{thread-block-unstructured}, 
the cells can be divided into $D$ blocks, where $N$ is the total number of cells. 
For each block, it contains some threads and 
the maximum number of threads in block is 1024. It is supposed that each block contains $T$ threads and $N=D\times T$ for example. 
If $N$ is not divisible by $T$, an extra block is needed.  The
variables "dimGrid" and "dimBlock" are defined as
\begin{align*}
{\rm dimGrid} &= {\rm dim3}(\text{D},1,1),\\
{\rm dimBlock} &={\rm dim3}(\text{T},1,1).
\end{align*}
Thus, the one-to-one correspondence for cell and thread ID  for parallel computation is  given by
\begin{equation*}
    i=\text{threadidx}\%\text{x}+\text{T}*(\text{blockidx}\%\text{x}-1).
\end{equation*} 
Main parts of the HGKS code are labeled blue in Algorithm.\ref{GPU-algorithm},  and the way to implement these parts in parallel on GPU is using kernels. 
The kernel is a subroutine, which executes at the same time by many threads on GPU, and 
the kernel for WENO reconstruction is given as an example to show the main idea of CUDA program. 
The Nvidia GPU is consisted of multiple streaming multiprocessors (SMs). Each block of grid is distributed to one SM, and the threads
of block are executed in parallel on SM. The executions are implemented automatically by GPU. Thus, the GPU code can be
implemented with specifying kernels and grids.

\begin{algorithm}[H]
\setstretch{1.25}
\begin{algorithmic} 
\STATE  
\textbf{\color{blue} Initial condition}
\WHILE {TIME $ \leq $ TSTOP}
\STATE \textbf{STEP 1} : \textbf{\color{blue} Calculation of time step}
\STATE \textbf{STEP 2} : \textbf{\color{blue} WENO reconstruction}
\STATE $\text{block}_N$ $=$ ${\rm N}/T$  \quad \textcolor{SpringGreen4}{\% ${\rm N}$ is the number of cell interfaces} 
\STATE dimGrid $=$ dim3($\text{block}_N$,1,1)
\STATE dimBlock $=$ dim3(T,1,1)  
\STATE \textbf{CALL} WENO$<<<$dimGrid,dimBlock$>>>$   \quad  \textcolor{SpringGreen4}{\%  kernel for reconstruction and flux calculation} 
\STATE istat=cudaDeviceSynchronize()  \quad \textcolor{SpringGreen4}{\% blocking the current device until all preceeding tasks have completed} 
\STATE \textbf{STEP 3} : \textbf{\color{blue} Update of flow variables} (Jacobi iteration or two-stage method)
\ENDWHILE \\
\textbf{\color{blue} Output of flow field}
\end{algorithmic} 
\caption{\label{GPU-algorithm} Program for unstructured HGKS}
\end{algorithm}

\section{Numerical tests}
In this section, numerical tests for both inviscid and viscous flows
will be presented to validate the current scheme. For the inviscid
flows, the collision time $\tau$ takes
\begin{align*}
\tau=\epsilon \Delta t+C\displaystyle|\frac{p_l-p_r}{p_l+p_r}|\Delta t,
\end{align*}
where $\epsilon=0.01$ and $C=1$. For the viscous flows, we have
\begin{align*}
\tau=\frac{\mu}{p}+C \displaystyle|\frac{p_l-p_r}{p_l+p_r}|\Delta t,
\end{align*}
where $p_l$ and $p_r$ denote the pressure on the left and right
sides of the cell interface,  $\mu$ is the dynamic viscous
coefficient and $p$ is the pressure at the cell interface. In smooth
flow regions, it will reduce to
\begin{equation*}
\tau=\frac{\mu}{p}
\end{equation*}
In the computation, the poly gas is used and the specific heat ratio
takes $\gamma=1.4$. For the numerical examples, the two-stage method
is used for the unsteady problems, and the LU-SGS method is used for
the steady problems.

\begin{figure}[!h]
\centering
\includegraphics[width=0.475\textwidth]{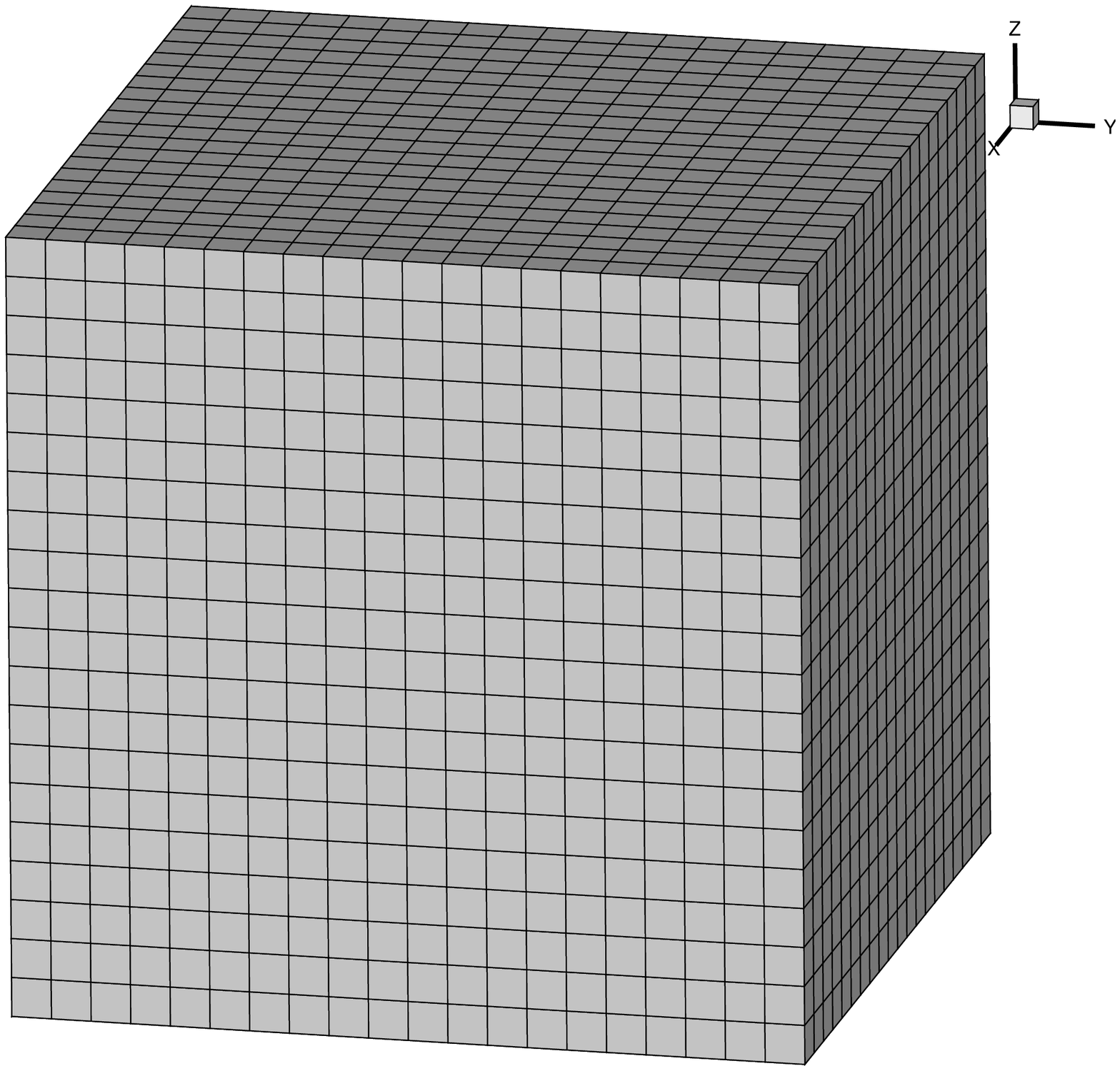}
\includegraphics[width=0.475\textwidth]{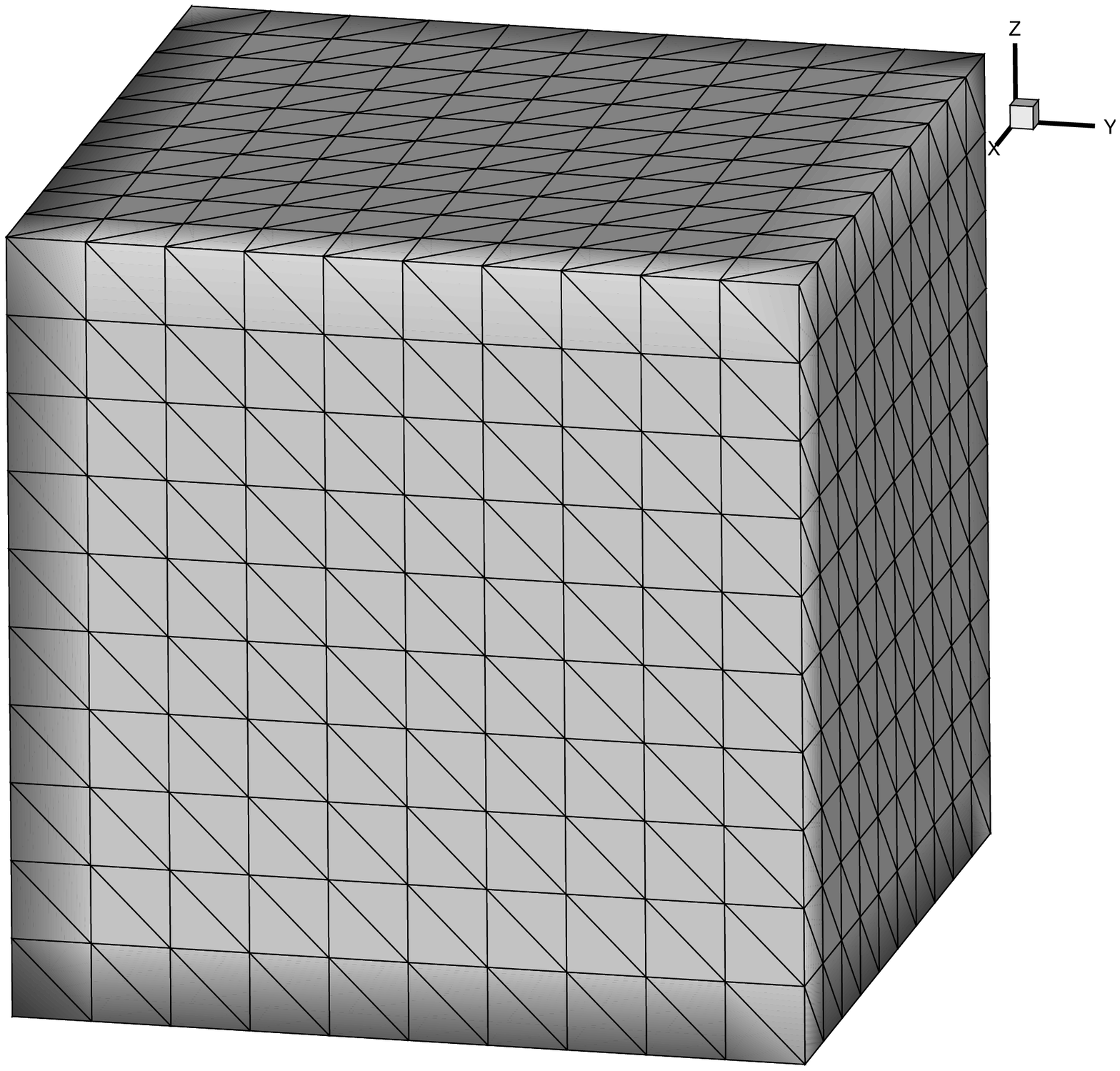}
\caption{\label{Accuracy-3d}  Accuracy test: the three-dimensional
hexahedron mesh with $20^3$ cells (left) and tetrahedron mesh
$10^3\times6$ cells (right) for accuracy test.}
\end{figure}

 \begin{table}[!h]
\begin{center}
\def\temptablewidth{0.75\textwidth}
{\rule{\temptablewidth}{1.0pt}}
\begin{tabular*}{\temptablewidth}{@{\extracolsep{\fill}}c|cc|cc}
mesh     & $L^1$ error  &    Order    &  $L^2$ error &  Order    \\
\hline
$10^3$     & 5.2689E-01  &      ~   &   2.0675E-01  &   ~            \\
$20^3$     & 1.0636E-01  &   2.3085 &   4.1798E-02   &   2.3064  \\
$40^3$     & 1.4729E-02  &   2.8521 &   5.7712E-03   &  2.8565 \\
$80^3$     & 1.8840E-03  &   2.9668 &   7.3678E-04   &  2.9695 \\
\end{tabular*}
{\rule{\temptablewidth}{1.0pt}}
\end{center}
\caption{\label{tab-3d-1} Accuracy test: errors and orders of
accuracy with uniform hexahedron meshes.}
\begin{center}
\def\temptablewidth{0.75\textwidth}
{\rule{\temptablewidth}{1.0pt}}
\begin{tabular*}{\temptablewidth}{@{\extracolsep{\fill}}c|cc|cc}
mesh     &   $L^1$ error  &Order &   $L^2$ error  &Order \\
\hline
$5^3\times6$    &    3.7642E-01   &    ~~        &  1.4740E-01 &   ~~      \\
$10^3\times6$  &    6.0508E-02   &   2.6371  &   2.4031E-02 &  2.6167 \\
$20^3\times6$  &    7.8292E-03   &   2.9502  &   3.1387E-03 &  2.9366 \\
$40^3\times6$  &    9.9153E-04   &   2.9811  &   4.0739E-04 &  2.9456 \\
\end{tabular*}
{\rule{\temptablewidth}{1.0pt}}
\end{center}
\caption{\label{tab-3d-2} Accuracy tests: errors and orders of
accuracy with tetrahedron meshes.}
\end{table}

\subsection{Accuracy test}
In this case, the three-dimensional advection of density
perturbation is used to test the order of accuracy with the
hexahedron and tetrahedron meshes. For this case, the computational
domain is $[0,2]\times[0,2]\times[0,2]$ and the initial condition is
given as follows
\begin{align*}
\rho_0&(x, y, z)=1+0.2\sin(\pi(x+y+z)),~p_0(x,y,z)=1,\\
&U_0(x,y,z)=1,~V_0(x,y,z)=1,~W_0(x,y,z)=1.
\end{align*}
The periodic boundary condition is applied on all boundaries, and the
exact solution is
\begin{align*}
\rho(x,y&,z,t)=1+0.2\sin(\pi(x+y+z-t)),~p(x,y,z,t)=1,\\
&U(x,y,z,t)=1,~V(x,y,z,t)=1,~W(x,y,z,t)=1.
\end{align*}
In smooth flow regions, $\tau=0$ and the gas-distribution function
reduces to
\begin{align*}
f(\boldsymbol{x}_{m_1,m_2},t,\boldsymbol{u},\varsigma)=g(1+At).
\end{align*}
For the hexahedron meshes, the uniform meshes with $\Delta x=\Delta
y=\Delta z=2/N$ are used. The $L^1$ and $L^2$ errors and orders of
accuracy at $t=2$ are presented in Tab.\ref{tab-3d-1}, where the
expected order of accuracy is achieved. For the tetrahedron meshes,
a series of meshes with $6 \times N^3$ cells are used, where every
cubic is divided into six tetrahedron cells. The $L^1$ and $L^2$
errors and orders of accuracy at $t=2$ are presented in
Tab.\ref{tab-3d-2}, where the expected third-order of accuracy is
also achieved. The mesh and the density distributions for the
hexahedron mesh with $N=32$, and for the tetrahedron mesh with
$N=20$  are given in Fig.\ref{Accuracy-3d}.

\begin{figure}[!h]
\centering
\includegraphics[width=0.75\textwidth]{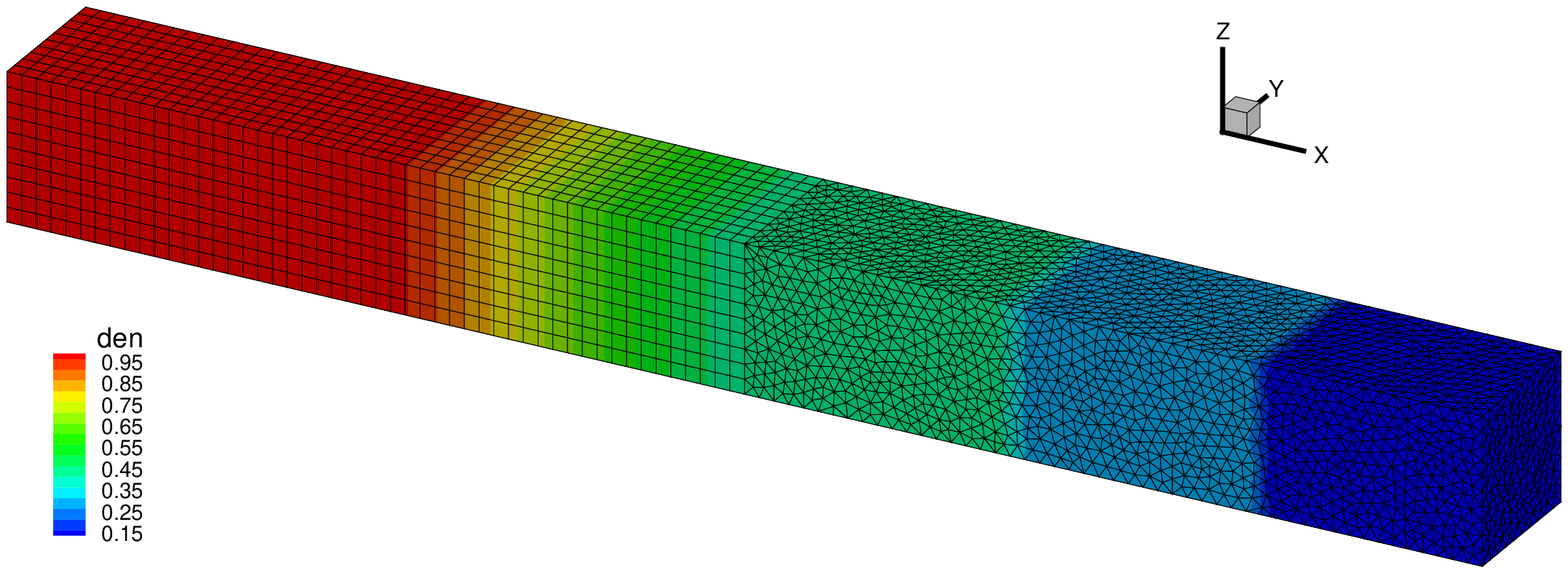}
\includegraphics[width=0.75\textwidth]{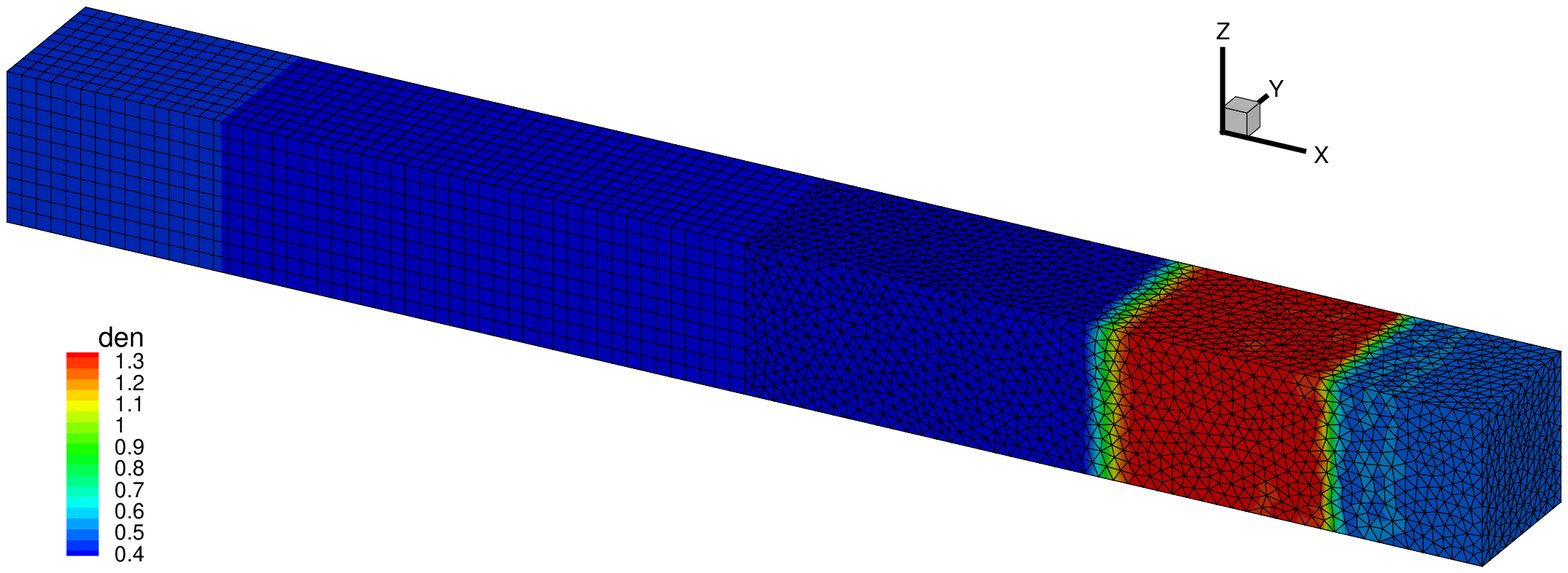}
\caption{\label{Riemann-3d} Riemann problem: the three-dimensional
mesh and density distributions for Sod problem (left) at $t=0.2$ and
Lax problem (right) at $t=0.16$ for the hybrid mesh.}
\end{figure}

\begin{figure}[!h]
\centering
\includegraphics[width=0.475\textwidth]{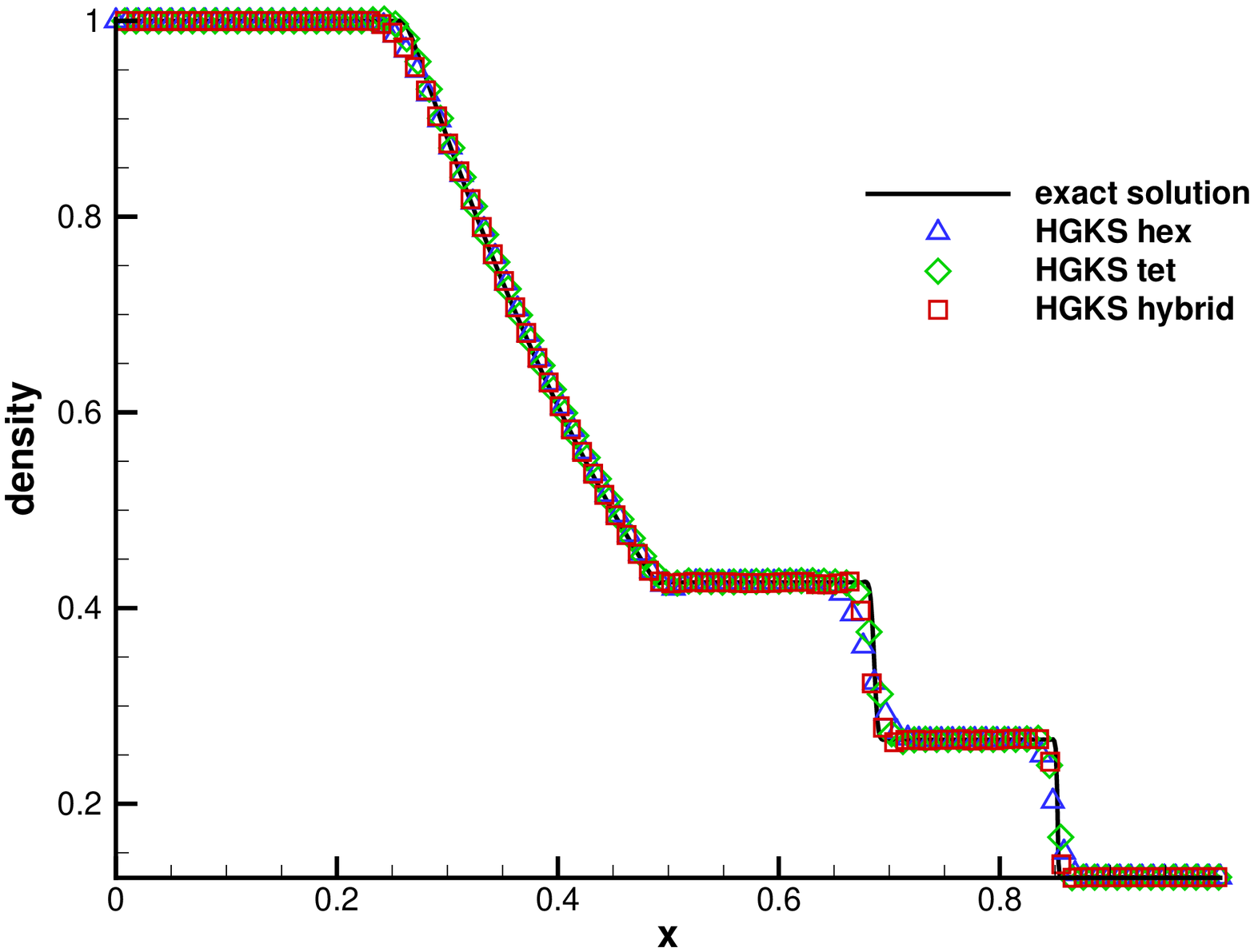}\includegraphics[width=0.475\textwidth]{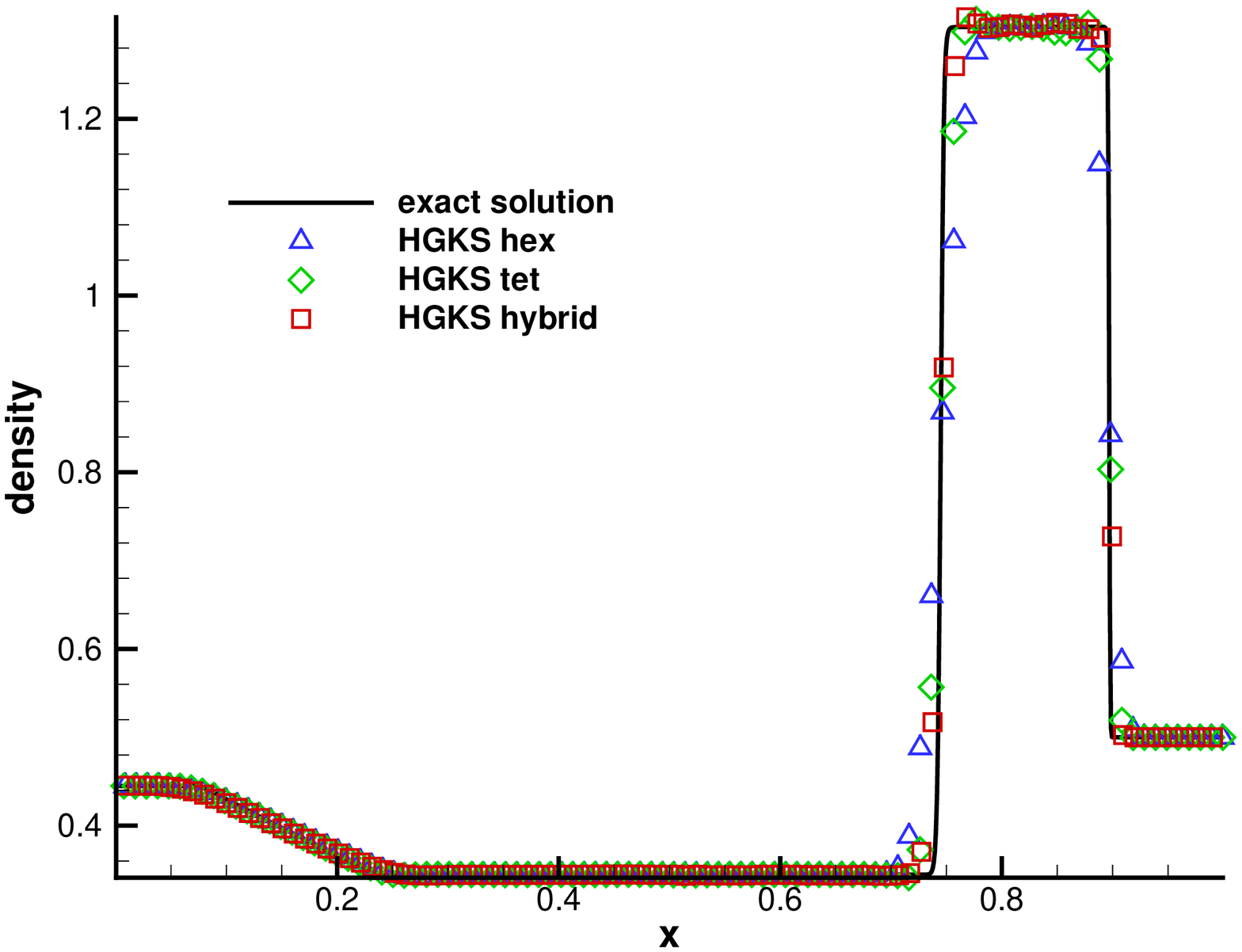}
\includegraphics[width=0.475\textwidth]{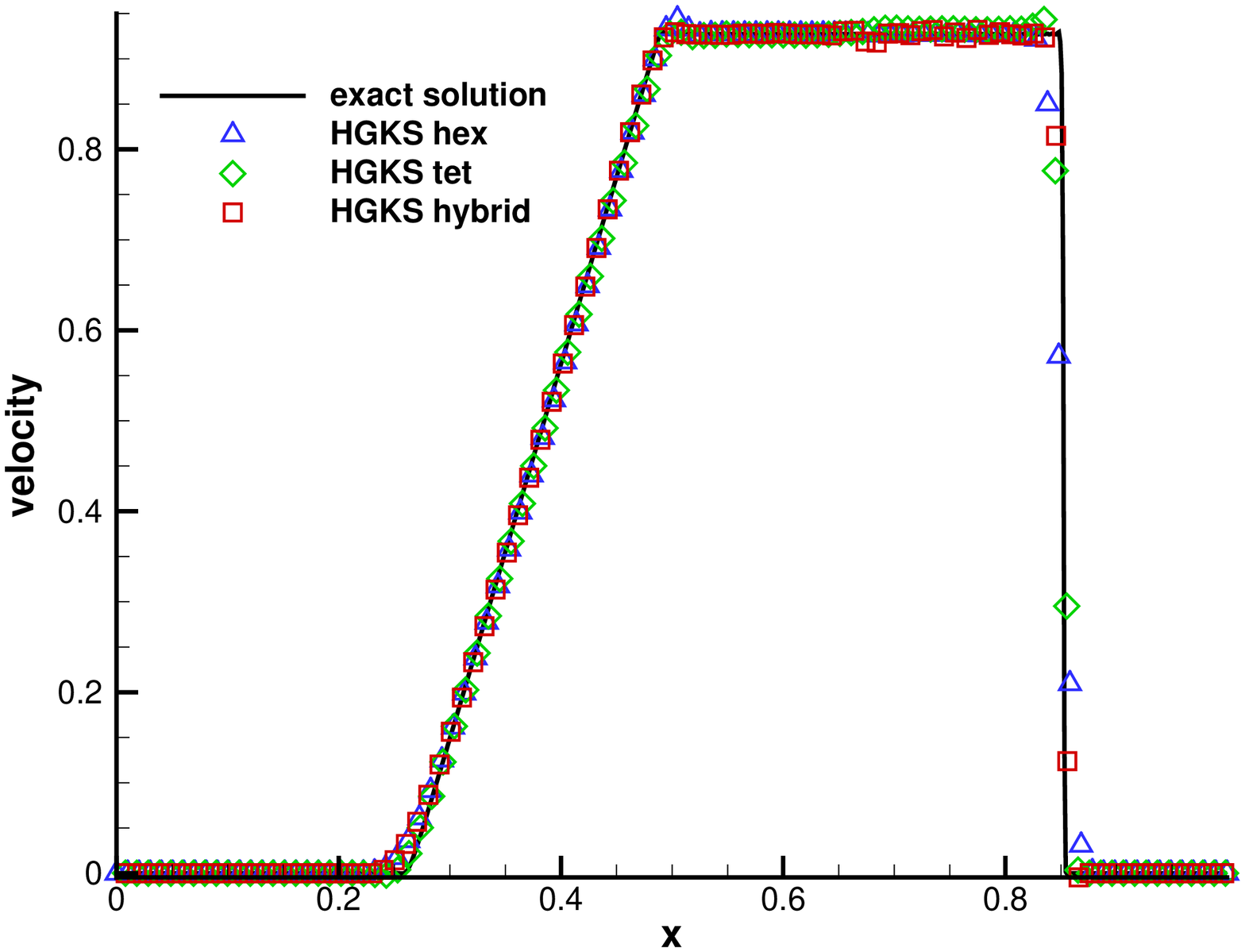}\includegraphics[width=0.475\textwidth]{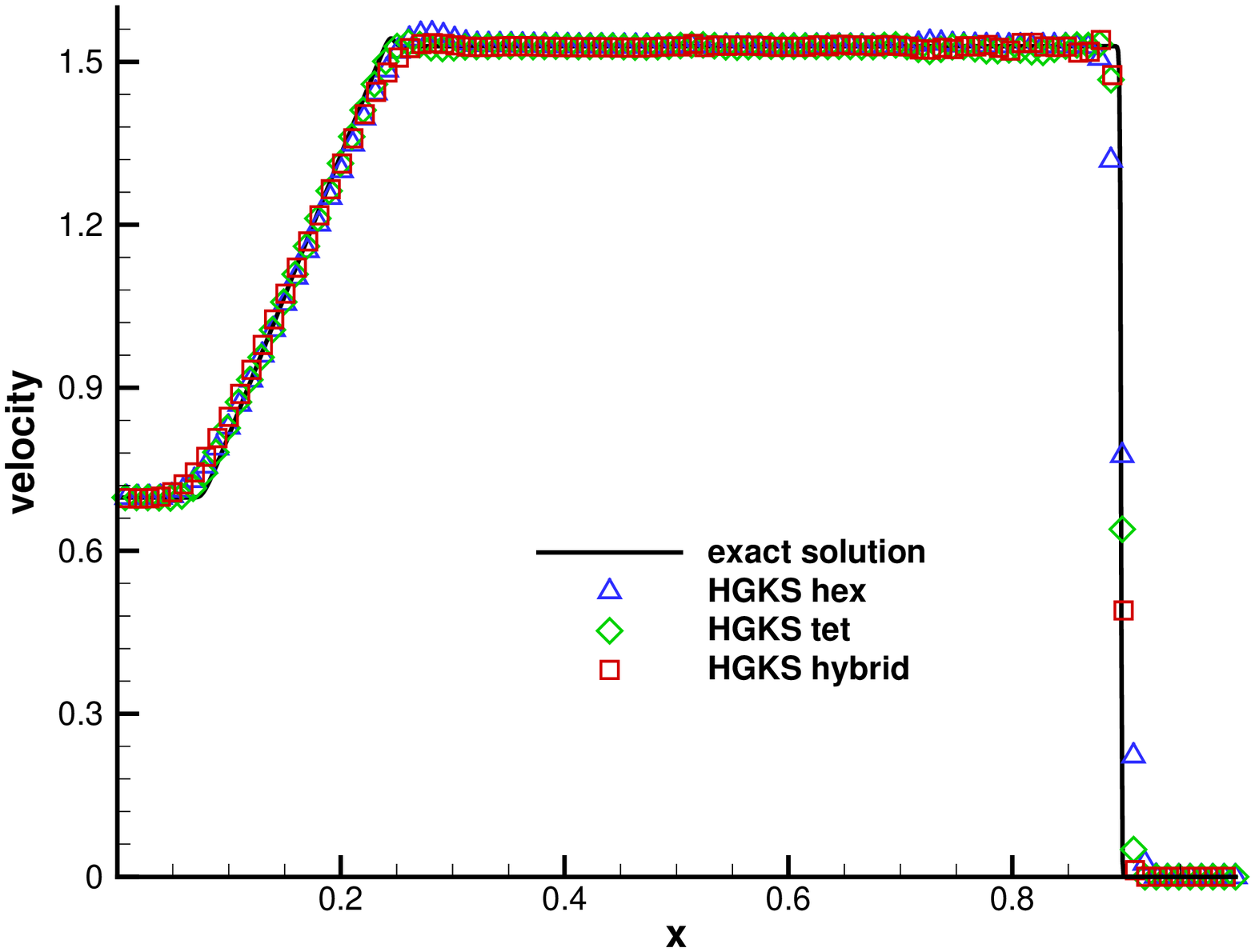}
\includegraphics[width=0.475\textwidth]{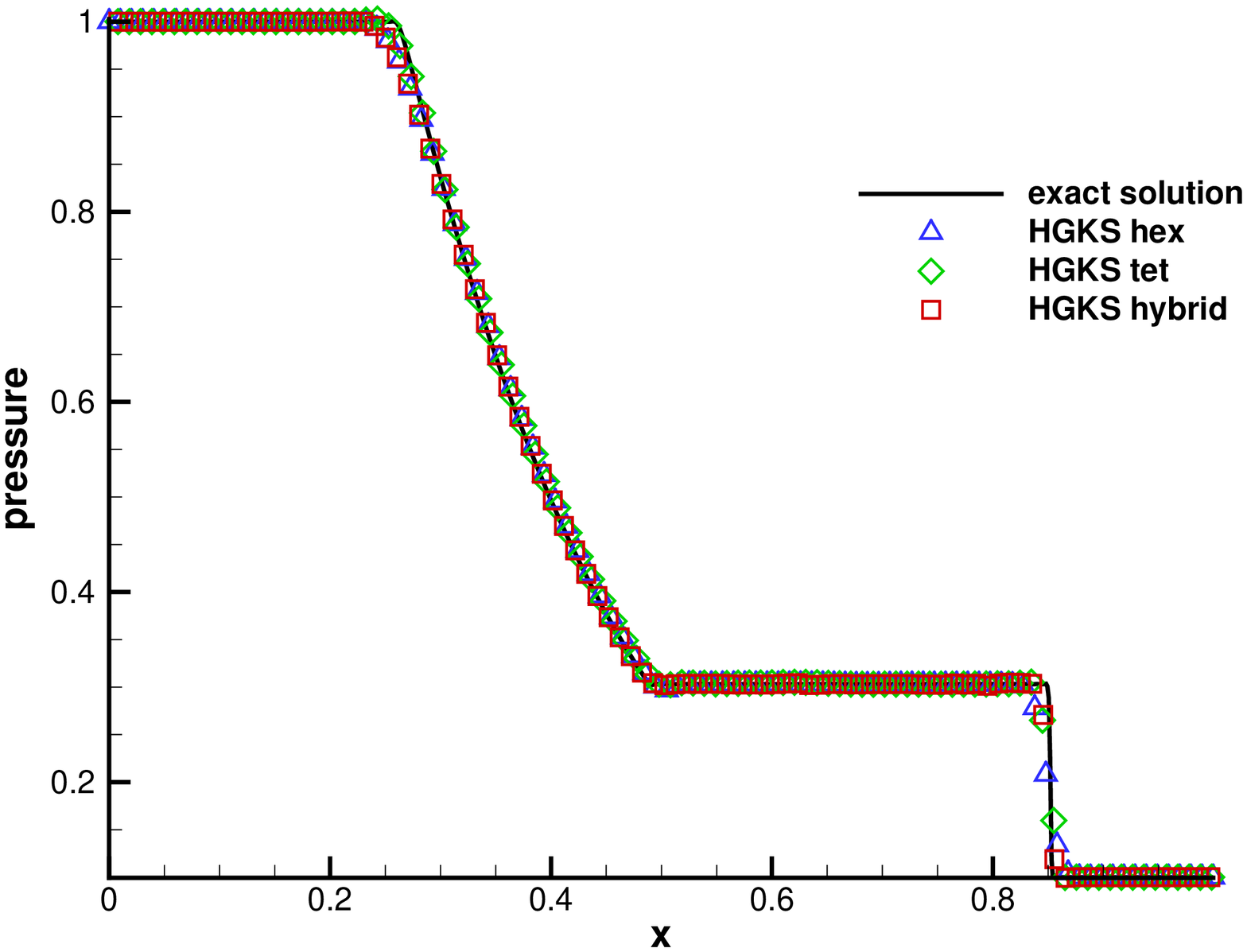}\includegraphics[width=0.475\textwidth]{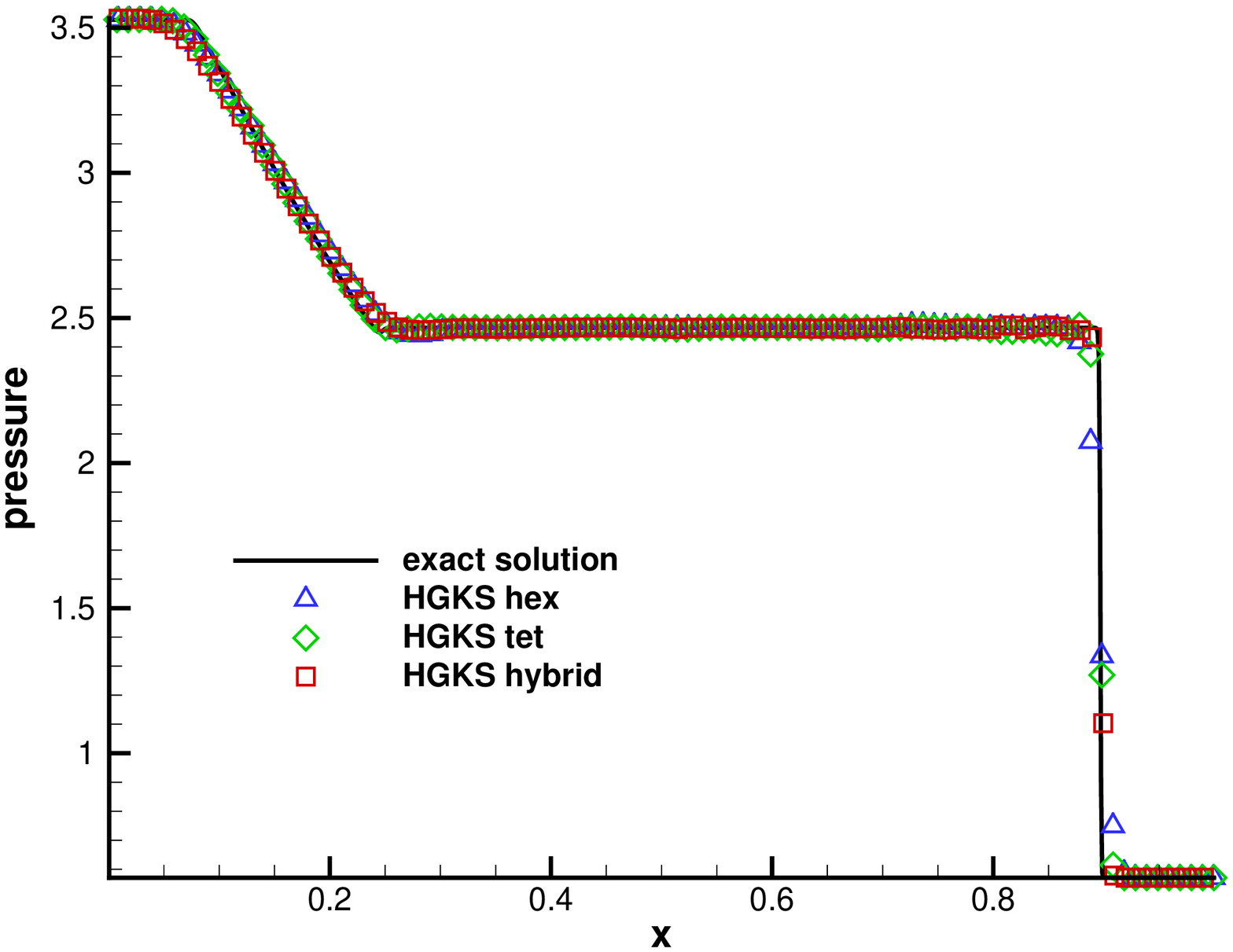}
\caption{\label{Riemann-1d} Riemann problem: density, velocity and
pressure distributions for Sod problem (left) at $t=0.2$ and Lax
problem (right) at $t=0.16$ at the center horizontal line of
hexahedral, tetrahedral and hybrid meshes.}
\end{figure}

\subsection{Riemann problems}
In this case, two one-dimensional Riemann problems are tested by the
third-order WENO scheme on the hybrid meshes which consist of
tetrahedral, pyramidal and hexahedral cells. The first one is the
Sod problem, and the initial condition is given as follows
\begin{equation*}
(\rho,U,V,W,p) = \begin{cases}
(1,0,0,0,1),  0\leq x<0.5,\\
(0.125, 0,0,0, 0.1), 0.5\leq x\leq1.
\end{cases}
\end{equation*}
The second one is the Lax problem, and the initial condition is
given as follows
\begin{equation*}
(\rho,U,V,W,p) = \begin{cases}
(0.445, 0.698,0,0, 3.528),   0\leq x<0.5,\\
(0.5, 0, 0,0,0.571), 0.5\leq x\leq1.
\end{cases}
\end{equation*}
For these two cases, the computational domain is
$[0,1]\times[0,0.1]\times[0.1]$. This case is tested by hybrid mesh
with total 44451 cells, including 5000 hexahedron cells, 100 prism
cells and 39351 tetrahedron cells. The mesh size is $h=1/100$. As
comparison, these two cases are also tested with hexahedral and
tetrahedral meshes. The uniform mesh with $100\times10\times10$
cells for the hexahedral mesh and 59003 cells with mesh size
$h=1/100$ for the tetrahedral mesh are used. Non-reflection boundary
condition is adopted at the left and right boundaries of the
computational domain, and reflection boundary condition is adopted
at other boundaries of the computational domain. The
three-dimensional mesh and density distributions for Sod problem at
$t=0.2$ and Lax problem at $t=0.16$ are given in
Fig.\ref{Riemann-3d} for hybrid mesh. The numerical results of
density, velocity and pressure for the Sod problem at $t=0.2$ and
for the Lax problem at $t=0.16$ are presented in
Fig.\ref{Riemann-1d} with $x\in[0,1], y=z=0$ for hexahedral,
tetrahedral and hybrid meshes. The exact solutions are also given.
The numerical results agree well with the exact solution, and the
discontinuities are well resolved by the current scheme. As
expected, the tetrahedral mesh contains more cells and resolves the
discontinuities better than the hexahedral mesh.

\begin{figure*}[!h]
\centering
\includegraphics[width=0.475\textwidth]{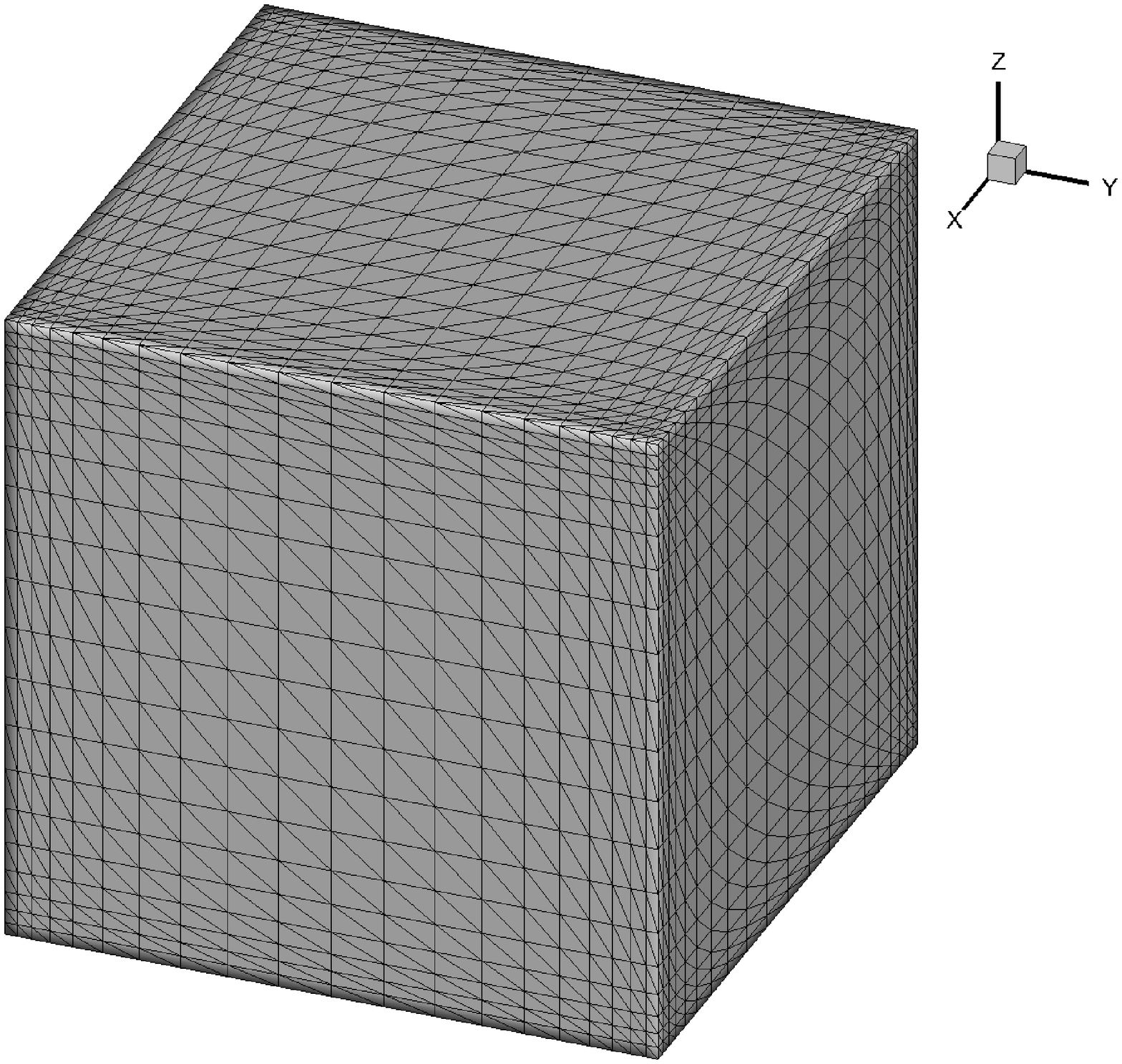}
\includegraphics[width=0.475\textwidth]{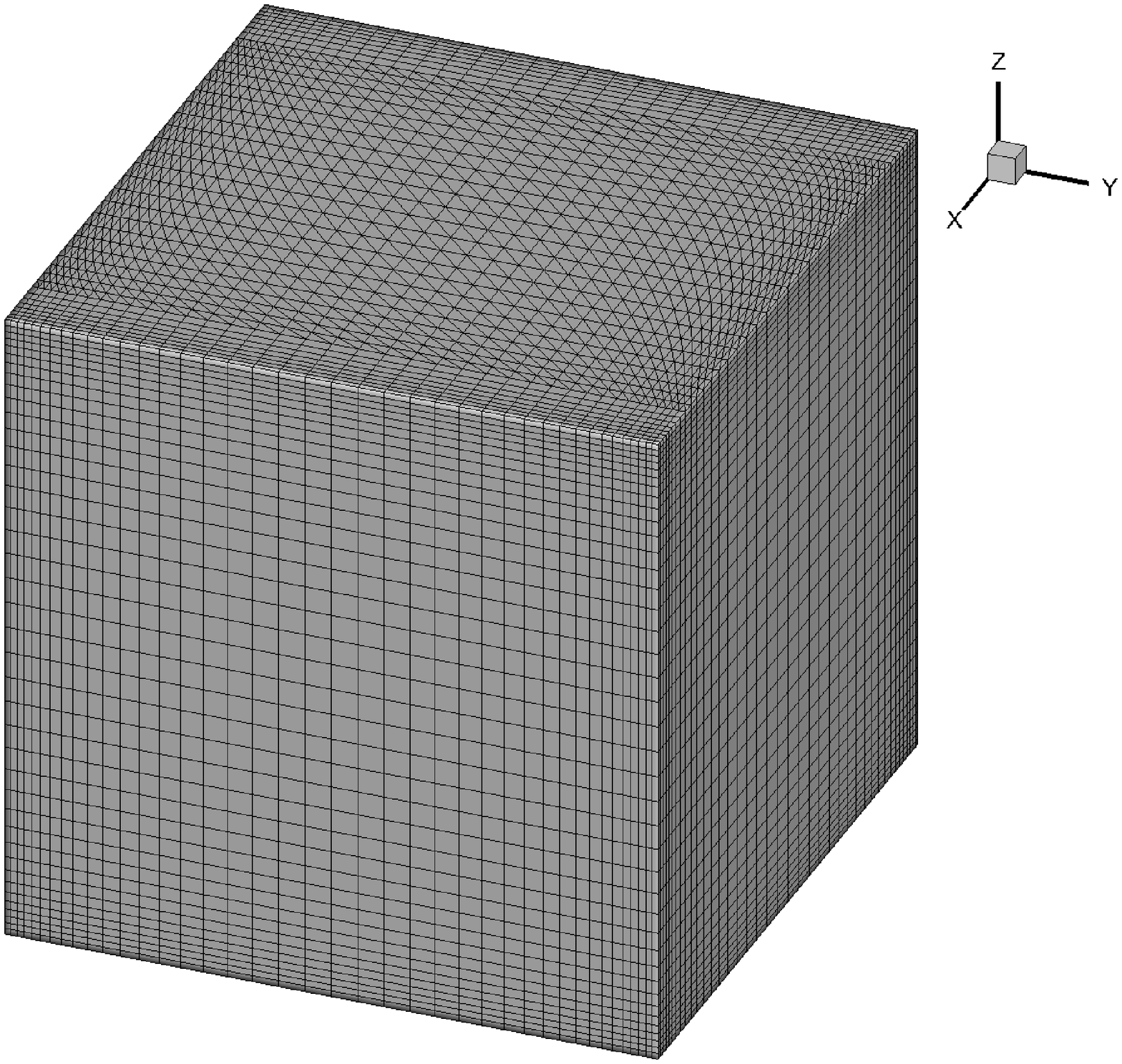}
\caption{\label{cavity-mesh} Lid-driven cavity flow: the local computational
mesh and pressure distribution.}
\end{figure*}

\begin{figure*}[!h]
\centering
\includegraphics[width=0.45\textwidth]{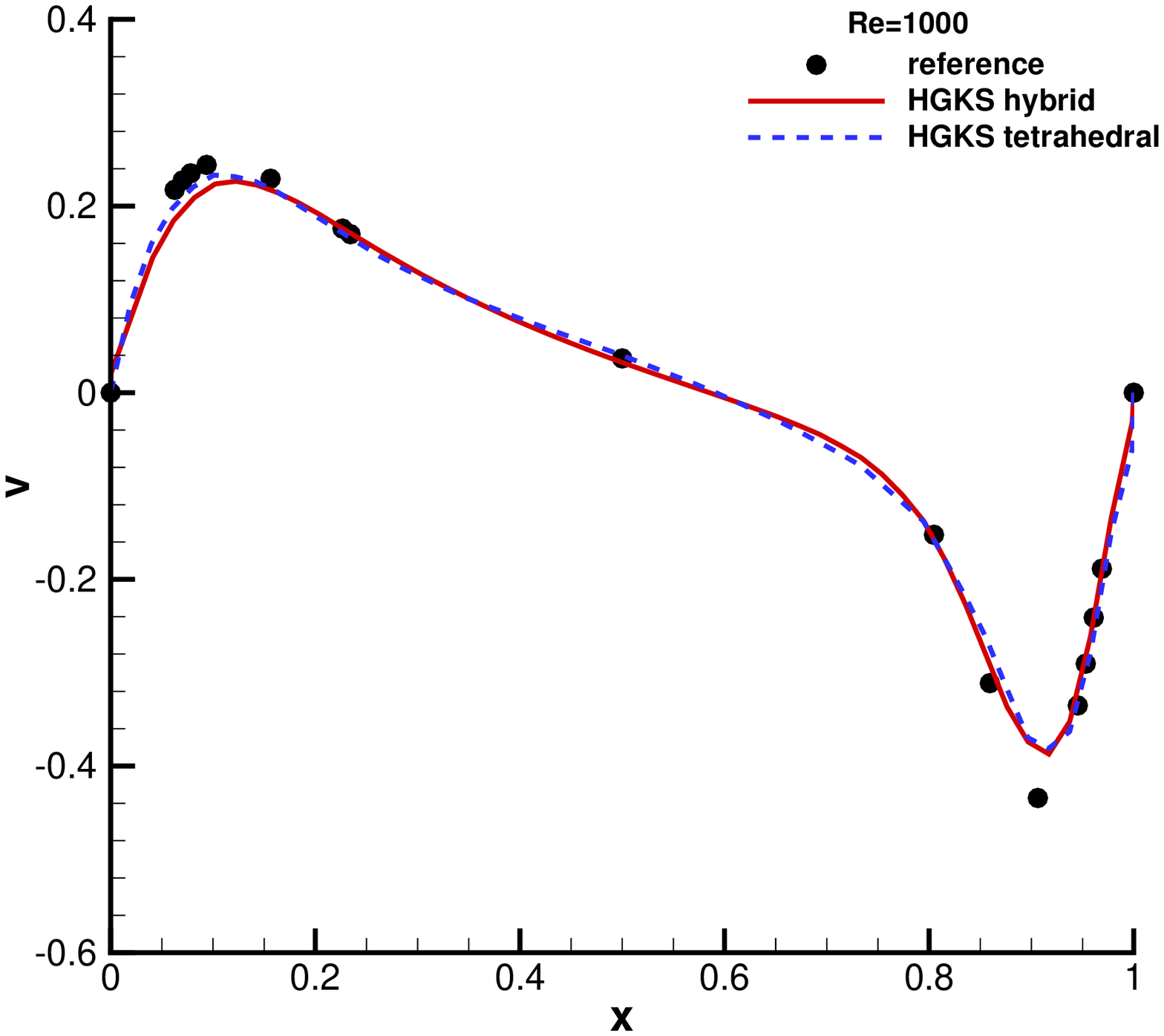}
\includegraphics[width=0.45\textwidth]{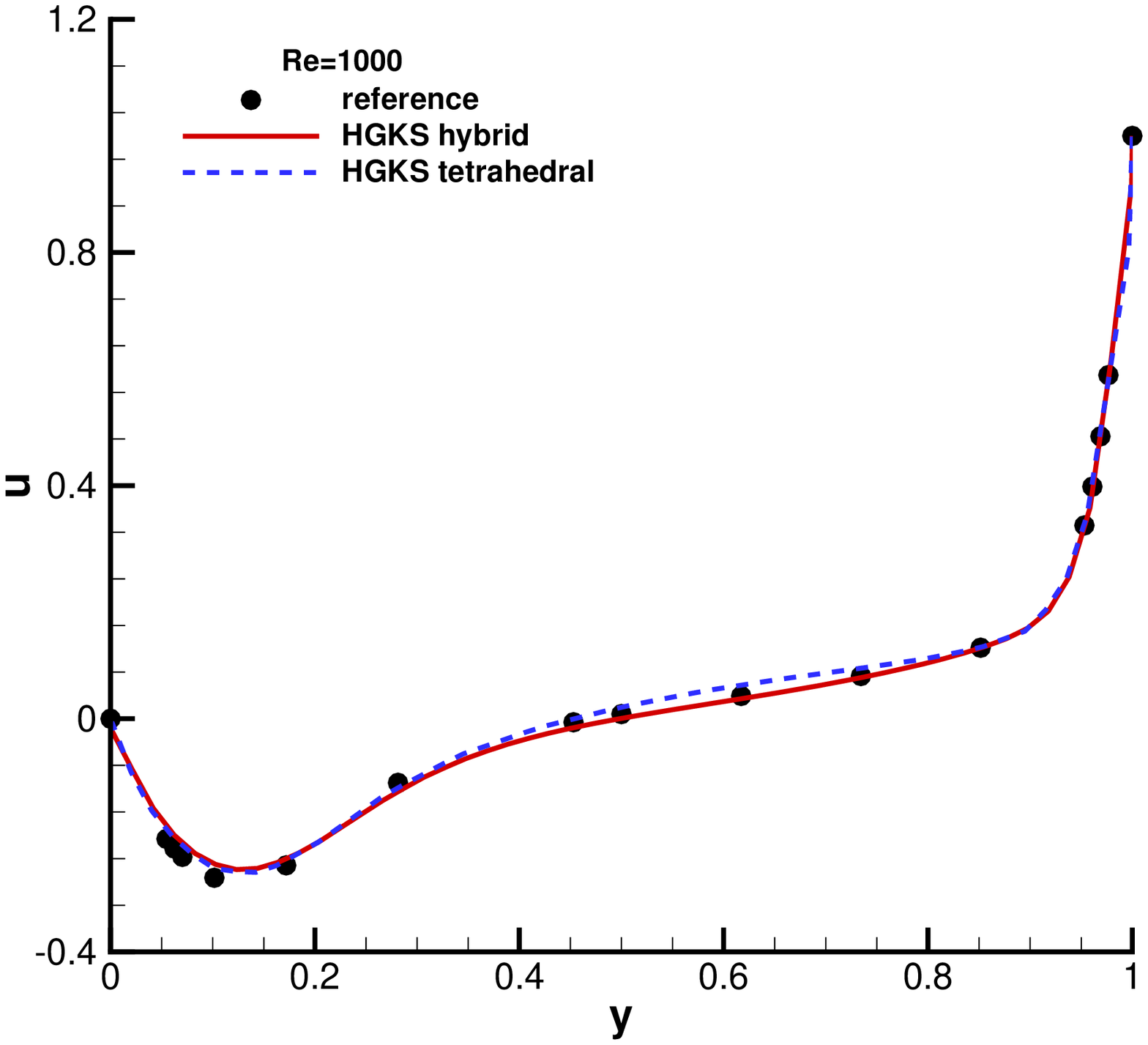}\\
\includegraphics[width=0.45\textwidth]{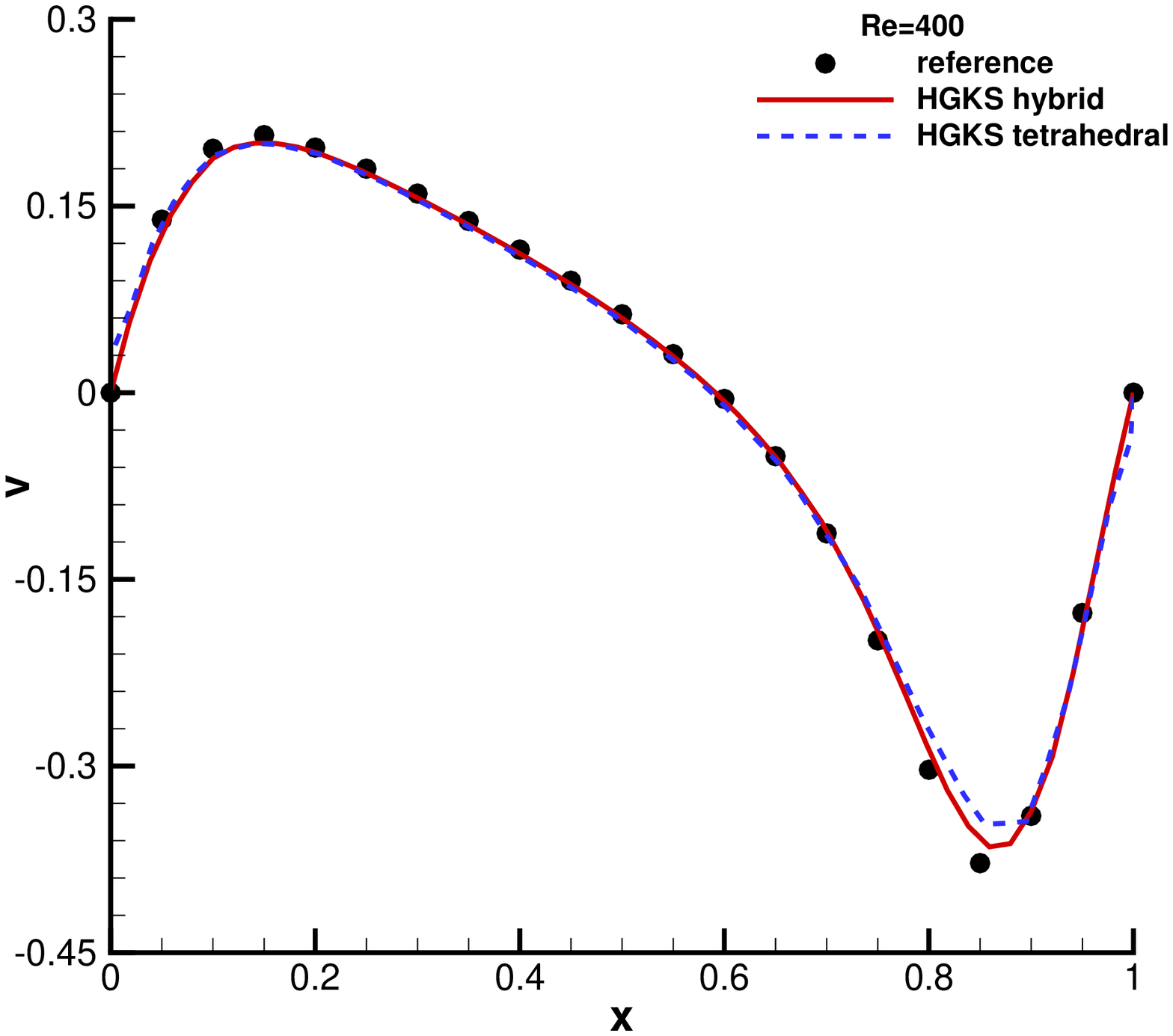}
\includegraphics[width=0.45\textwidth]{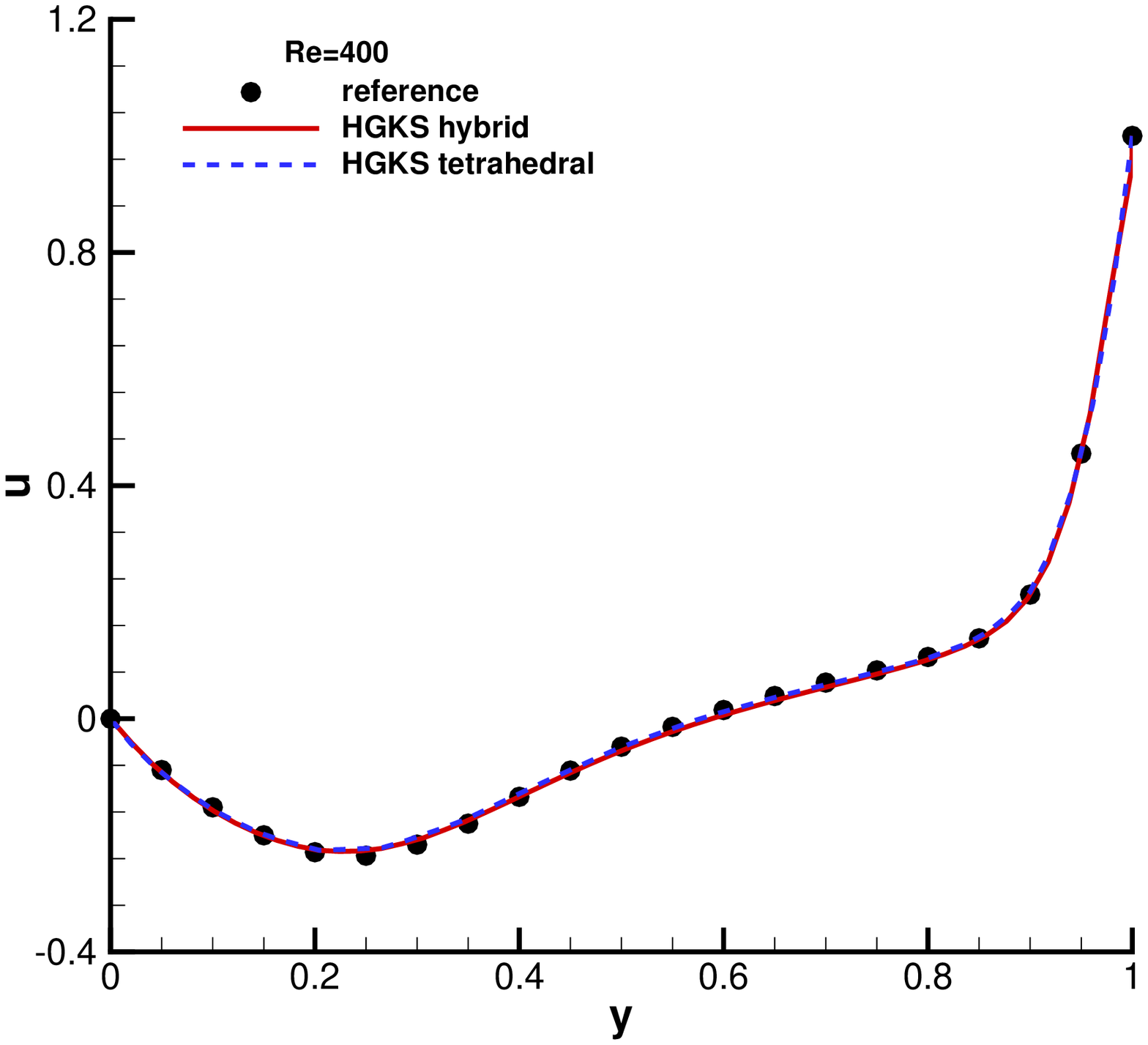} \\
\includegraphics[width=0.45\textwidth]{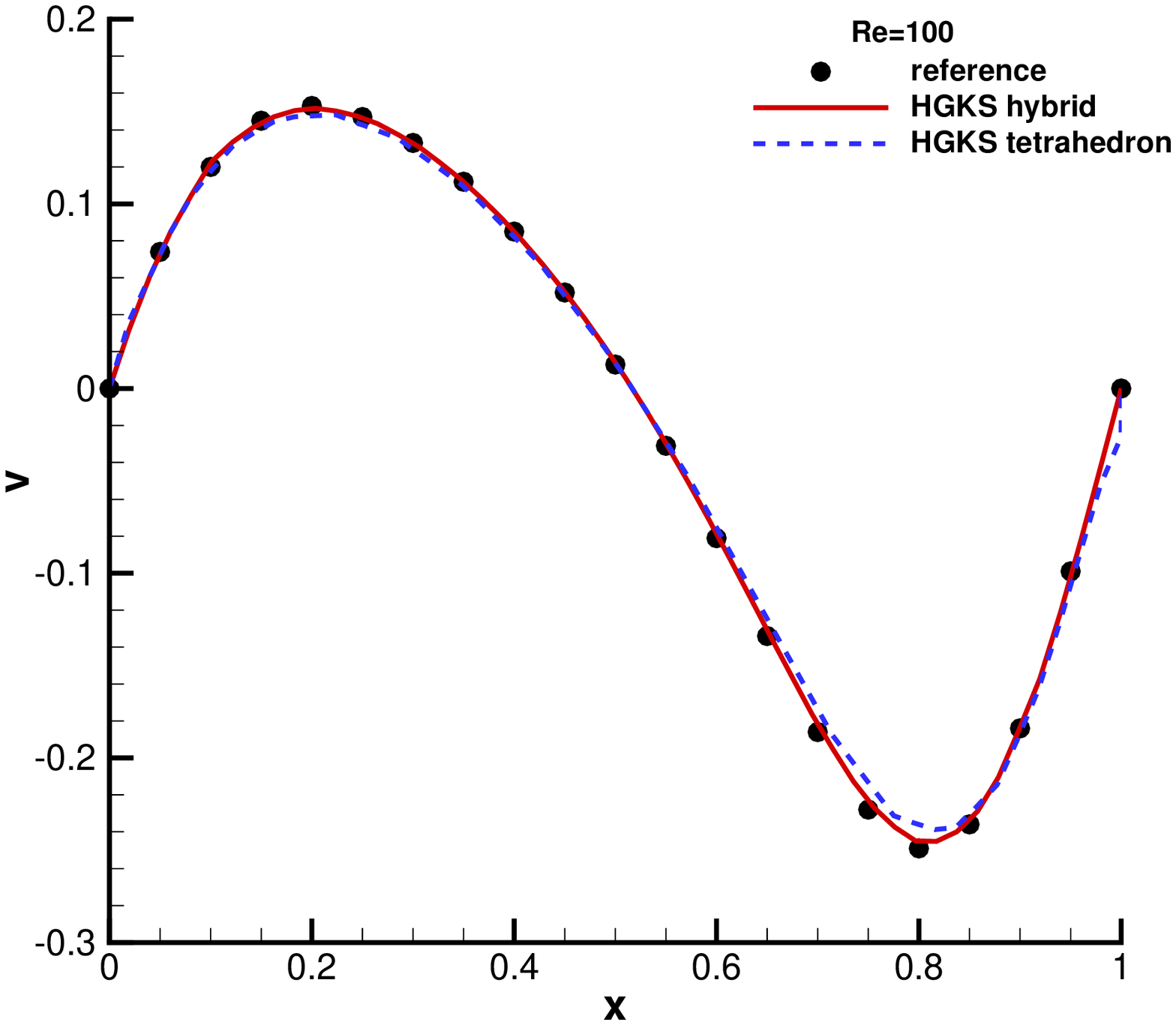}
\includegraphics[width=0.45\textwidth]{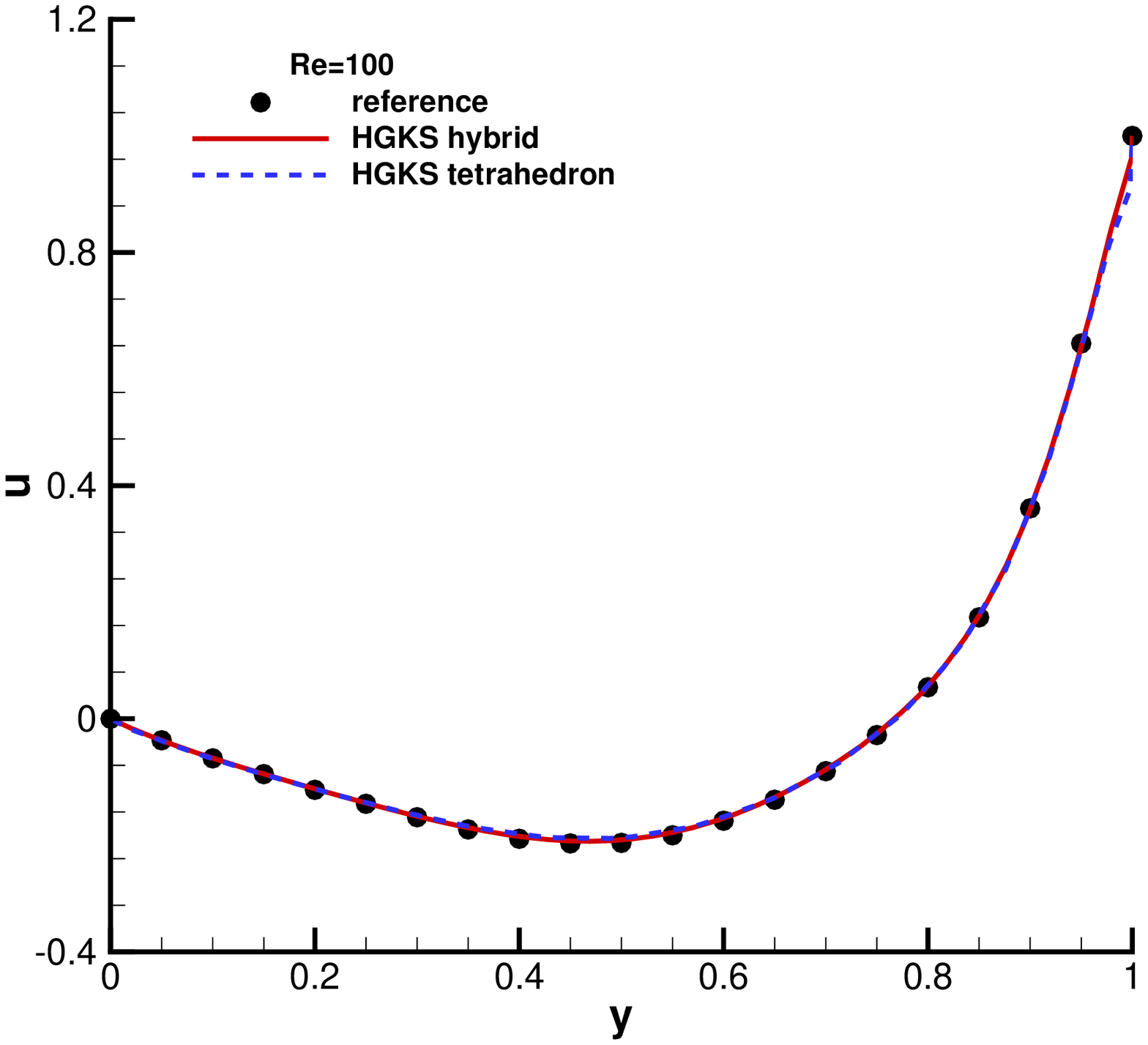}
\caption{\label{cavity-velocity} Lid-driven cavity flow: the steady
state U-velocity profiles along the vertical centerline, V-velocity
profiles along the horizontal centerline for $Re=1000,400$ and
$100$.}
\end{figure*}

\subsection{Lid-driven cavity flow}
The lid-driven cavity problem is one of the most important
benchmarks for numerical Navier-Stokes solvers. The fluid is bounded
by a unit cubic $[0, 1]\times[0, 1]\times[0, 1]$ and driven by a
uniform translation of the top boundary with $y=1$.  In this case,
the flow is simulated with Mach number $Ma=0.15$ and all the
boundaries are isothermal and nonslip. Numerical simulations are
conducted with the Reynolds numbers of $Re=1000$, $400$ and $100$.
This case is performed by both hybrid and tetrahedral meshes. The
tetrahedral mesh contains $6\times20^3$ cells,  in which every
cuboid is divided into six tetrahedron cells. The hybrid mesh
contains $1.6\times40^3$ cells, including $76800$ prisms and $25600$
hexahedrons. To improve the resolution, the mesh near the well is
refined and both meshes are shown in Fig.\ref{cavity-mesh}. The
three cases correspond to the convergent solutions, and the LU-SGS
method is used for the temporal discretization. The $U$-velocity
profiles along the vertical centerline line, $V$-velocity profiles
along the horizontal centerline in the symmetry $x-y$ plane and the
benchmark data for $Re=1000$ \cite{Case-Albensoeder}, $Re=400$ and
$Re=100$ \cite{Case-Shu} are shown in Fig.\ref{cavity-velocity}. The
numerical results agree well with the benchmark data. For this case,
a coarse mesh is used, especially for the tetrahedral meshes. The
agreement between them shows that current HGKS is capable of
simulating three-dimensional laminar flows.

\begin{figure}[!h]
\centering
\includegraphics[width=0.66\textwidth]{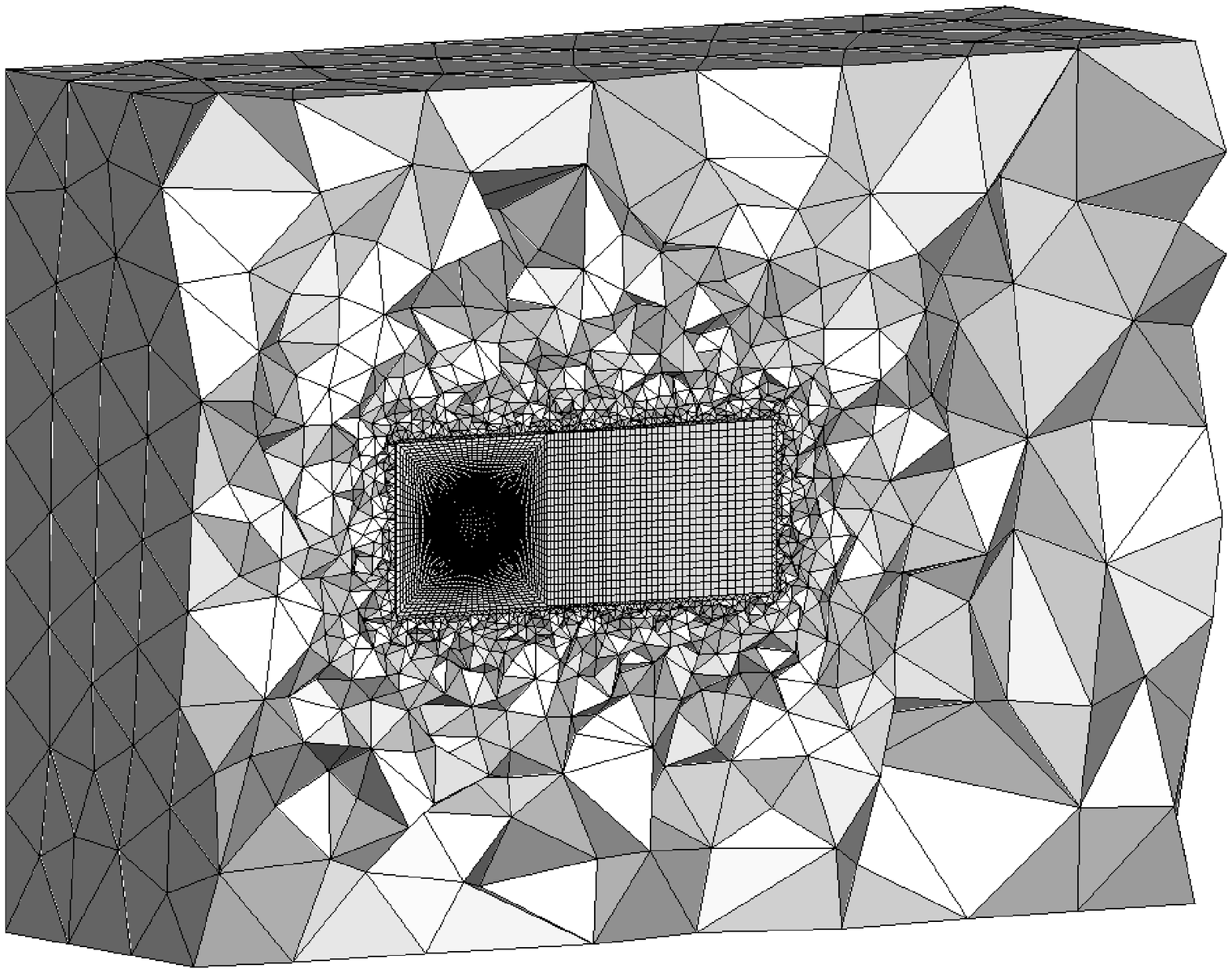}
\caption{\label{sphere-mesh} Flows passing through a sphere: the
mesh distribution.}
\includegraphics[width=0.66\textwidth]{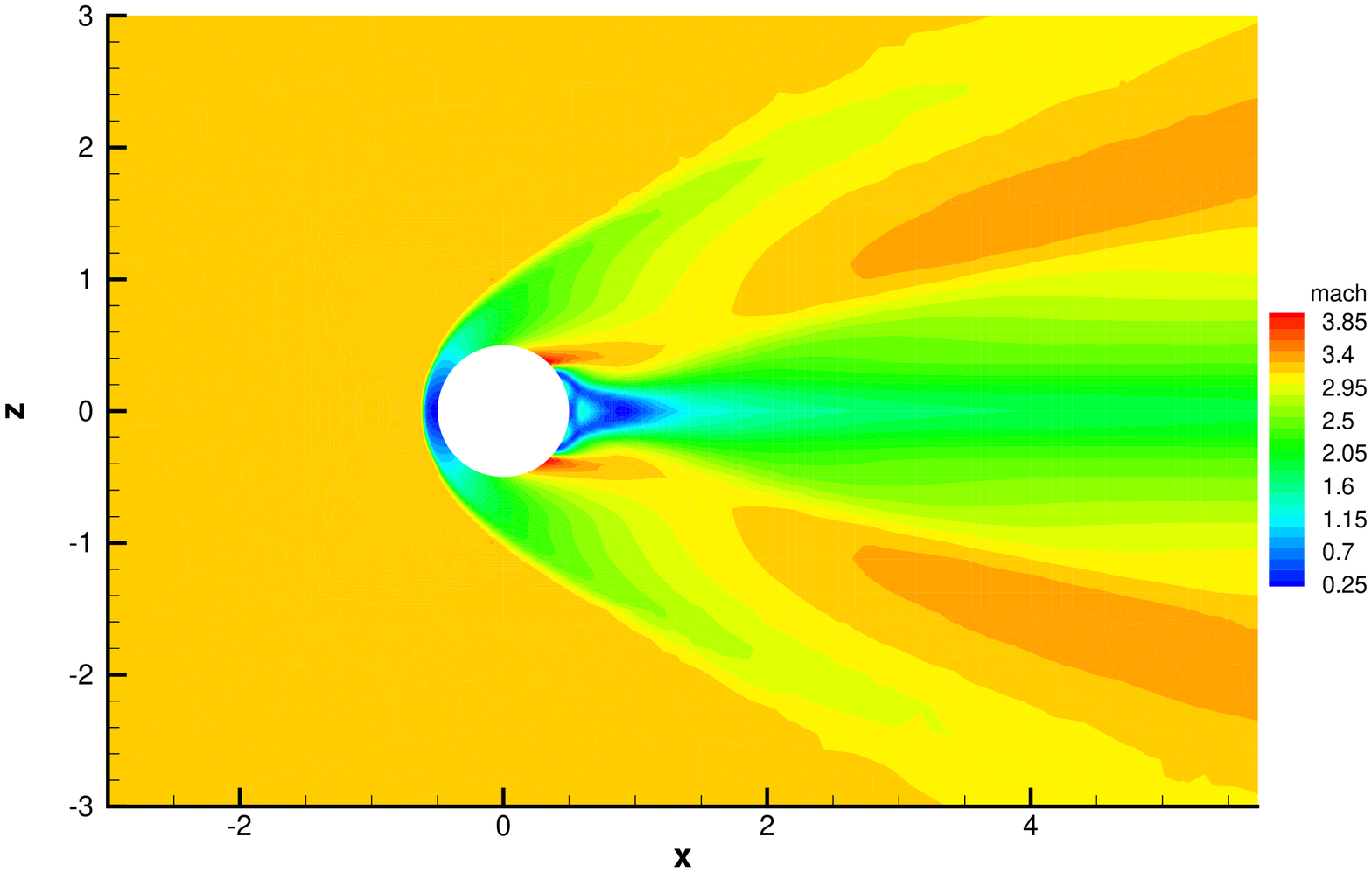}
\caption{\label{sphere-inviscid} Flows passing through a sphere: the
Mach number distribution  at vertical centerline planes with
$Ma_\infty=3.2$ for the inviscid flow.}
\end{figure}

\begin{figure}[!h]
\centering
\includegraphics[width=0.495\textwidth]{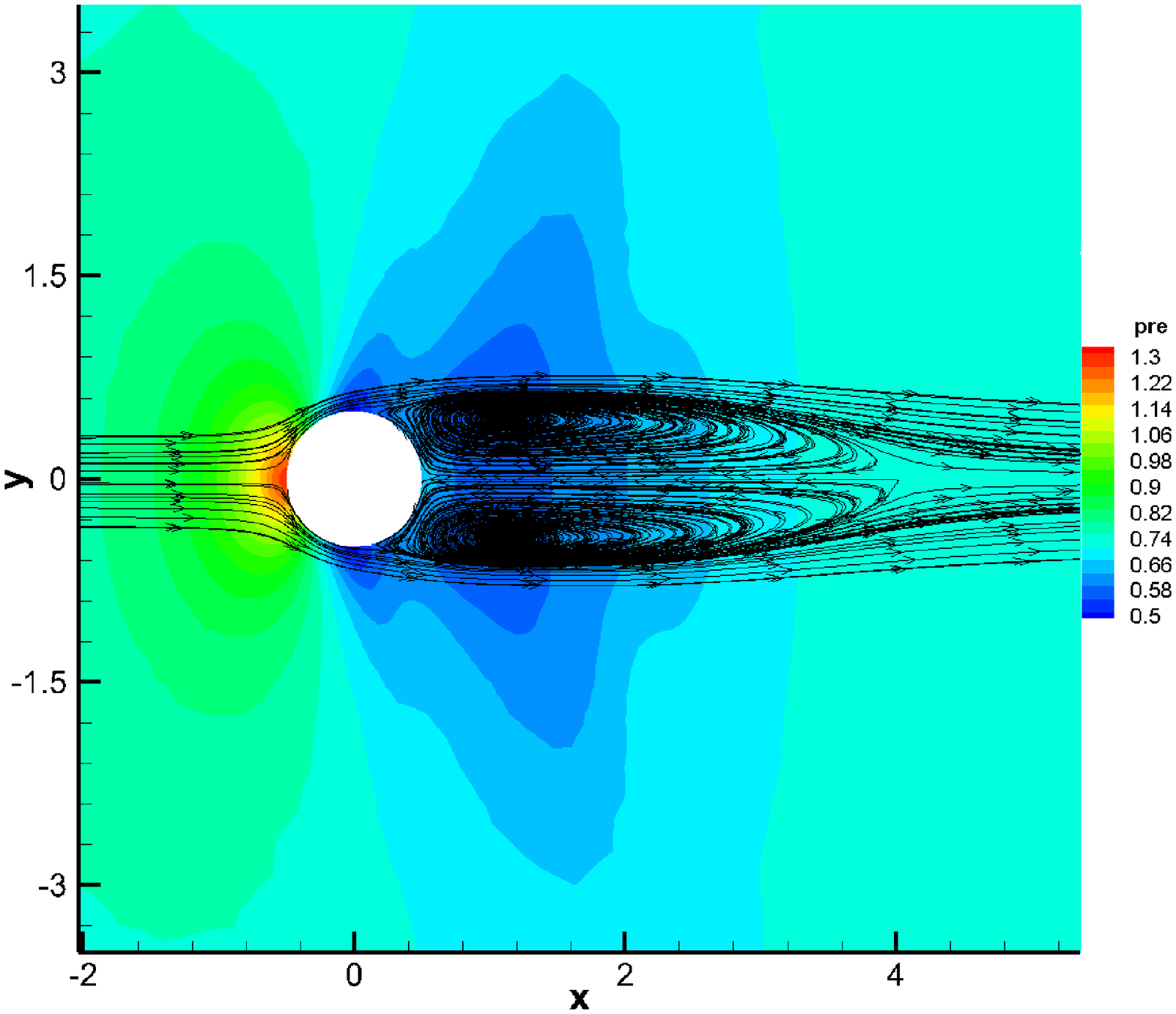}
\includegraphics[width=0.495\textwidth]{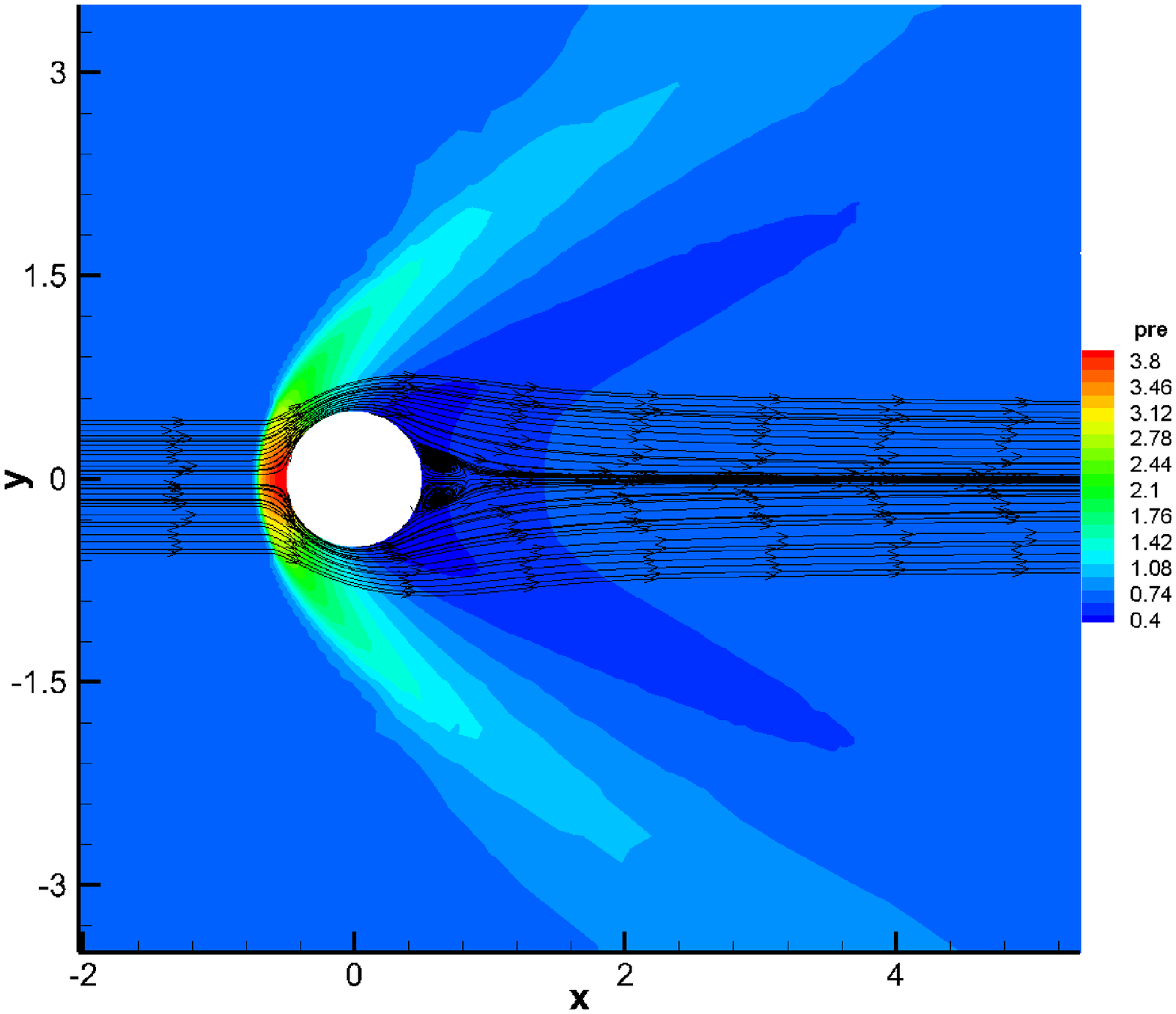}
\caption{\label{sphere-1} Flows passing through a sphere: the
pressure  and streamline  distributions at vertical centerline
planes with $Ma_\infty=0.95$ and $2.0$.}
\includegraphics[width=0.495\textwidth]{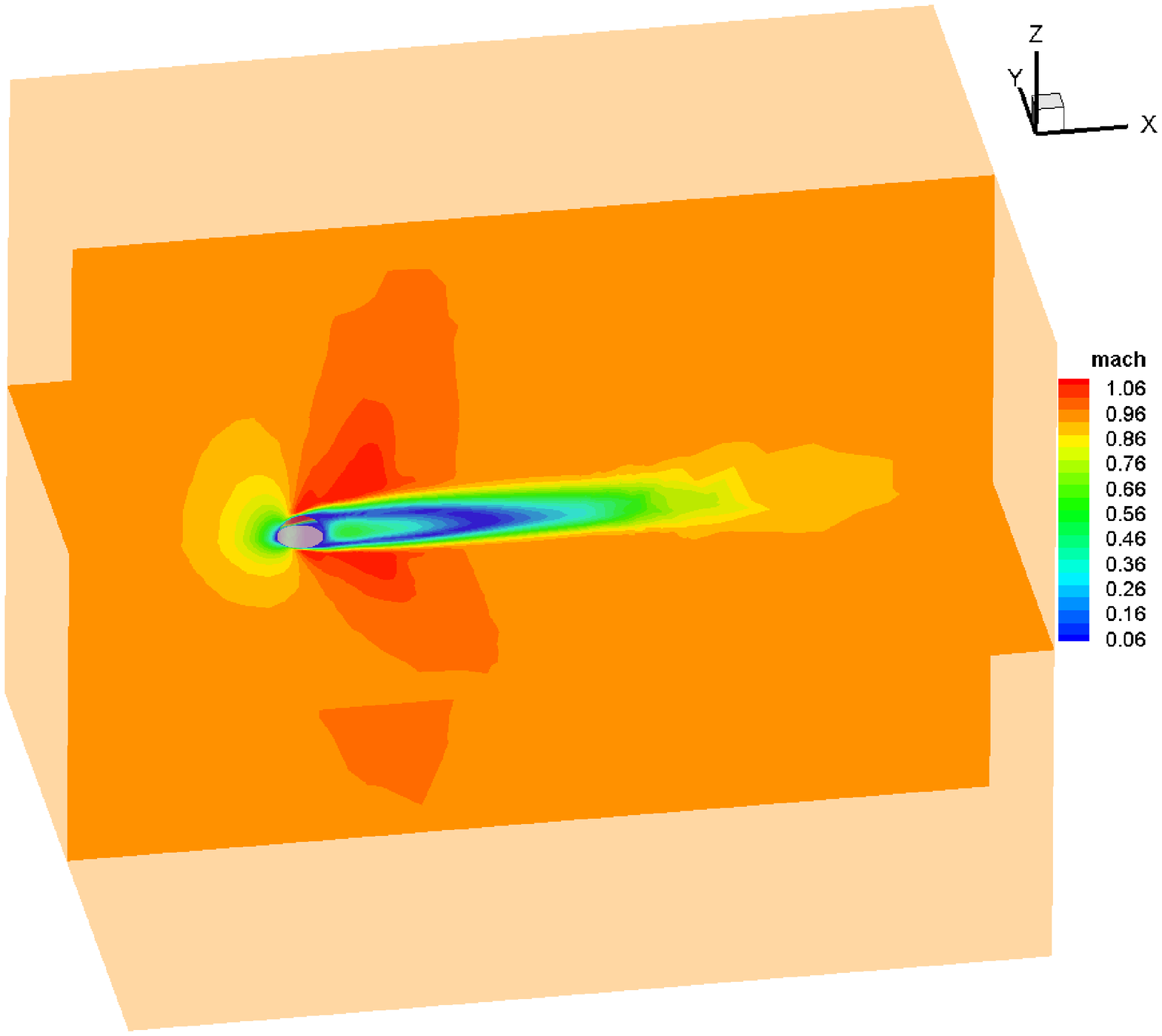}
\includegraphics[width=0.495\textwidth]{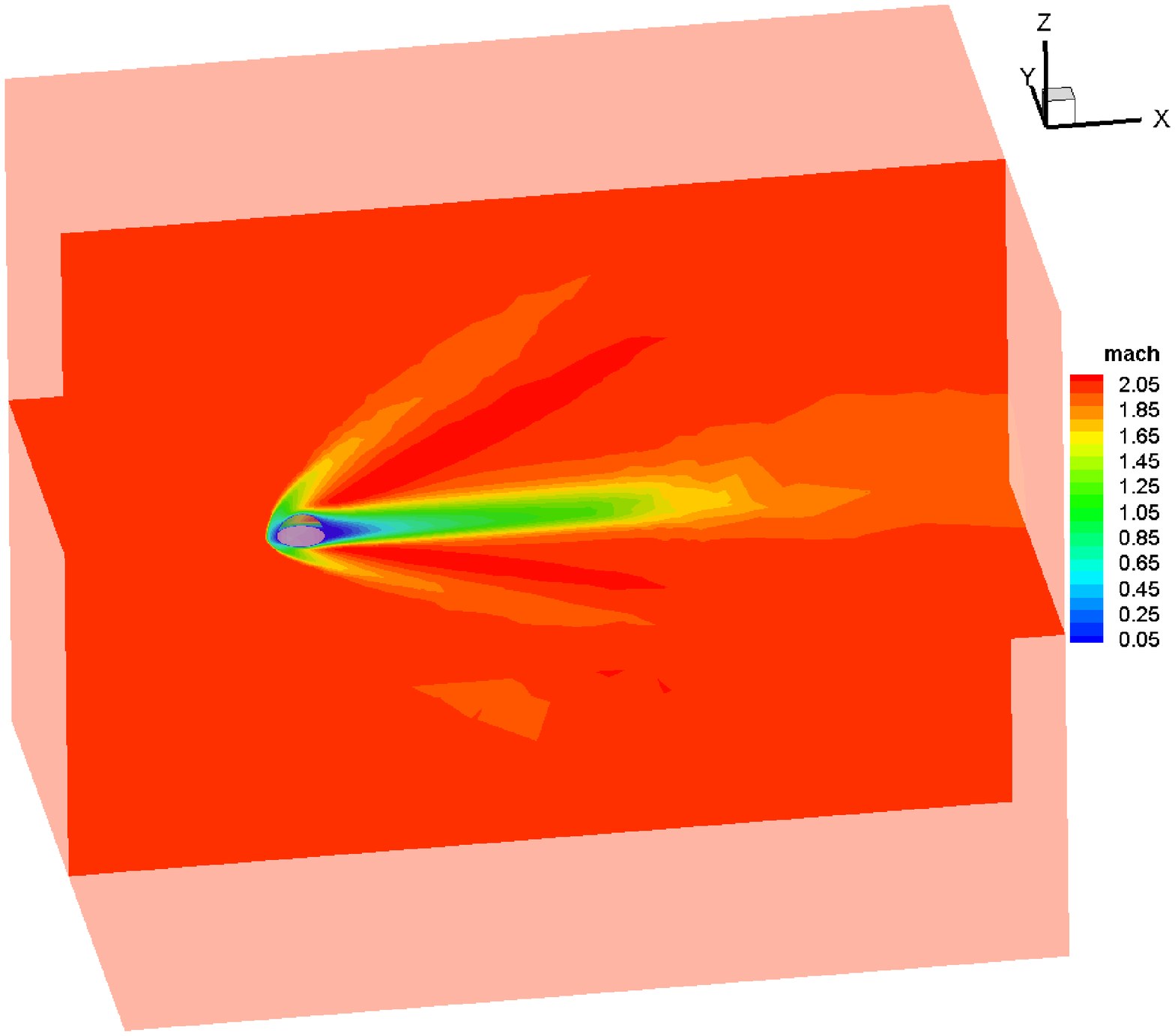}
\caption{\label{sphere-2} Flows passing through a sphere: the Mach
number distributions at vertical centerline planes with
$Ma_\infty=0.95$ and $2.0$.}
\end{figure}

\begin{table}[!h]
\begin{center}
\def\temptablewidth{0.55\textwidth}
{\rule{\temptablewidth}{1.0pt}}
\begin{tabular*}{\temptablewidth}{@{\extracolsep{\fill}}c|cc}
scheme                                    & $\phi$ &   $L$     \\
\hline
WENO6 \cite{Case-Nagata}                     &  111.5 & 3.48   \\
HGKS-compact scheme  \cite{GKS-high-2}    &  112.7 & 3.30\\
Current scheme                                              &  107.0       &   3.49
\end{tabular*}
{\rule{\temptablewidth}{1.0pt}}
\end{center}
\caption{\label{sphere-tran} Flows passing a sphere: quantitative
comparisons of separation angle $\phi$ and closed wake length $L$
for $Ma_\infty=0.95$ and $Re=300$.}
\begin{center}
\def\temptablewidth{0.75\textwidth}
{\rule{\temptablewidth}{1.0pt}}
\begin{tabular*}{\temptablewidth}{@{\extracolsep{\fill}}c|ccc}
scheme                                    & $\phi$ &   $L$   &Shock stand-off  \\
\hline
WENO6 \cite{Case-Nagata}                  &  150.9      &   0.38  & 0.21 \\
HGKS-compact scheme  \cite{GKS-high-2}    &  148.5      &   0.45  & 0.28-0.31\\
Current scheme                            &  149.8      &   0.34  &  0.22\\
\end{tabular*}
{\rule{\temptablewidth}{1.0pt}}
\end{center}
\caption{\label{sphere-super} Flows passing a sphere: quantitative
comparisons  of separation angle $\phi$, closed wake length $L$ and
shock stand-off for $Ma_\infty=2$ and $Re=300$.}
\end{table}

\subsection{Flows passing a sphere}
This case is used to test the capability in resolving the low-speed
to hypersonic flows, and the initial condition is given as a free
stream condition
\begin{align*}
(\rho,U,V,W,p)_{\infty} = (1, Ma_\infty,0,0, 1/\gamma),
\end{align*}
where $\gamma= 1.4$. The computation domain is  $[-5,15] \times
[-7.5, 7.5] \times [-7.5,7.5]$. As shown in Fig.\ref{sphere-mesh},
the hybrid mesh with $462673$ cells is used. The inlet and outlet
boundary conditions are given according to Riemann invariants, the
slip adiabatic boundary condition is used for inviscid flows and the
non-slip adiabatic boundary condition is imposed for viscous flows
on the surface of sphere.

\begin{figure*}[!h]
\centering
\includegraphics[width=0.45\textwidth]{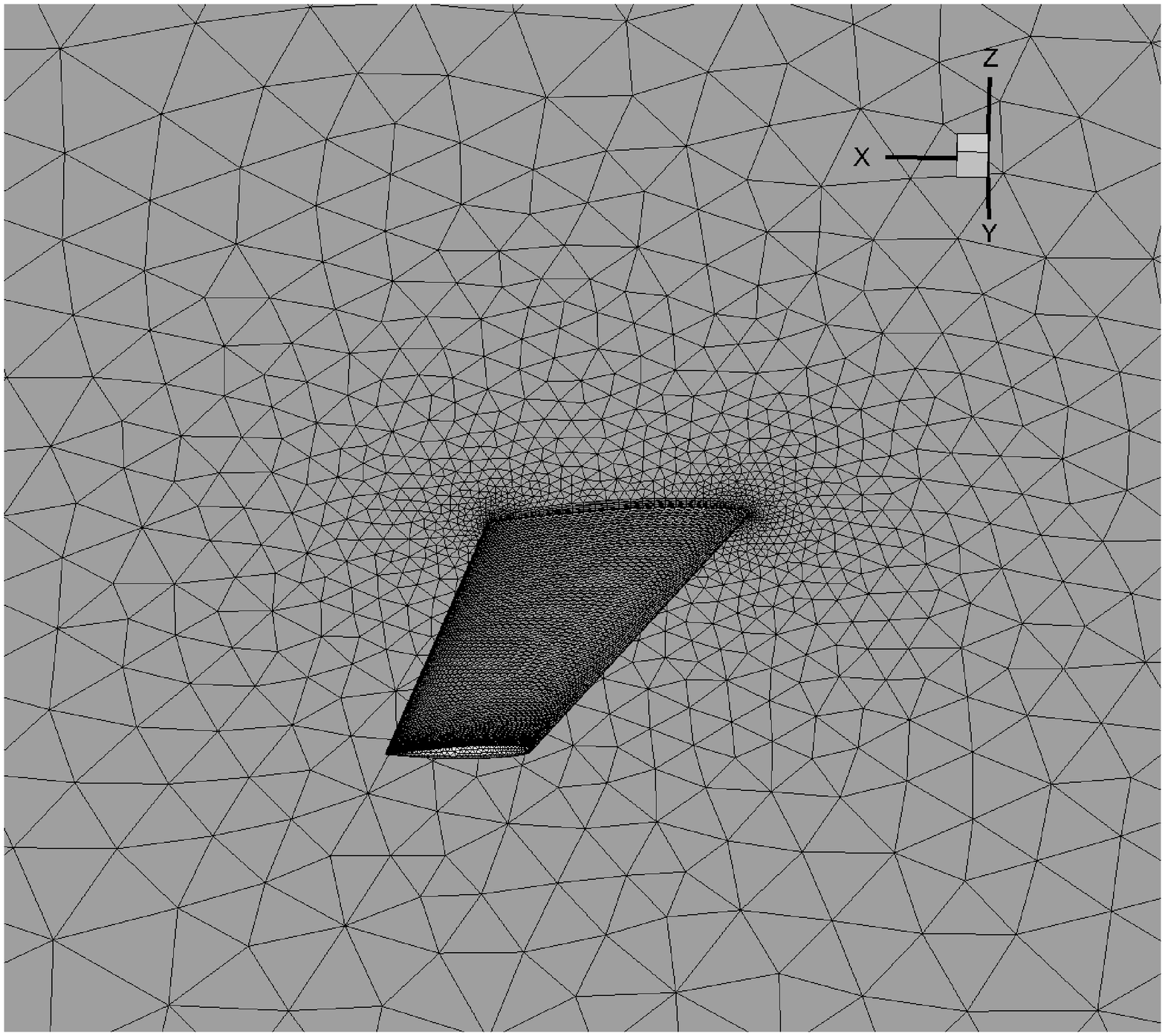}
\includegraphics[width=0.45\textwidth]{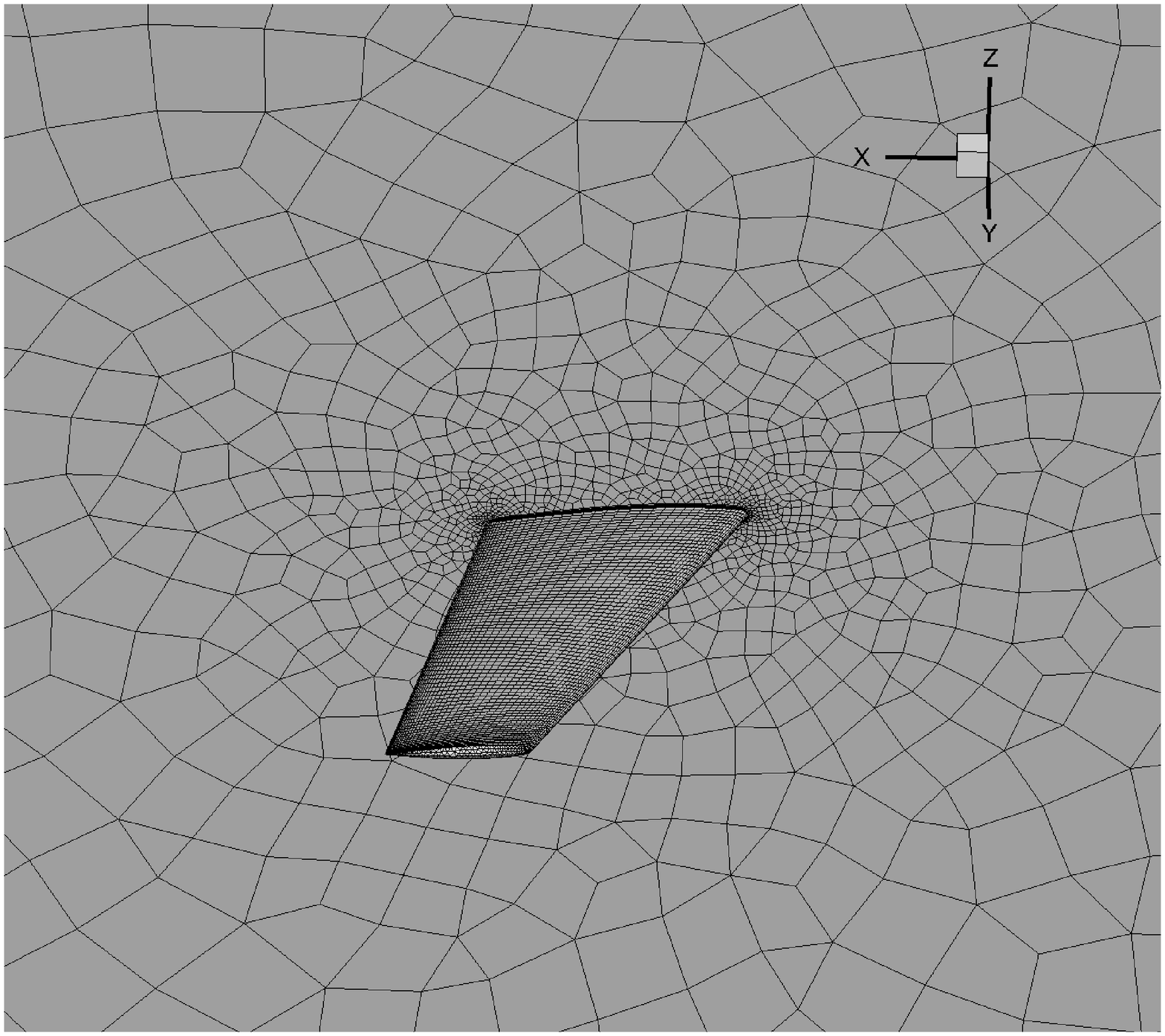}\\
\includegraphics[width=0.45\textwidth]{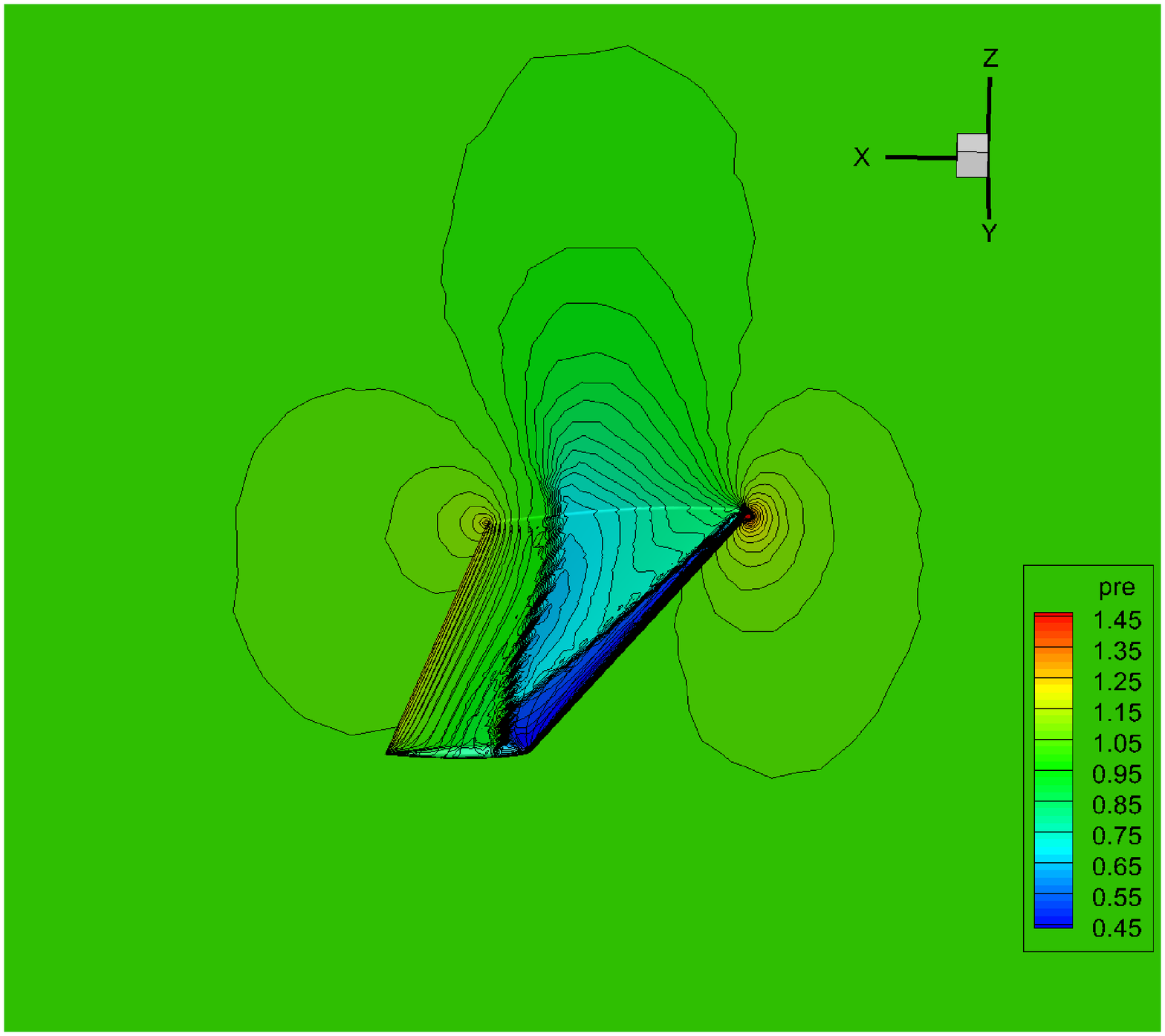}
\includegraphics[width=0.45\textwidth]{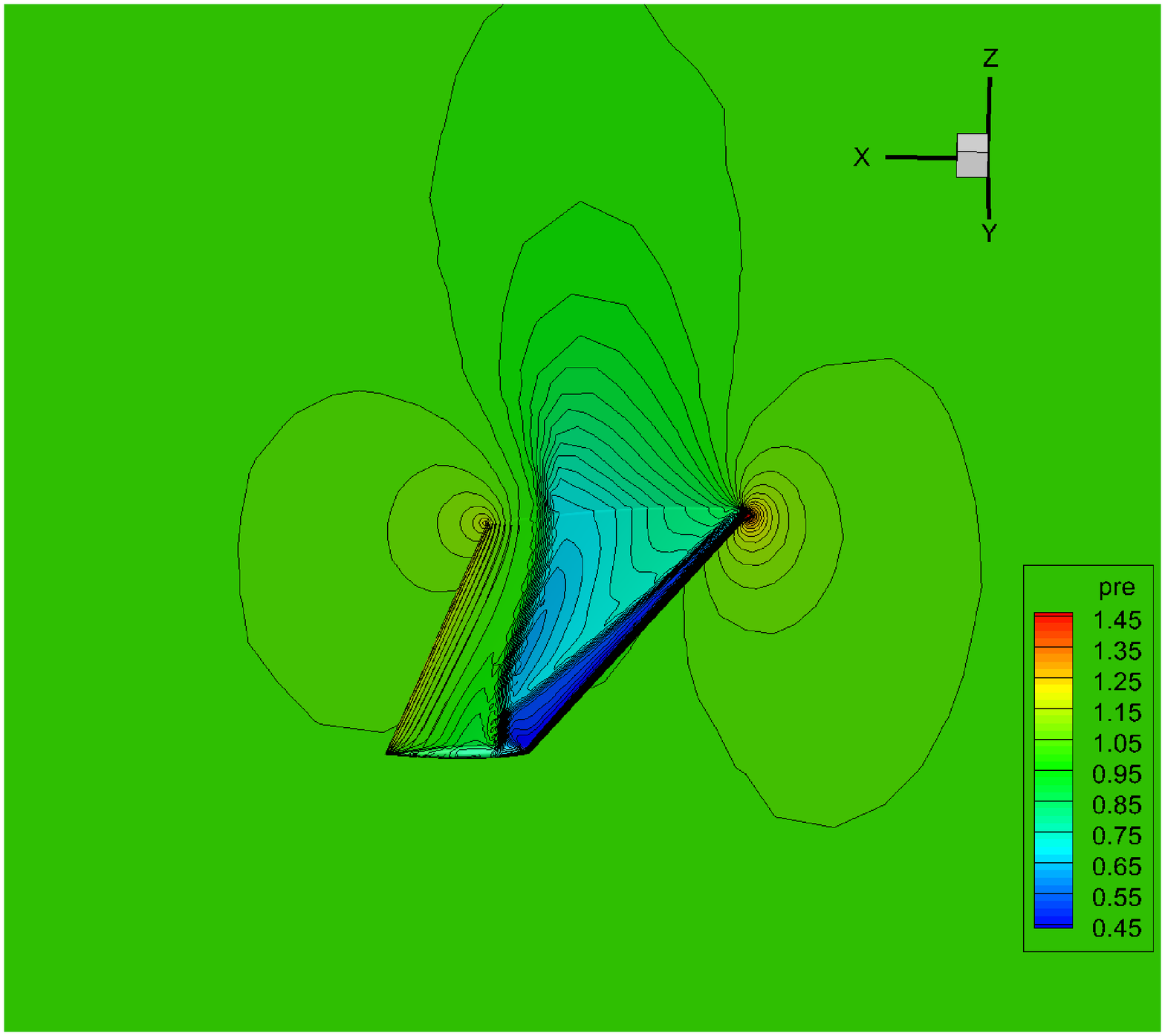}\\
\includegraphics[width=0.45\textwidth]{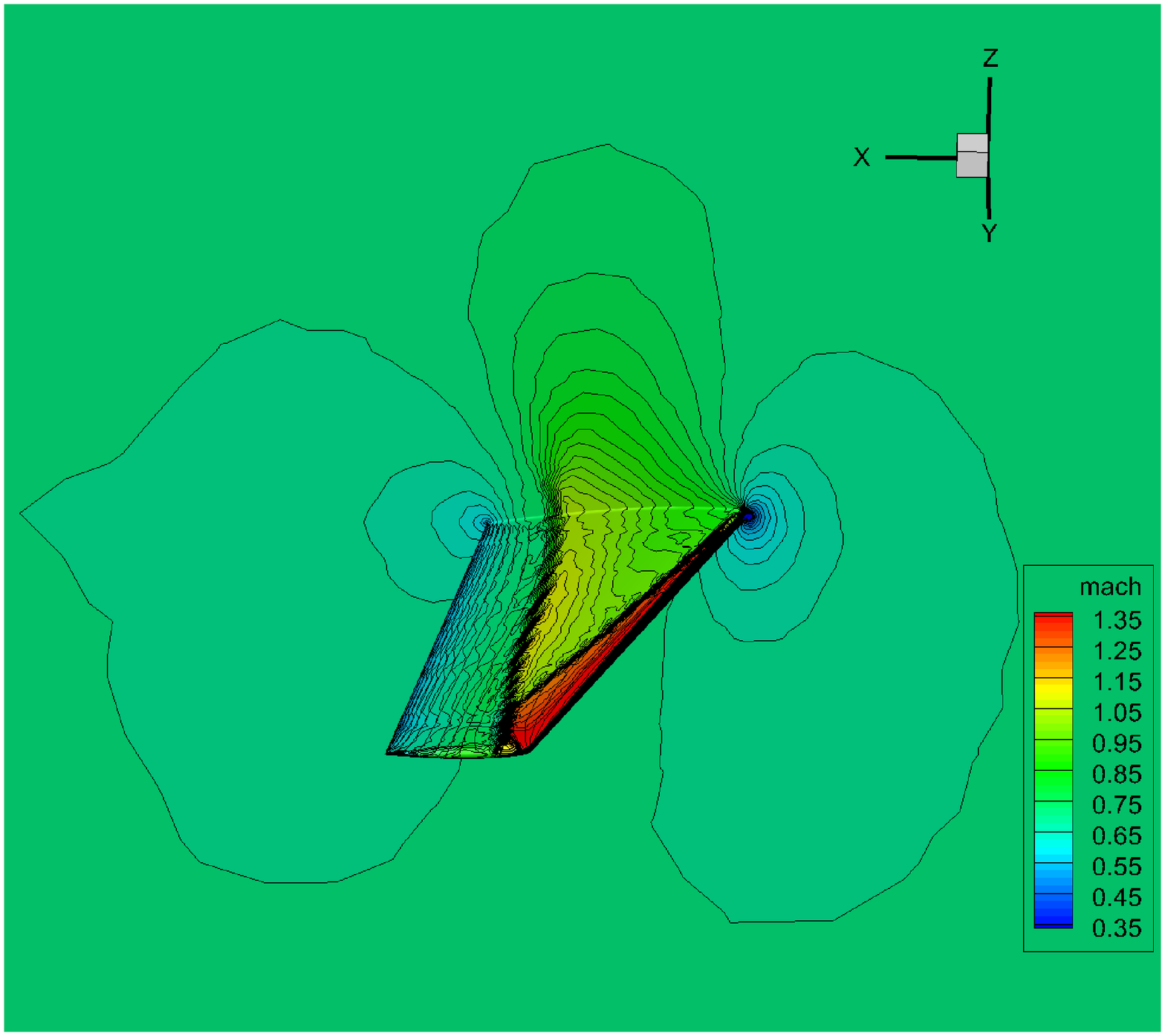}
\includegraphics[width=0.45\textwidth]{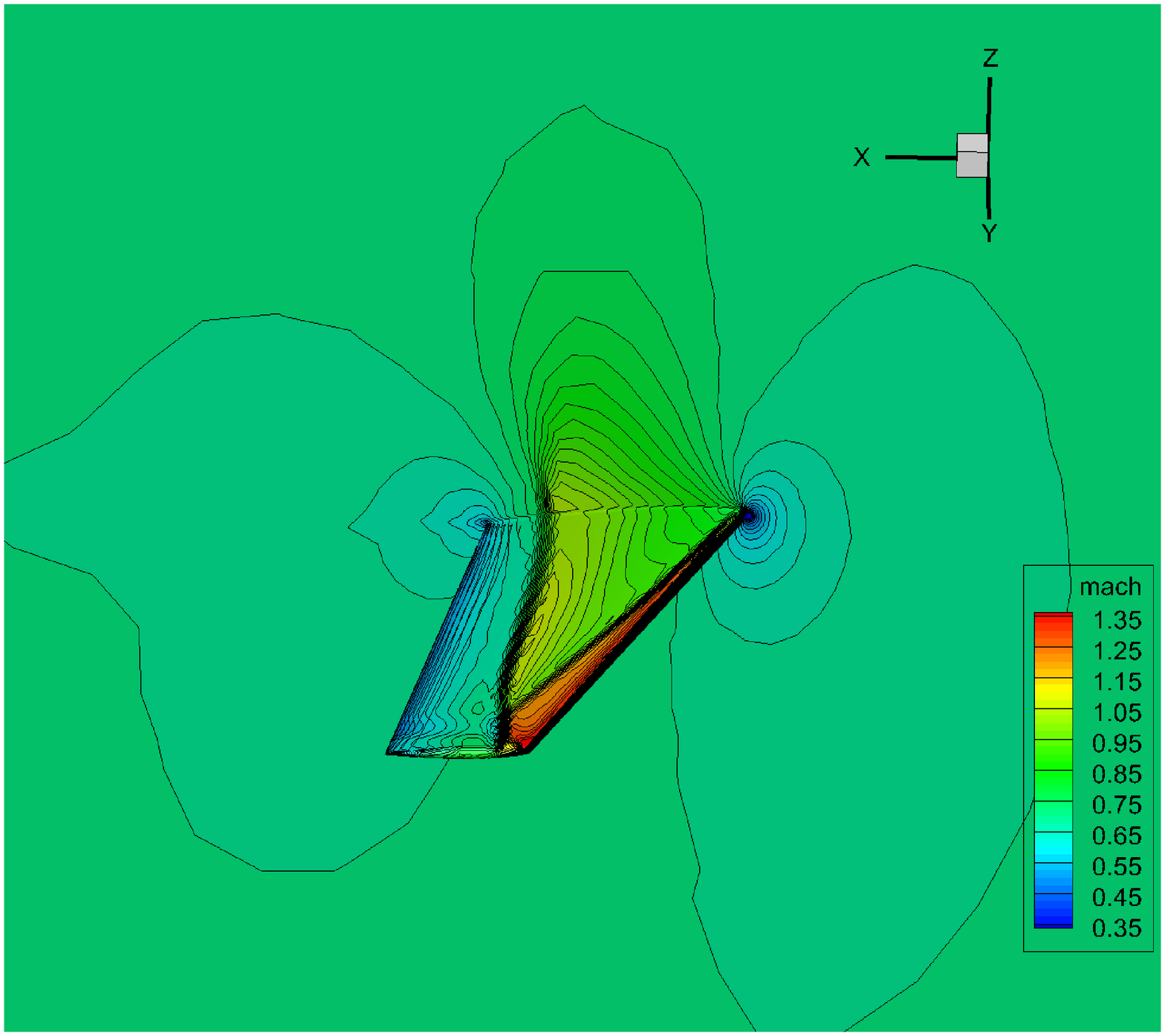}
\caption{\label{M6-wing-1} ONERA M6 wing: the local mesh, pressure
and Mach number distribution  for the inviscid flow with tetrahedral
and hybrid meshes.}
\end{figure*}

\begin{figure*}[!h]
\centering
\includegraphics[width=0.45\textwidth]{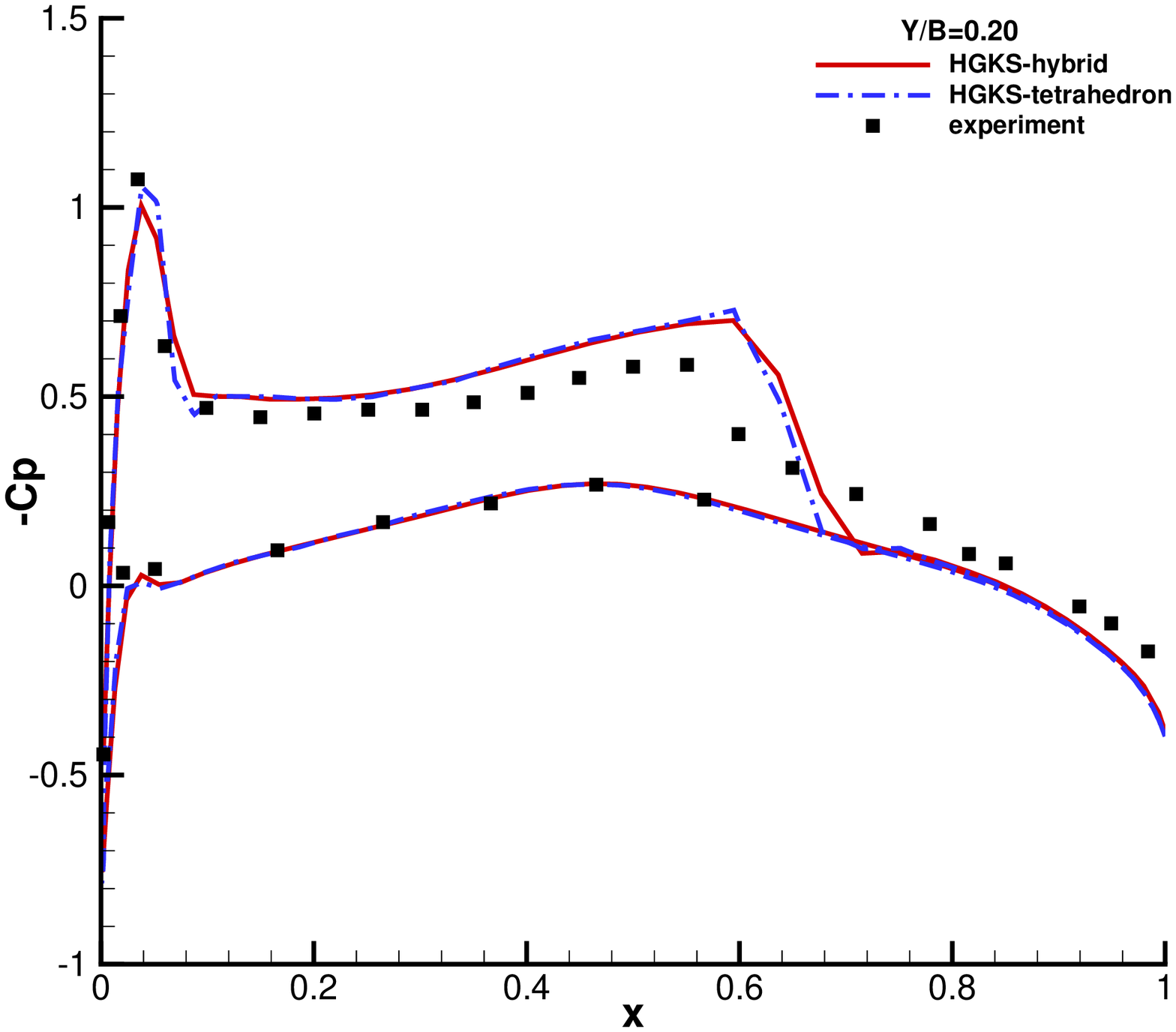}
\includegraphics[width=0.45\textwidth]{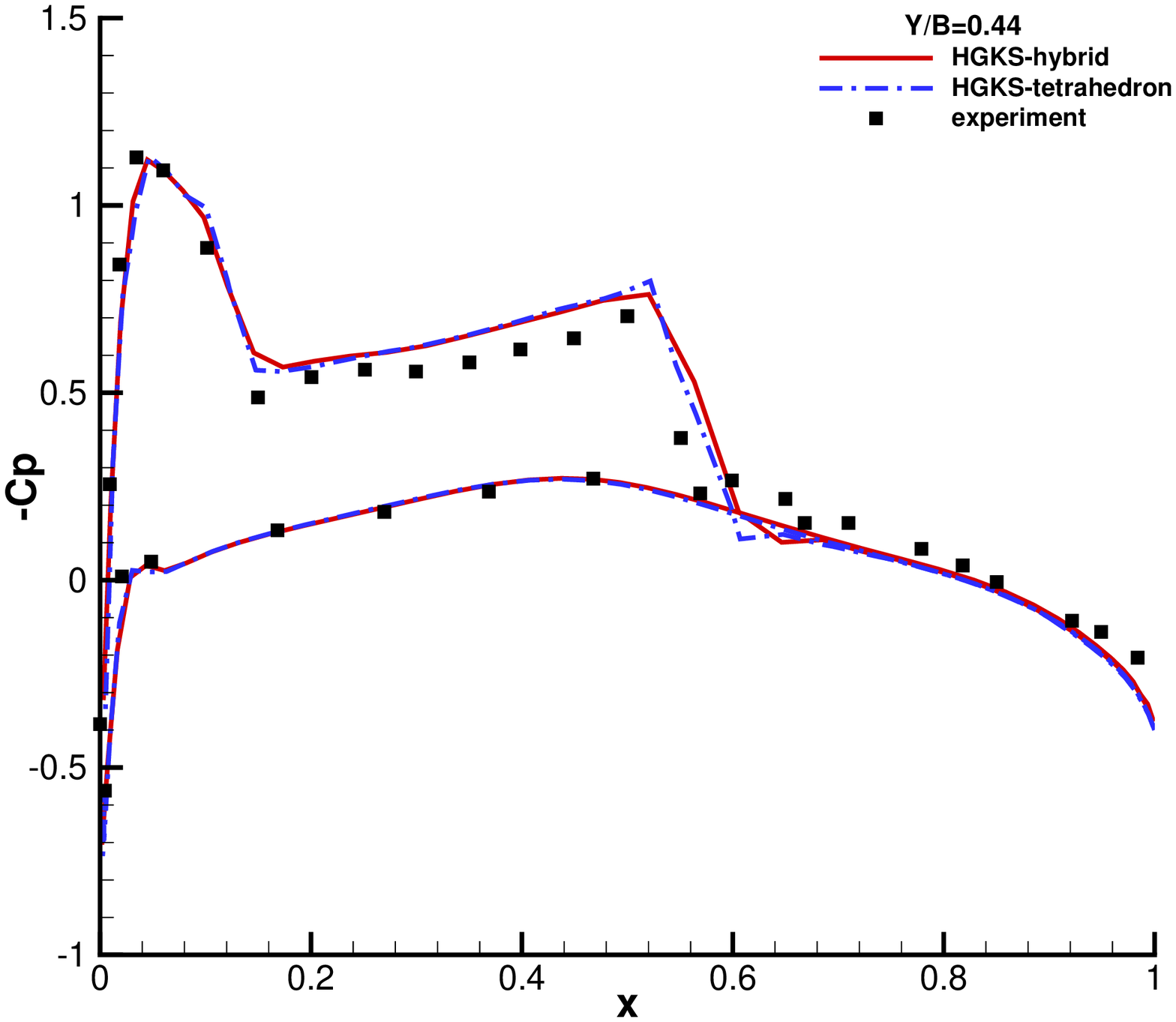}\\
\includegraphics[width=0.45\textwidth]{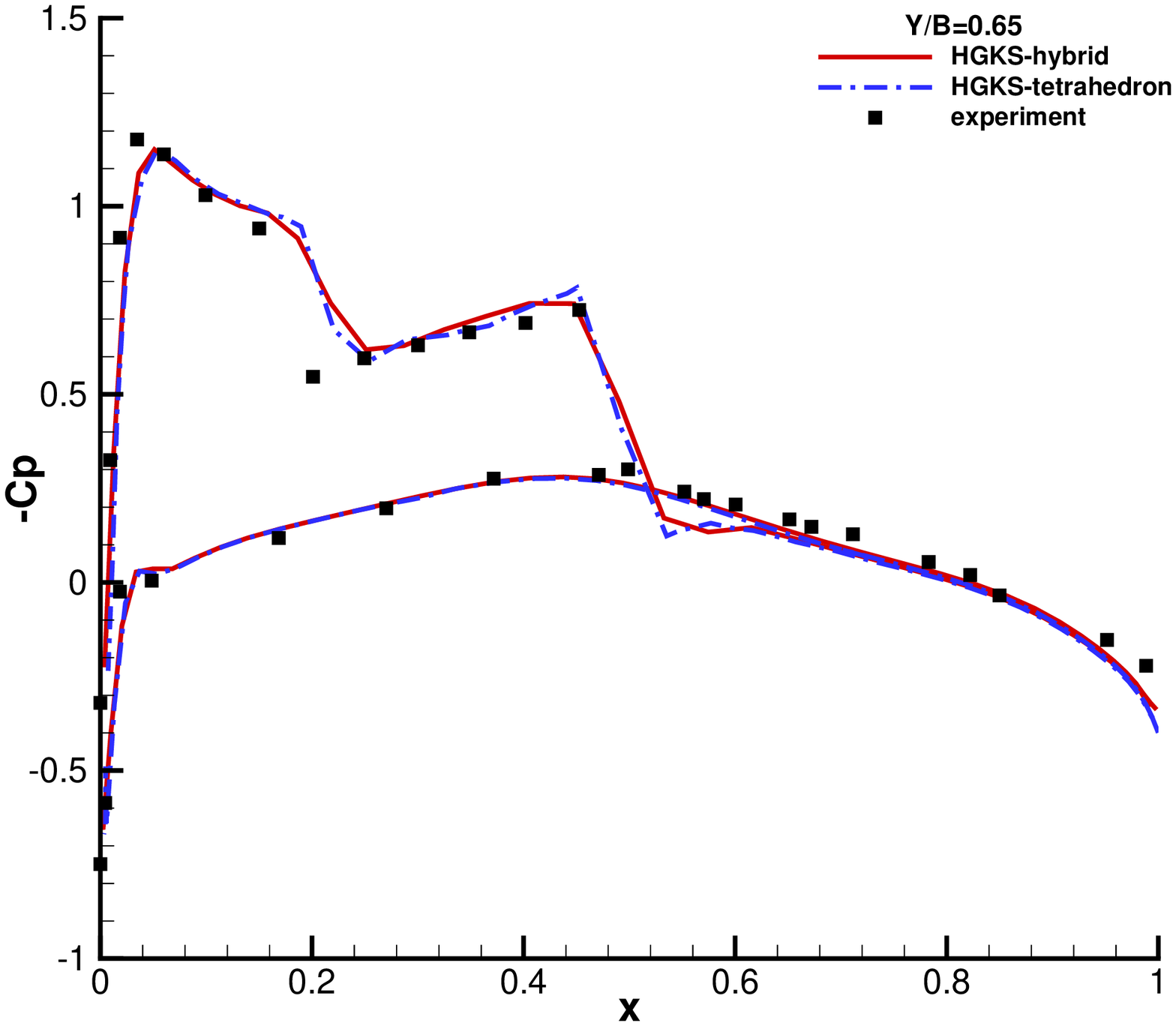}
\includegraphics[width=0.45\textwidth]{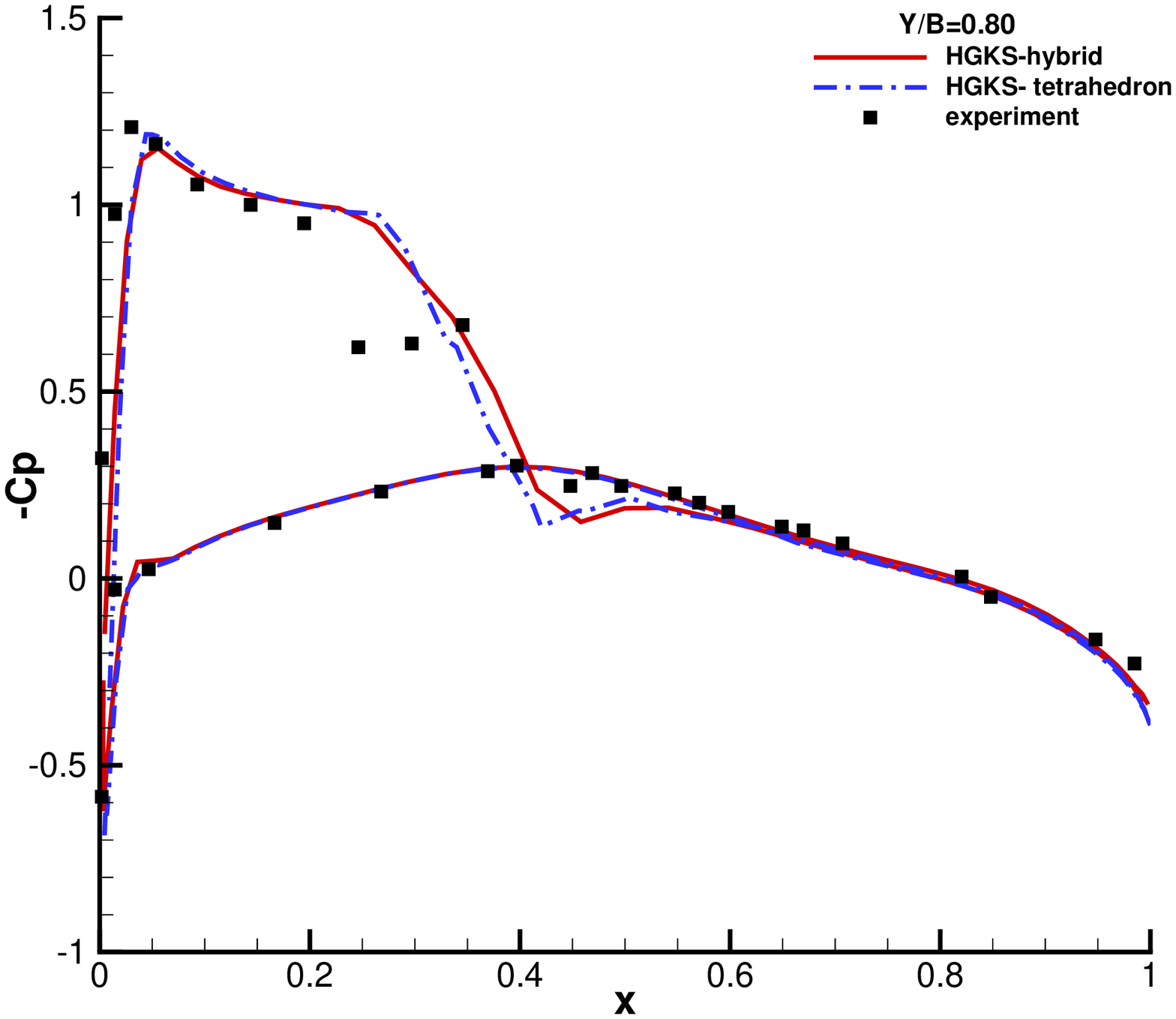}\\
\includegraphics[width=0.45\textwidth]{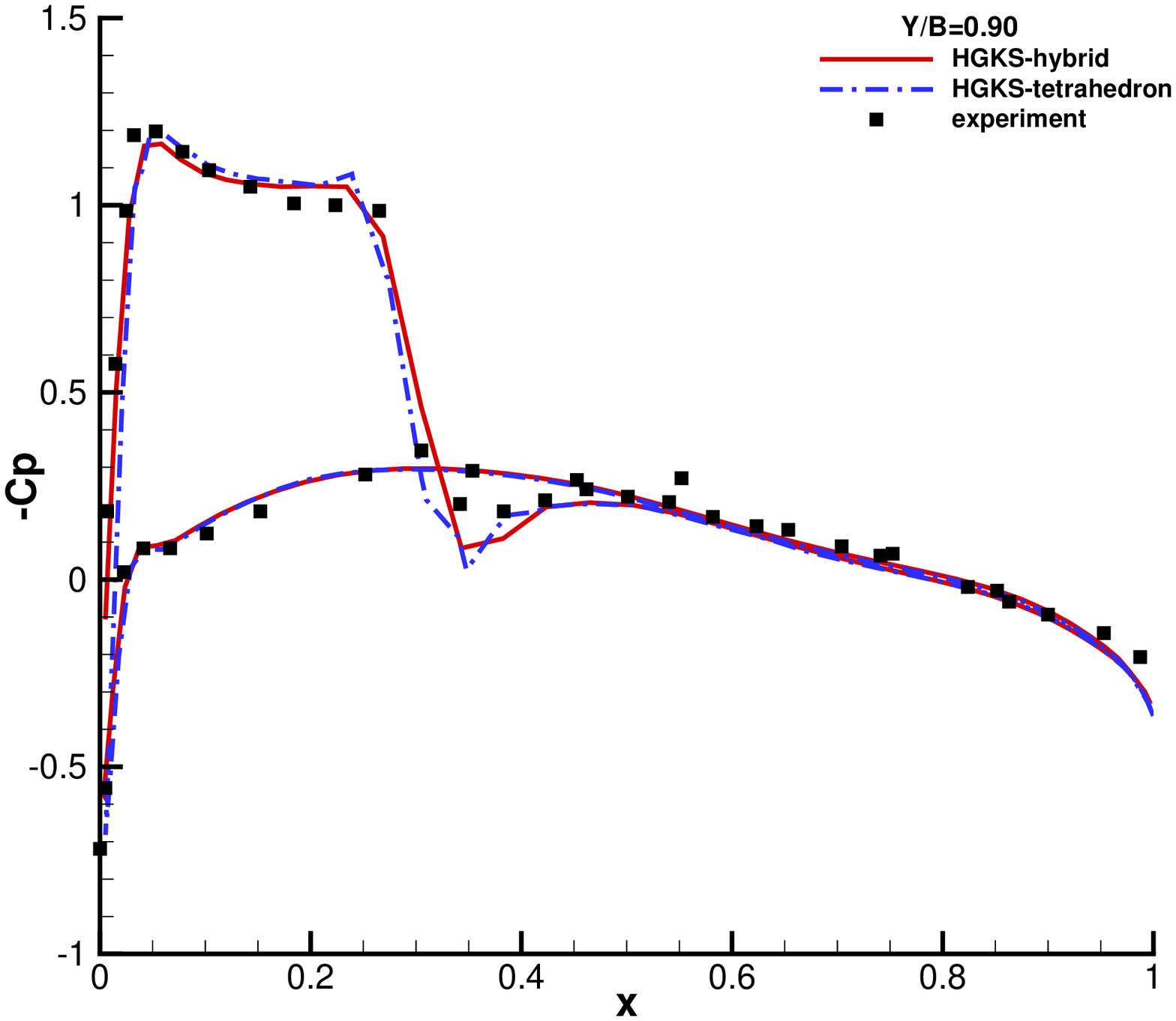}
\includegraphics[width=0.45\textwidth]{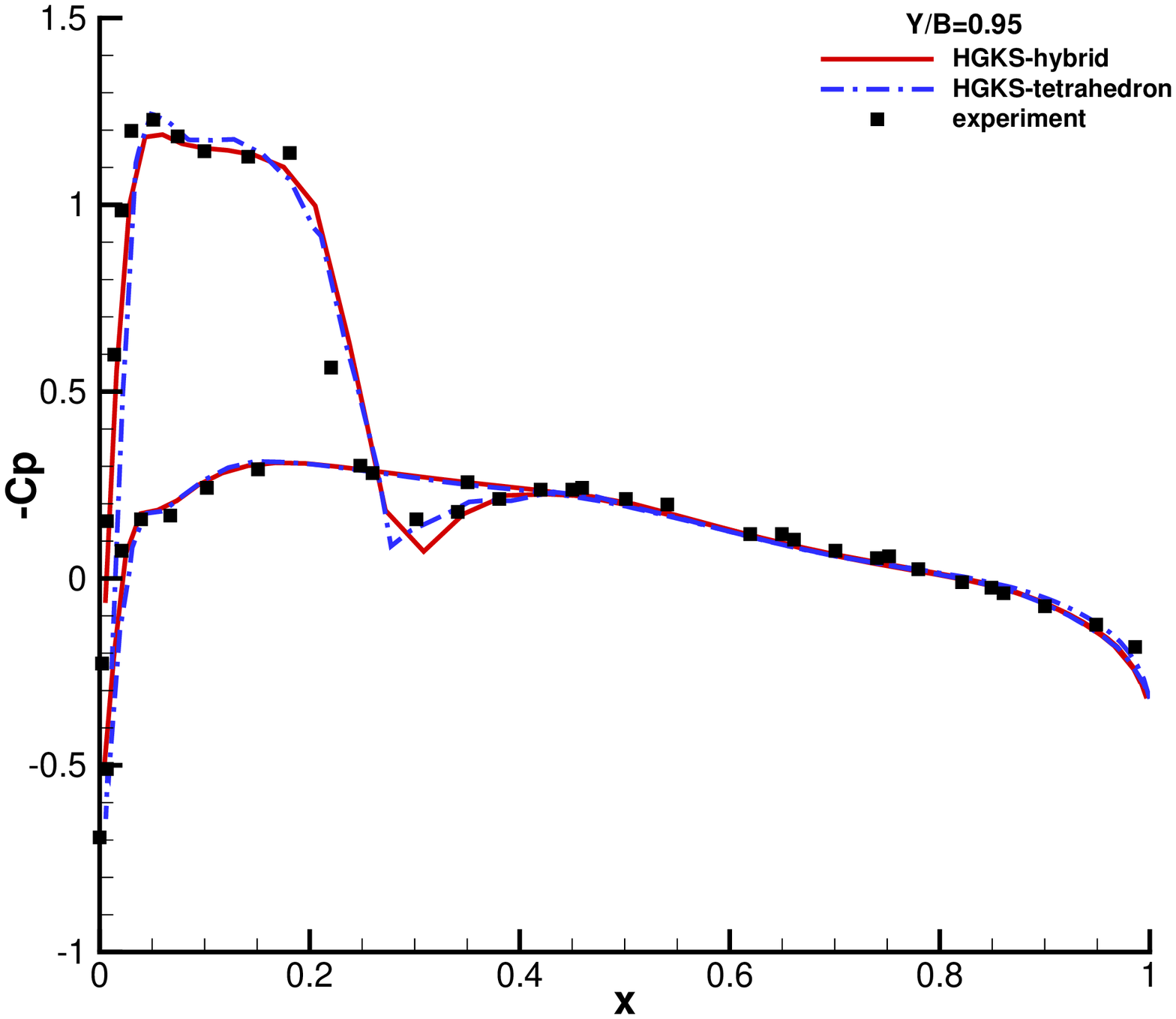}
\caption{\label{M6-wing-2} ONERA M6 wing: the pressure coefficient
distributions at $Y/B = 0.20$, $0.44$, $0.65$, $0.80$, $0.90$ and
$0.95$ for the inviscid flow with tetrahedral and hybrid meshes.}
\end{figure*}

For the hyperbolic inviscid flows, this case is tested with the
explicit and implicit schemes. For the explicit scheme, the Mach
number of inflow can be tested only up to $Ma_\infty=2.2$, and the
codes blow up due to the vacuum state forms at the leeward side of
the sphere.  The implicit scheme is more robust than the explicit
scheme, and the Mach number of inflow can be tested up to
$Ma_\infty=3.2$. The Mach number distribution at vertical centerline
planes is given in Fig.\ref{sphere-inviscid}, where the largest Mach
number is 3.97 and no special treatment is needed for
reconstruction. For the viscous flows,  two cases with $Re=300$ are
tested with unstructured hybrid mesh, including supersonic flow with
$Ma_\infty=2.0$ and transsonic flow with $Ma_\infty=0.95$. For these
two cases, the computation will converge to steady states,  and the
LU-SGS method is used for the temporal discretization. For the
hypersonic flow, the dynamic viscosity is given by
\begin{align*}
\mu=\mu_\infty(\frac{T}{T_\infty})^{0.7},
\end{align*}
where $T_\infty$ and $\mu_\infty$ are free stream temperature and
viscosity. The pressure, Mach number and streamline distributions
are presented in Fig.\ref{sphere-1} and Fig.\ref{sphere-2} for the
two cases, and the robustness of the current scheme is validated.
The quantitative results of separation angle $\phi$ and closed wake
length $L$ are given in Table.\ref{sphere-tran} and
Table.\ref{sphere-super}. For the supersonic case, the position of
shock stand-off is also given in Table.\ref{sphere-super}.
Quantitative results agree well with the benchmark solutions
\cite{Case-Nagata}, and the slight deviation of compact gas-kinetic
scheme \cite{GKS-high-2} might caused by the coarser mesh.

\subsection{Transonic flow around ONERA M6 wing}
The transonic flow around the ONERA M6 wing is a standard benchmark
for engineering simulations. Besides the three-dimensional geometry,
the flow structures are complex including the interaction of shock
and turbulent boundary. Thus, it is a good candidate to test the
performance of the extended BGK model and high-order gas-kinetic
scheme. The inviscid flow around the wing is first tested, which
corresponds to a rough prediction of the flow field under a very
high Reynolds number. The incoming Mach number and angle of attack
are given by
\begin{align*}
Ma_{\infty} = 0.8395, ~AoA = 3.06^{\circ}.
\end{align*}
This case is performed by both hybrid and tetrahedral meshes, which
are given in Fig.\ref{M6-wing-1}. The tetrahedral includes 294216
cells and the hybrid mesh includes 201663 cells. The subsonic inflow
and outflow boundaries are all set according to the local Riemann
invariants, and the adiabatic and slip wall condition is imposed on
the solid wall. The local pressure and Mach number distributions are
also shown in Fig.\ref{M6-wing-1}, and the $\lambda$ shock is well
resolved by the current schemes. The comparisons on the pressure
distributions at the semi-span locations $Y/B = 0.20$, $0.44$,
$0.65$, $0.80$, $0.90$ and $0.95$ of the wing are given in
Fig.\ref{M6-wing-2}. The numerical results quantitatively agree well
with the experimental data \cite{Case-Schmitt}.

\begin{figure*}[!h]
\centering
\includegraphics[width=0.6\textwidth]{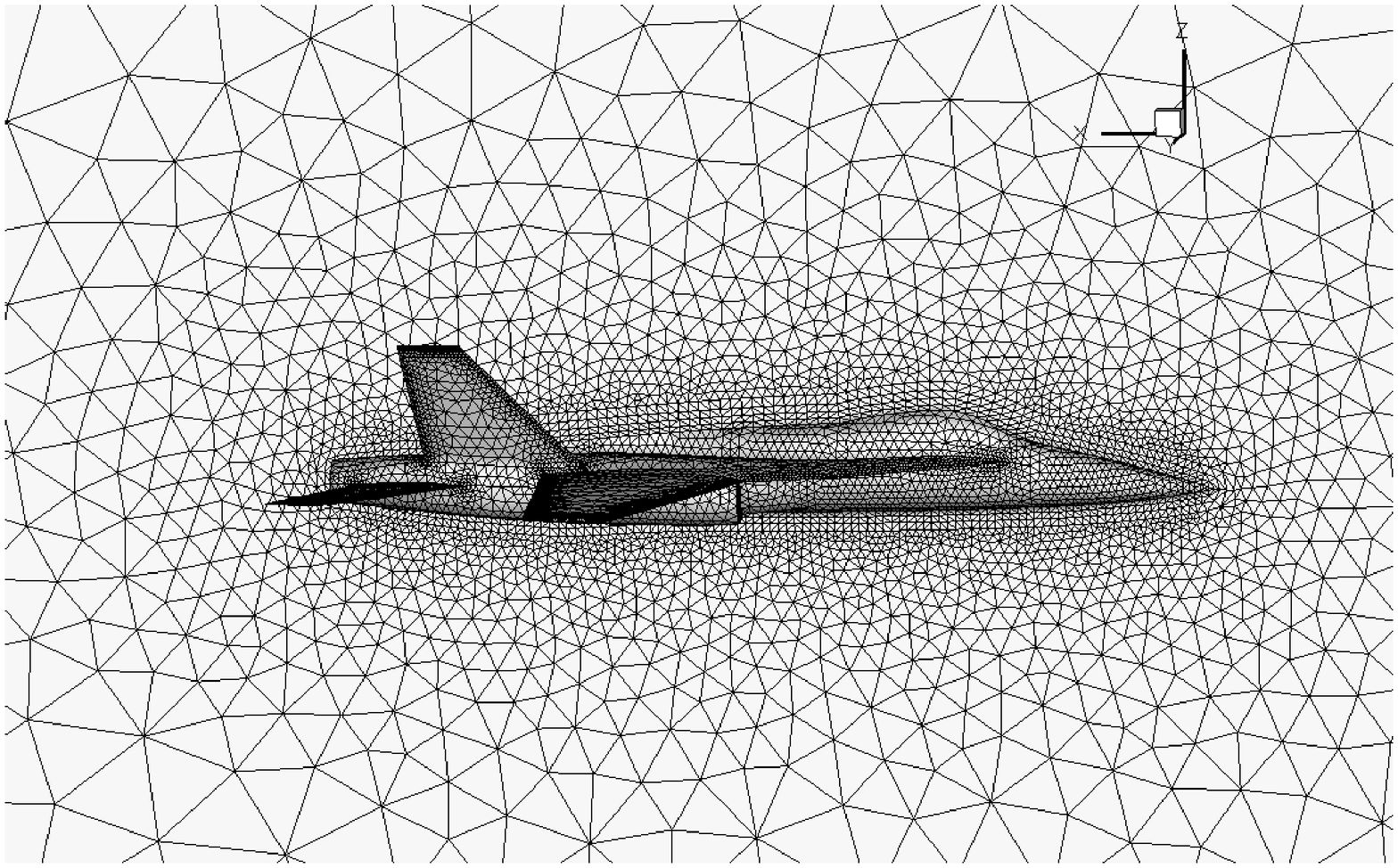}
\caption{\label{YF-17-mesh} Flow over a $YF-17$ fighter: the local
mesh distribution.}
\centering
\includegraphics[width=0.475\textwidth,angle=-90]{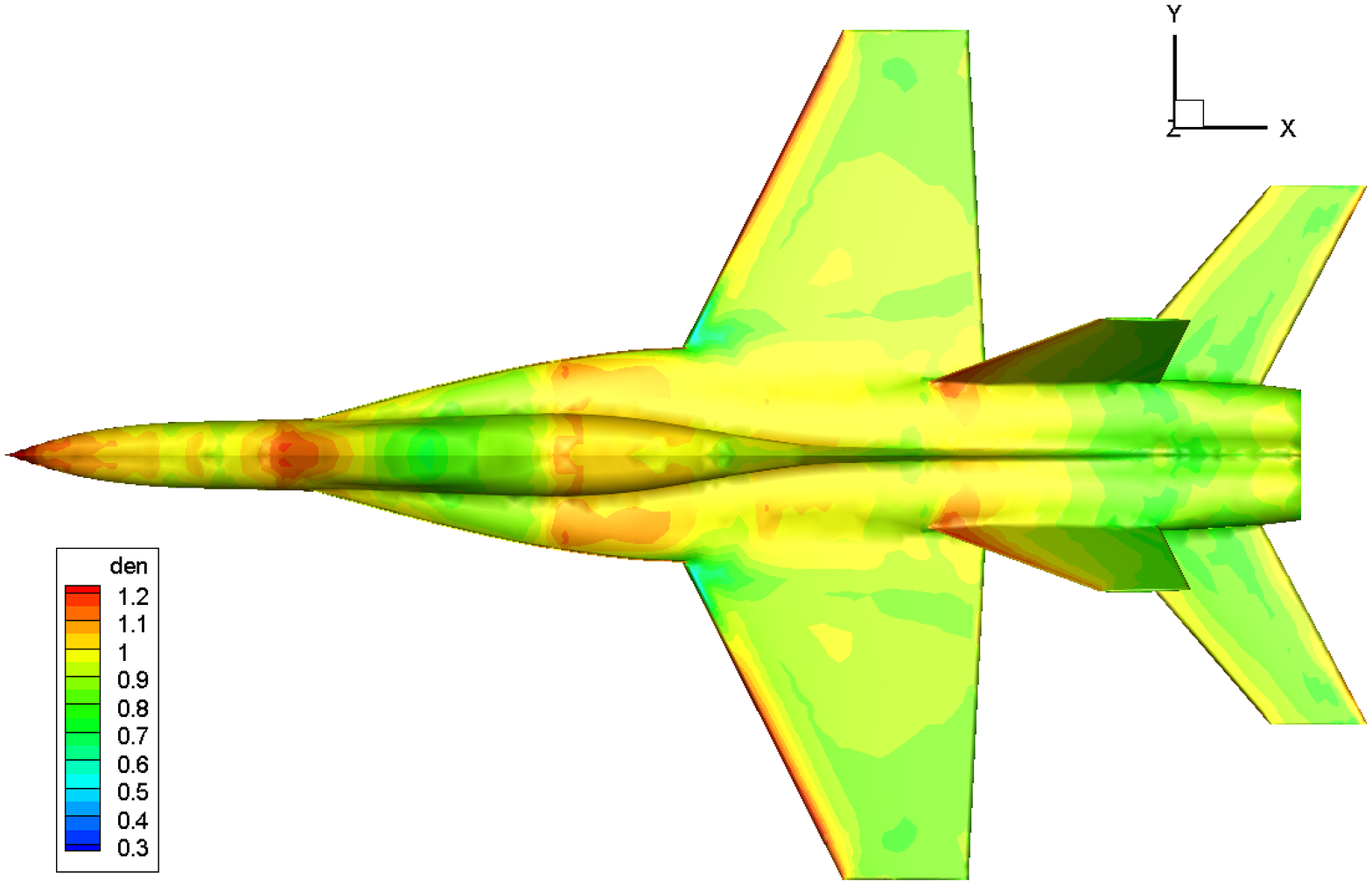}{a}
\includegraphics[width=0.475\textwidth,angle=-90]{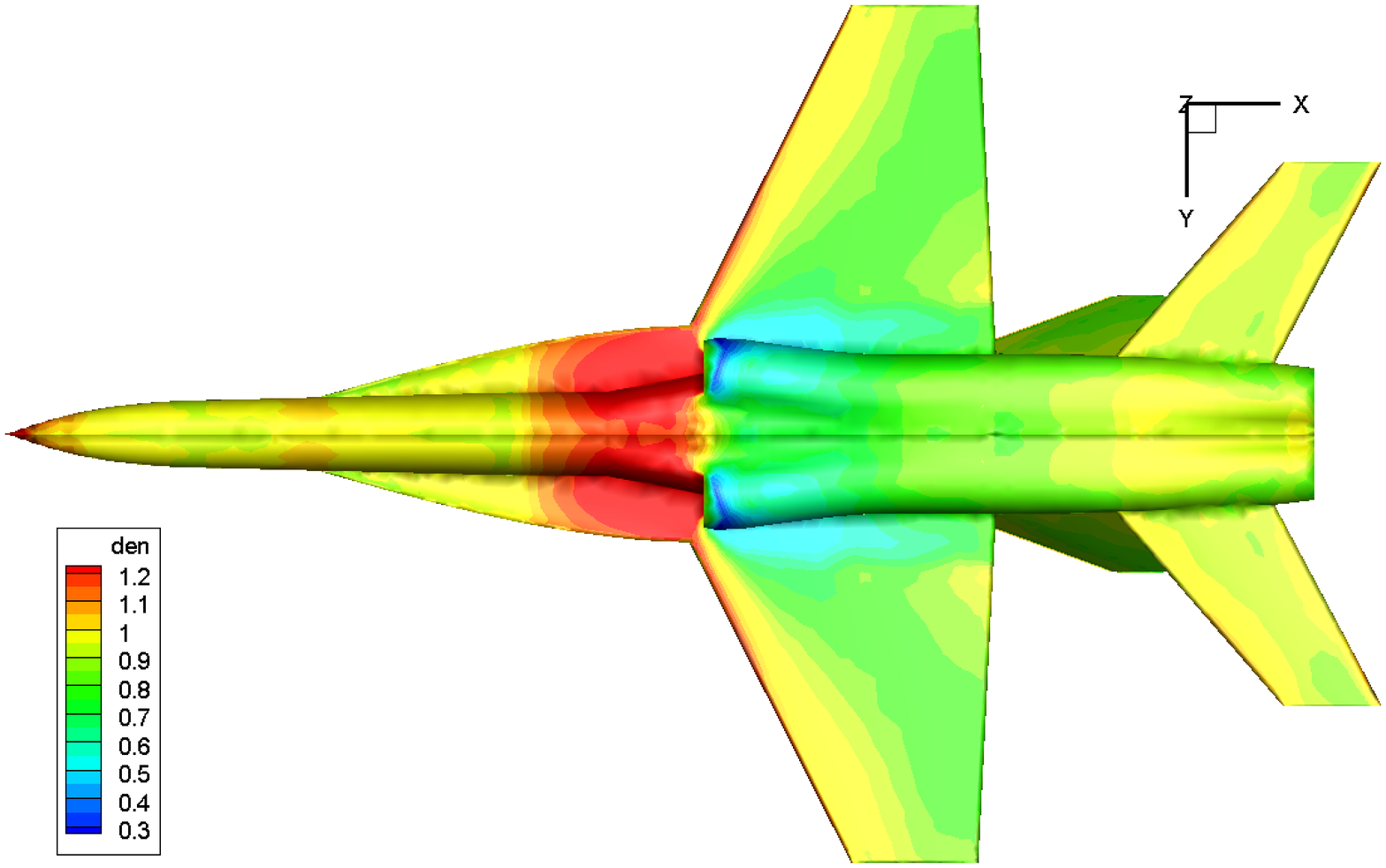} {b}
\includegraphics[width=0.575\textwidth]{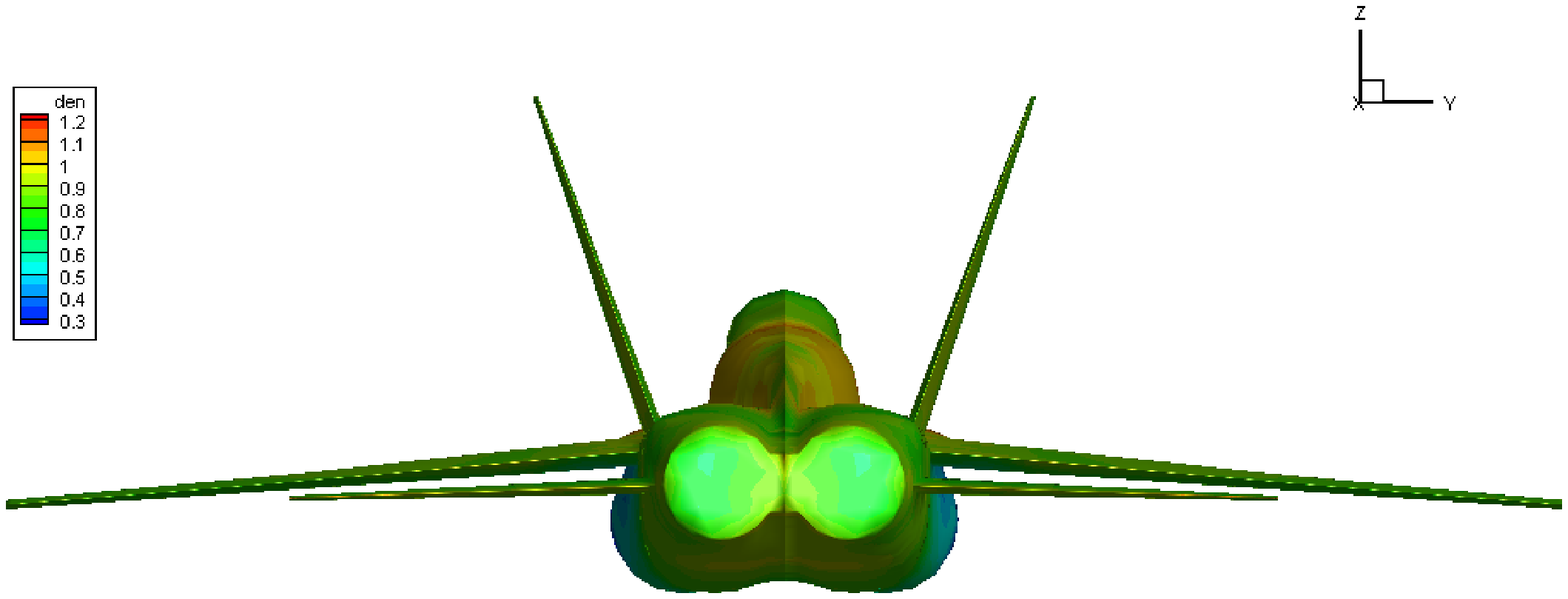} {c}
\includegraphics[width=0.575\textwidth]{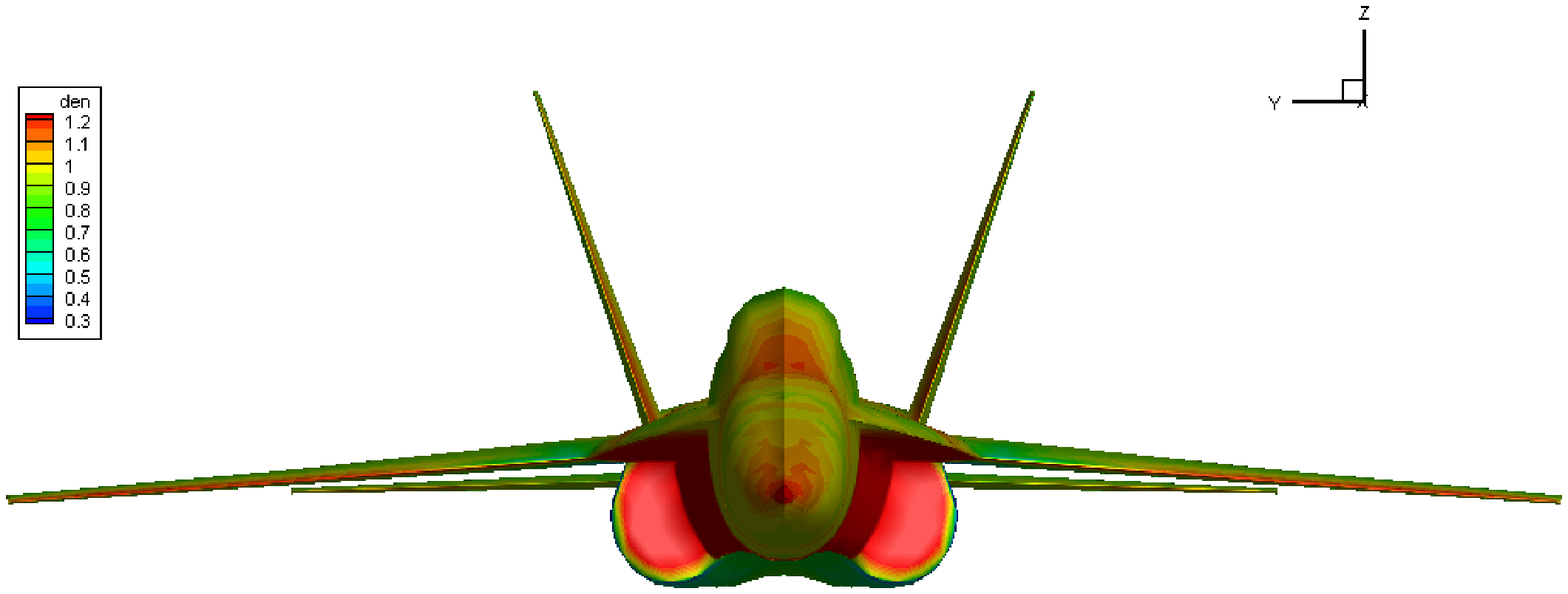} {d}
\caption{\label{YF-17-den} Flow over a $YF-17$ fighter: the steady
state density distributions from the top (a), bottom (b), forward
(c) and backward (d) views.}
\end{figure*}

\subsection{Flow over a $YF-17$ fighter}
The inviscid supersonic flow passing through a complete aircraft
model is computed. The computational mesh for a YF-17 Cobra fighter
model is shown in Fig.\ref{YF-17-mesh}, which is provided at
"https:// cgns.github. io/CGNSFiles. html". The mesh includes 174127
tetrahedral cells. The incoming Mach number and angle of attack are
given by
\begin{align*}
Ma_{\infty} =1.05, ~AoA =0^{\circ}.
\end{align*}
The surface density distributions from the top, bottom, forward and
backward views  are given in Fig.\ref{YF-17-den}. Complicated shocks
appear in the locations including the nose, cockpit-canopy wing,
horizontal stabilizer, and vertical stabilizer. This cases validate
the capability of current scheme to handle complicated geometry,
such as the mesh skewness near the wing tips and the lack of
neighboring cell for the cell near boundary corners.

 \begin{table}[!h]
    \begin{center}
    \def\temptablewidth{0.8\textwidth}{\rule{\temptablewidth}{1.0pt}}
    \begin{tabular*}{\temptablewidth}{@{\extracolsep{\fill}}c|c|c|c|c}
    Case    &       cell & face   &   Jacobi+CPU   &  LUSGS+CPU    \\
    \hline
    Cavity  &       102400 &  274560 &     568 s        &     703 s         \\
    \hline
    Sphere  &     462673 &  1174016  &      6590 s         &     8210 s    
    \end{tabular*}
    {\rule{\temptablewidth}{1.0pt}}
    \end{center}
    \caption{\label{GPU-CPU-A} Efficiency comparison: the computational time per 1000 steps for Jacobian and LUSGS iterations with CPU code.}
    \begin{center}
    \def\temptablewidth{0.7\textwidth}{\rule{\temptablewidth}{1.0pt}}
    \begin{tabular*}{\temptablewidth}{@{\extracolsep{\fill}}c|c|c|c}
    Case    &  Jacobi+GPU  &  Jacobi+CPU  &  Speedup   \\
    \hline
    Cavity   &      97 s      &      568 s     &   5.86       \\
    \hline
    Sphere   &     1093 s     &     6590 s     &   6.03   
    \end{tabular*}
    {\rule{\temptablewidth}{1.0pt}}
    \end{center}
    \caption{\label{GPU-CPU-B} Efficiency comparison: the computational time per 1000 steps and speedup for GPU and CPU code with Jacobian iteration.}
\end{table}

\begin{figure*}[!h]
\centering
\includegraphics[width=0.6\textwidth]{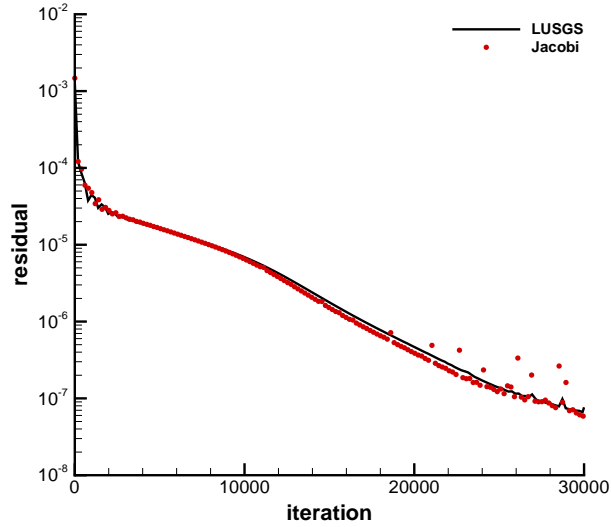}
\caption{\label{LU-GPU-res} Efficiency comparison: the comparison of convergence history with Jacobian and LUSGS iterations for the cavity case.}
\end{figure*}

\subsection{Efficiency comparison of CPU and GPU}
The efficiency comparison of CPU and GPU codes is provided with the current scheme. 
The CPU code is run with Intel i7-11700 CPU using Intel Fortran compiler with OpenMP directives, 
while Nvidia TITAN RTX is used for GPU code with Nvidia CUDA and NVFORTRAN compiler. 
The clock rates of GPU and CPU are 1.77 GHz and 2.50 GHz respectively, and the double precision is used in computation. 
The lid-driven cavity flow with $Re=1000$ and the hyperbolic inviscid flow passing a sphere with $Ma_{\infty}=2$ 
on three dimensional hybrid unstructured meshes are used to test the efficiency. The HGKS with Jacobi iteration is implemented with GPU and CPU with $k_{\max} = 3$. 
The HGKS with LUSGS iteration is only implemented with CPU. In the computation, the update part 
is implemented sequentially and the left parts are implemented in parallel. 
For these two cases, the computational time of Jacobi and LUSGS  implicit iterations  with CPU code are given in Table.\ref{GPU-CPU-A}. 
Due to the smooth GKS flux Eq.\eqref{flux-2} used in cavity case and the full GKS flux Eq.\eqref{flux-1} used in sphere case, and more control volumes, 
the sphere case costs more computational times. Meanwhile, because of the sequential update part, the LUSGS method costs more computational times. 
 The computational times of GPU code and speedups for the two cases are given in Table.\ref{GPU-CPU-B}, where speedup is defined as
\begin{equation*}
\displaystyle \text{speedup}=\frac{\text{computational time of CPU}}{\text{computational time of GPU}}.
\end{equation*}
As shown in Table.\ref{GPU-CPU-B}, 5x-6x speedup is achieved for Jacobian GPU code compared with Jacobi CPU code. 
As expected, the acceleration effect is more obvious as the number of cells increases. 
The comparison of the convergence histories of LUSGS and Jacobi implicit iterations for the cavity case is given in Fig.\ref{LU-GPU-res}. 
For both the Jacobian and LUSGS iterations, it can be observed that the residuals converge to the order of $10^{-8}$.
However,  the convergence rate of Jacobian iteration is not as stable as that of LUSGS iteration. 
In the future, the more efficient parallel implicit algorithm will be considered, and more challenging problems for compressible flows will be investigated.

\section{Conclusion}
In the paper, a third-order gas-kinetic scheme is developed on
three-dimensional arbitrary unstructured meshes for compressible
inviscid and viscous flows. In the classical WENO schemes, the
high-order of accuracy is obtained by the non-linear combination of
lower order polynomials from the candidate stencils. To achieve the
spatial accuracy, a weighted essentially non-oscillatory (WENO)
reconstruction is developed on the three-dimensional arbitrary
unstructured meshes. However, Due to the topology complicated of the
arbitrary unstructured meshes with tetrahedral, pyramidal, prismatic
and hexahedral cells, great difficulties are introduced for the
classical WENO schemes. A simple strategy of selecting stencils for
reconstruction is adopted and the topology independent linear
weights are used for tetrahedron, pyramid, prism and hexahedron.
Incorporate with the two-stage fourth-order temporal discretization,
the explicit high-order gas-kinetic schemes are developed for
unsteady problems. With lower-upper symmetric Gauss-Seidel (LU-SGS)
method, the implicit high-order gas-kinetic schemes are developed
for steady problems. Various three-dimensional numerical
experiments, including the inviscid flows and laminar flows, are
presented. The results validate the accuracy and robustness of the
proposed scheme. In the future, the high-order gas-kinetic scheme on
arbitrary unstructured meshes will be applied for the engineering
turbulent flows with high-Reynolds numbers.

\section*{Acknowledgements}
The current research of L. Pan is supported by National Natural
Science Foundation of China (11701038) and the Fundamental Research
Funds for the Central Universities. The work of K. Xu is supported
by National Natural Science Foundation of China (11772281,
91852114), Hong Kong research grant council (16206617), and the
National Numerical Windtunnel project.


\begin{thebibliography}{10}
\bibitem{ENO-un} R. Abgrall, On essentially non-oscillatory schemes on unstructured meshes: analysis and implementation, J. Comput. Phys. 114 (1994) 45-58.

\bibitem{Case-Albensoeder} S. Albensoeder, H.C. Kuhlmann,  Accurate three-dimensional lid-driven cavity flow, J. Comput. Phys.  206 (2005) 536-558.

\bibitem{WENO-ao-1} D.S. Balsara, S. Garain, C.-W. Shu. An efficient class of WENO schemes with adaptive order, J. Comput. Phys.  326 (2016) 780-804.

\bibitem{WENO-ao-2} D.S. Balsara, S. Garain, V. Florinski,  W. Boscheri,  An efficient class of WENO schemes with adaptive order for unstructured meshes, J. Comput. Phys. 404 (2020) 109062.

\bibitem{BGK-1} P.L. Bhatnagar, E.P. Gross, M. Krook, A Model for Collision Processes in Gases I: Small Amplitude Processes in Charged and Neutral One-Component Systems, Phys. Rev. 94 (1954) 511-525.

\bibitem{WENO-Z} R. Borges, M. Carmona, B. Costa, W. S. Don, An improved weighted essentially non-oscillatory scheme for hyperbolic conservation laws, J. Comput. Phys. 227 (2008) 3191-3211.

\bibitem{GKS-high-1} G.Y. Cao, L. Pan, K. Xu, High-order gas-kinetic scheme with parallel computation for direct numerical simulation of turbulent flows, J. Comput. Phys.  448 (2022) 110739.

\bibitem{GKS-rans-1} G.Y. Cao, H.M. Su, J.X Xu, K. Xu, Implicit high-order gas kinetic scheme for turbulence simulation, Aerospace Science and Technology 92 (2019) 958-971.

\bibitem{BGK-2} S. Chapman, T.G. Cowling, The Mathematical theory of non-uniform gases, third edition, Cambridge University Press, (1990).

\bibitem{DDG} J. Cheng, X.Q. Yang, X.D Liu, T.D. Liu, H. Luo, A direct discontinuous Galerkin method for the compressible Navier-Stokes equations on arbitrary grids, J. Comput. Phys.  327 (2016) 484-502.

\bibitem{DG-A}   B. Cockburn, C.-W. Shu, The Runge-Kutta discontinuous Galerkin method for conservation laws V: multidimensional systems, Journal of Computational Physics 141 (1998) 199-224.

\bibitem{GRP-high-2} Z.F. Du, J.Q. Li, A Hermite WENO reconstruction for fourth order temporal accurate schemes based on the GRP solver for hyperbolic conservation laws, J. Comput. Phys. 355 (2018) 385-396.

\bibitem{DG-B}  M. Dumbser, D. Balsara, E. Toro, C.D. Munz, A unified framework for the construction of one-step finite volume and discontinuous Galerkin schemes on unstructured meshes, J. Comput. Phys. 227 (2008) 8209-8253.

\bibitem{WENO-un4}  M. Dumbser, M. K\"{a}ser, Arbitrary high order non-oscillatory finite volume schemes on unstructured meshes for linear hyperbolic systems, J. Comput. Phys.  221 (2007) 693-723.

\bibitem{WENO-un1} C. Hu, C. W. Shu, Weighted essentially non-oscillatory schemes on triangular meshes, J. Comput. Phys. 150 (1999) 97-127.

\bibitem{LU-SGS-1} A. Jameson, S. Yoon, Lower-upper implicit schemes with multiple grids for the Euler equations, AIAA J. 25 (7) (1987) 929-935.

\bibitem{GKS-high-6} X. Ji, F. Zhao, W. Shyy, K. Xu, A HWENO reconstruction based high-order compact gas-kinetic scheme on unstructured mesh, J. Comput. Phys. 410 (2020) 109367.

\bibitem{WENO-JS} G.S. Jiang, C.-W. Shu, Efficient implementation of weighted ENO schemes, J. Comput. Phys. 126 (1996) 202-228.

\bibitem{GKS-rans-2} J. Jiang, Y.H. Qian, Implicit gas-kinetic BGK scheme with multigrid for 3d stationary transonic high-Reynolds number flows, Comput. Fluids 66 (2012) 21-28.

\bibitem{CWENO3}  O. Kolb, On the full and global accuracy of a compact third order WENO scheme, SIAM J. Numer. Anal. 52 (2014) 2335-2355.

\bibitem{CWENO1}  D. Levy, G. Puppo, G. Russo, Central WENO schemes for hyperbolic systems of conservation laws, Math. Model. Numer. Anal. 33 (1999) 547-571.

\bibitem{CWENO2}  D. Levy, G. Puppo, G. Russo, Compact central WENO schemes for multidimensional conservation laws, SIAM J. Sci. Comput. 22 (2000) 656-672.

\bibitem{GRP-high-1} J.Q. Li, Z.F. Du, A two-stage fourth order time-accurate discretization for Lax-Wendroff type flow solvers I. hyperbolic conservation laws, SIAM J. Sci. Computing, 38 (2016) 3046-3069.

\bibitem{WENO-Liu} X.D. Liu, S. Osher, T. Chan, Weighted essentially non-oscillatory schemes, J. Comput. Phys. 115 (1994) 200-212.

\bibitem{Case-Nagata}  T. Nagata, T. Nonomura, S. Takahashi, Y. Mizuno,  K. Fukuda. Investigation on subsonic to supersonic flow around a sphere at low Reynolds number of
between 50 and 300 by direct numerical simulation.  Physics of Fluids, 28 (2016) 056101.

\bibitem{GKS-high-1} L. Pan, K. Xu, Q.B. Li, J.Q. Li, An efficient and accurate two-stage fourth-order gas-kinetic scheme for the Navier-Stokes equations, J. Comput. Phys. 326 (2016) 197-221.

\bibitem{GKS-high-2} L. Pan, J. Li, K. Xu, A few benchmark test cases for higher-order Euler solvers, Numer. Math., Theory Methods Appl. 10 (4) (2017) 711-736.

\bibitem{GKS-high-5} L. Pan, G.Y. Cao, K. Xu, Fourth-order gas-kinetic scheme for turbulence simulation with multi-dimensional WENO reconstruction, Computers and Fluids 221 (2021) 104927.

\bibitem{HWENO} J.X. Qiu, C.-W. Shu, Hermite WENO schemes and their application as limiters for Runge-Kutta discontinuous Galerkin method, III: unstructured meshes, J. Sci. Comput. 39 (2009) 293-321.

\bibitem{Case-Schmitt} V. Schmitt, F. Charpin,  Pressure distributions on the ONERA-M6-wing at transonic Mach numbers, Experimental Data Base for Computer Program Assessment,
Report of the Fluid Dynamics Panel Working Group 04, AGARD AR 138, 1979.

\bibitem{WENO-un2} J. Shi, C. Hu, C.W. Shu, A Technique of Treating Negative Weights in WENO Schemes, J. Comput. Phys. 175 (2002) 108-127.

\bibitem{Case-Shu}  C. Shu, L. Wang, Y. T. Chew, Numerical Computation of Three-dimensional Incompressible Navier.Stokes Equations in Primitive Variable form by DQ Method, International Journal for Numerical Methods in Fluids 43 (2003) 345-368.

\bibitem{GKS-implicit-2} S. Tan, Q.B. Li. Time-implicit gas-kinetic scheme. Computers $\&$ Fluids, 144 (2017) 44-59.

\bibitem{WENO-un5} V. Titarev, P. Tsoutsanis, D. Drikakis, WENO schemes for mixed-element unstructured meshes, Commun. Comput. Phys. 8 (2010) 585-609.

\bibitem{WENO-un6} P. Tsoutsanis, V. Titarev, D. Drikakis, WENO schemes on arbitrary mixed-element unstructured meshes in three space dimensions, J. Comput. Phys. 230 (2011) 1585-1601.

\bibitem{WENO-un7} P. Tsoutsanis, A.F. Antoniadis, D. Drikakis,  WENO schemes on arbitrary unstructured meshes for laminar, transitional and turbulent flows,  J. Comput. Phys. 256 (2014) 254-276.

\bibitem{SV-A} Z.J. Wang, Spectral (finite) volume method for conservation laws on unstructured grids: basic formulation, J. Comput. Phys. 178 (2002) 210-251.

\bibitem{CPR-A} Z.J. Wang,  H. Gao,  A unifying lifting collocation penalty formulation including the discontinuous Galerkin, spectral volume/difference methods for conservation laws on mixed grids. J. Comput. Phys. 228  (2009) 8161-8186.


\bibitem{GKS-GPU} Y.H. Wang, L. Pan, Three-dimensional discontinuous Galerkin based high-order gas-kinetic scheme and GPU implementation,  arXiv: 2202.13821.


\bibitem{GKS-Xu1} K. Xu, Direct modeling for computational fluid dynamics: construction and application of unfied gas kinetic schemes, World Scientific (2015).

\bibitem{GKS-Xu2} K. Xu, A gas-kinetic BGK scheme for the Navier-Stokes equations and its connection with artificial dissipation and Godunov method, J. Comput. Phys. 171 (2001) 289-335.

\bibitem{DG-C}  Z. Xu, Y. Liu, C.-W. Shu, Hierarchical reconstruction for discontinuous Galerkin methods on unstructured grids with a WENO-type linear reconstruction and partial neighboring cells, J. Comput. Phys. 228 (2009) 2194-2212.

\bibitem{GKS-high-4} Y.Q. Yang, L. Pan,  K. Xu,  High-order gas-kinetic scheme on three-dimensional unstructured meshes for compressible flows, Physics of Fluids 33 (2021) 096102.

\bibitem{LU-SGS-2}  S. Yoon, A. Jameson, Lower-upper symmetric-Gauss-Seidel method for the Euler and Navier-Stokes equations, AIAA J. 26 (9) (1988) 1025-1026.

\bibitem{WENO-un3} F.X. Zhao, L. Pan, S.H. Wang, Weighted essentially non-oscillatory scheme on unstructured quadrilateral and triangular meshes for hyperbolic conservation laws, J. Comput. Phys. 374 (2018) 605-624.

\bibitem{GKS-high-7} F.X. Zhao, X. Ji , W. Shyy, K. Xu,  A compact high-order gas-kinetic scheme on unstructured mesh for acoustic and shock wave computations, J. Comput. Phys.  449 (2022) 110812

\bibitem{WENO-simple-1} J. Zhu, J.X. Qiu, A new fifth order finite difference weno scheme for solving hyperbolic conservation laws. J. Comput. Phys. 318 (2016) 110-121.

\bibitem{WENO-simple-2} J. Zhu, J.X. Qiu, New finite volume weighted essentially non-oscillatory scheme on triangular meshes, SIAM J. Sci. Computing, 40 (2018) 903-928.

\bibitem{GKS-implicit-1} Y.J. Zhu, C.W. Zhong, K. Xu.Implicit unified gas-kinetic scheme for steady state solutions in all flow regimes. Journal of Computational Physics, 315 (2016) 16-38.
\end{thebibliography}
\end{document}